\documentstyle[11pt,twoside]{article}

\newcommand{\ncm}{\newcommand}
\ncm{\rncm}{\renewcommand}

\ncm{\lb}{\label}

\newcommand{\be}{\begin{equation}}
\newcommand{\ee}{\end{equation}}
\newcommand{\bea}{\begin{eqnarray}}
\newcommand{\ea}{\end{eqnarray}}

\newcommand{\om}{\omega}
\newcommand{\Om}{\Omega}

\newtheorem{theorem}{Theorem }[subsection]
\newtheorem{lemma}[theorem]{Lemma }
\newtheorem{proposition}[theorem]{Proposition }
\newtheorem{corollary}[theorem]{Corollary }
\newtheorem{definition}[theorem]{Definition }

\ncm{\Theorem}[2]{\begin{theorem} \lb{Thm #1} {\sl #2} \end{theorem}}

\ncm{\thm}[1]{Theorem \ref{Thm #1}}

\ncm{\Definition}[2]{\begin{definition} \lb{Def #1}{\rm #2}\end{definition}}

\ncm{\defi}[1]{Definition \ref{Def #1}}

\ncm{\Lemma}[2]{\begin{lemma} \lb{Lem #1} {\sl #2} \end{lemma}}

\ncm{\lem}[1]{Lemma \ref{Lem #1}}

\ncm{\Proposition}[2]{\begin{proposition}\lb{Prop #1}{\sl #2}
                      \end{proposition}}

\ncm{\prop}[1]{Proposition \ref{Prop #1}}

\ncm{\Corollary}[2]{\begin{corollary}\lb{Cor #1} {\sl #2} \end{corollary}}

\ncm{\cor}[1]{Corollary \ref{Cor #1}}

\ncm{\Thm}{\Theorem}
\ncm{\Def}{\Definition}
\ncm{\Lem}{\Lemma}
\ncm{\Prop}{\Proposition}
\ncm{\Cor}{\Corollary}

\renewcommand{\theequation}{\mbox{\arabic{section}.\arabic{equation}}}

\ncm{\setc}[1]{\setcounter{equation}{#1}}
\rncm{\sec}{\setc{0}\section}
\ncm{\no}[1]{(\ref{#1})}
\ncm{\Eq}[1]{Eq.(\ref{#1})}

\ncm{\beq}{\begin{equation}}
\ncm{\bleq}[1]{\beq\lb{#1}}

\ncm{\eeq}{\end{equation}}

\ncm{\beanon}{\begin{eqnarray*}}

\ncm{\eea}{\end{eqnarray}}

\ncm{\eeanon}{\end{eqnarray*}}

\ncm{\qed}{\hspace*{0.4cm}\rule{0.24cm}{0.24cm}}

\ncm{\A}{{\cal A}}
\ncm{\hA}{\hat{\cal A}}
\ncm{\B}{{\cal B}}

\ncm{\G}{{\cal G}}
\ncm{\J}{{\cal J}}
\ncm{\I}{{\cal I}}
\rncm{\L}{{\cal L}}
\rncm{\H}{{\cal H}}
\ncm{\C}{{\cal C}}
\ncm{\D}{{\cal D}}
\ncm{\F}{{\cal F}}
\ncm{\K}{{\cal K}}
\ncm{\T}{{\cal T}}
\ncm{\N}{{\cal N}}
\ncm{\M}{{\cal M}}
\ncm{\R}{{\cal R}}
\ncm{\Z}{{\cal Z}}
\ncm{\CC}{{\bf C}}
\ncm{\ZZ}{{\bf Z}}
\ncm{\RR}{{\bf R}}

\ncm{\one}{{\bf 1}}

\ncm{\NN}{{\bf N}}
\ncm{\eee}{{\bf e}}
\def\lcros{\raise1.5pt\hbox{$\scriptstyle\triangleright$}\!
           \raise1.9pt\hbox{$\scriptscriptstyle < \,$}}
\def\<cros{\lcros}
\def\Lcros{\Lambda_{cros}}
\def\cros{\raise1.9pt\hbox{$\scriptscriptstyle  > $}\!
          \raise1.5pt\hbox{$\scriptstyle\triangleleft\,$}}
\def\>cros{\cros}
\def\crosAd{\cros\!_{Ad}\,}

\def\la{\rightharpoonup}
\def\ra{\leftharpoonup}
\def\arr{\la}
\def\arl{\ra}

\def\lef{{\,\hbox{$\textstyle\triangleright$}\,}}
\def\re{\lef}

\def\reli{\raise1.5pt\hbox{$\scriptstyle\triangleright$}\!
           \raise1.5pt\hbox{$\scriptstyle\triangleleft\,$}}

\def\span{\mbox{span}\,}
\def\Rep{\mbox{Rep}\,}
\def\Ker{\mbox{Ker}\,}
\def\End{\mbox{End}\,}
\def\Aut{\mbox{Aut}\,}
\def\Out{\mbox{Out}\,}
\def\pr{\mbox{pr}\,}
\def\prr{\mbox{pr}_{\Out\M}}
\def\Ad{\mbox{Ad}\,}
\def\Ind{\mbox{Ind}\,}
\def\Index{\mbox{Ind}\,}

\def\Hom{\mbox{Hom}\,}
\def\id{\mbox{id}\,}
\def\bra{\langle}
\def\ket{\rangle}
\ncm{\bsn}{\bigskip\noindent}
\ncm{\proof}{\bsn{\bf Proof:\ }}
\def\o{\otimes}
\def\e{\varepsilon}
\def\x{\times}
\def\al{\alpha}
\def\be{\beta}
\def\r{\rho}
\def\s{\sigma}
\def\Del{\Delta}
\ncm{\1}{_{(1)}}
\ncm{\2}{_{(2)}}
\ncm{\3}{_{(3)}}
\ncm{\4}{_{(4)}}
\ncm{\ba}{\begin{array}}
\rncm{\ea}{\end{array}}
\rncm{\l}{\lambda}

\begin{document}

\title{ \bf Weak Hopf Algebras and \\
Reducible Jones Inclusions of Depth 2. \\
I: From Crossed Products to Jones Towers}

\date{April 3, 1998}

\author{\sc F. Nill$^*$, K. Szlach\'anyi$^\dagger$, H.-W. Wiesbrock$^*$}

\maketitle

{\rncm{\thefootnote}{\fnsymbol{footnote}}
\footnotetext[1]
{Freie Universit\"at Berlin, Institut f\"ur theoretische Physik,
Arnimallee 14, D-14195 Berlin, Germany.\\
email: nill@physik.fu-berlin.de\quad wiesbroc@physik.fu-berlin.de\\
Supported by DFG, SFB 288 {\em Differentialgeometrie und
  Quantenphysik}}
\footnotetext[2]
{Central Research Institute for Physics,
H-1525  Budapest 114, P.O.B. 49, Hungary.\\
email: szlach@rmki.kfki.hu}
}

\vspace{-7cm}
\rightline{\normalsize q-alg/9806130}
\vspace{7cm}

\begin{abstract}
We apply the theory of finite dimensional weak $C^*$-Hopf algebras
$\A$ as developed by G. B\"ohm,  F. Nill and  K. Szlach\'anyi
[BSz,Sz,N1,N2,BNS] to
study reducible inclusion triples of von-Neumann algebras
$\N\subset\M\subset\M\cros\A$, where $\M$ is an $\A$-module algebra with left $\A$-action $\lef:\A\o\M\to\M$, $\N\equiv\M^\A$ is the fixed
point algebra and $\M\cros\A$ is the crossed
product extension. Here ``weak'' means that
the coproduct $\Delta$ on $\A$ is non-unital, requiring various
modifications of the standard definitions for Hopf 
(co-)actions and crossed products. We show that
normalized positive and nondegenerate left integrals $l\in\A$ give
rise to faithful conditional expectations $E_l:\M\to\N$ via
$E_l(m):=l\lef m$, where under certain regularity conditions this
correspondence is one-to-one. Associated
with such left integrals 
we construct ``Jones projections'' $e_l\in\A$  obeying for all
$m\in\M$ the Jones relations $e_lme_l=E_l(m)e_l=e_lE_l(m)$ as an
identity in $\M\cros\A$. We also
present a concept of Plancherel-duality ("p-duality") for positive nondegenerate left
integrals, where by definition the p-dual
$\l\in\hat\A$ of $l\in\A$ is the
unique solution of $\l\arr e_l\equiv\hat E_\l(e_l)=\one$. Here
$\arr$ denotes the canonical left action of the dual weak Hopf
algebra $\hat\A$ on $\A$.
Finally, we prove that $\N\subset\M$ always has finite index and depth
2 and that the basic
Jones construction for $\N\subset\M$ is given  by the ideal
$\M_1:=\M e_l\M\subset\M\cros\A$.
In this way  p-dual left integrals precisely correspond to
Haagerup-dual conditional expectations. We
give appropriate regularity conditions (such as outerness of the
$\A$-action) guaranteeing $\M_1=\M\cros\A$. Under these
conditions the standard invariant for $\N\subset\M$ is
given by the Jones triple $\A_L\subset\A\subset\A\cros\hat\A$, where
$\A_L:=\one_\A\arl\hA
\subset\A$
is a nontrivial subalgebra isomorphic to $\A\lef\one_\M\subset\M$.
As a particular example we discuss crossed products by
partly inner group actions.

In a subsequent paper we will show that converseley any reducible
 finite index and depth-2 Jones tower of von-Neumann factors arises in
 this way, where the inclusions are irreducible if and only if $\A$
 and $\hat\A$ are ordinary Kac algebras. This generalizes the well
 known results on irreducible depth-2 inclusions obtained by
various people after a proposal by A. Ocneanu.

\end{abstract}


\newpage
\footnotesize
\tableofcontents

\normalsize


\sec{Introduction}
\subsection{Motivation}

There is a lot of work concerning the classification of irreducible Jones
inclusions of factors, see [Po] for example. In the special case of inclusions
$\N\subset \M$ of depth 2 in the sense of Ocneanu [Oc1] it is
known that the subfactor $\N \equiv \M^\A$ always arises as the fixed point
algebra under an outer action of a Kac algebra $\A$ on $\M$, see
[Da, Lo1, Szy] for the finite index case and [EN, NW] for the infinite
index case. Converseley, any such outer action gives rise to a Jones triple
$\N\equiv \M^\A\subset\M \subset\M \cros\A$ obeying Ocneanu's depth 2
condition and  the ``symmetry algebra" can be recovered
as $\A=\N' \cap (\M\cros\A)$.

With this work we initiate a project where we want to
generalize these ideas to reducible finite index inclusions of
von-Neumann algbras.
As a motivation we point out that there are many natural situations,
where inclusions are not irreducible. Let us give two examples:

\bsn
{\em 1. A partly inner group action:}

\smallskip
Let
$
\alpha :G\to \Aut\M
$
be an action of a finite group $G$ on a von-Neumann factor $\M$
and let $\N:=\M^G\subset\M$ be the fixed point algebra and
$\M\subset\M \cros G$ the crossed product extension.
%
Denote
$\pr_{\Out\M}:\Aut\M\to\Out\M$ the canonical projection onto
the group of outer automorphisms of $\M$.
Then $H:= \Ker (\pi \circ \alpha)$ is a normal
subgroup of $G$ and all implementers $u(h)\in\M,\ h\in H,$
(i.e. satisfying $\alpha_h = \Ad u(h)$) commute with $\N$,
$$
u(h)\in\N' \cap \M\quad\forall h\in H
$$
Hence, if $G$ does not act purely outerly, we have a non-trivial relative
commutant. By definition this means that the inclusion $\N\subset\M$
is reducible. Similarly,
$\M\subset\M \cros G$ is irreducible if and only if $G$
acts outerly on $\M$.
We will frequently come back to this example within our approach,
showing that for a non-outer Galois $G$-action on $\M$
\footnote{The notion of a Galois action has been introduced in
[CS]. Equivalently, this means that
$\N\subset\M\subset\M\cros G$ is a Jones triple, i.e. $\M\cros G=\M h\M$, where $h\in\CC G$ is the normalised Haar integral, see also Appendix A.}
the second relative commutant
$\A:=\N'\cap\M\cros G$ naturally acquires the structure of a
{\em weak Hopf algebra}.

\bsn
{\em 2. The observable algebra in the reduced field bundle:}

\smallskip
Let $\A$ denote the algebra of quasilocal observables of a
quantum field theory
in the algebraic frame work [Haa]. Let us assume that the theory has
only finitely many irreducible sectors $\rho_i,\ i=1\dots n$, with
finite statistics.
We define the ``master" endomorphism $\rho:=
\mathop{\oplus}\limits^n_{i=1} \rho_i$.
Then the inclusion
$$
\rho(\A) \subset\A
$$
is by construction reducible. Taking the minimal conditional expectation
$E_0:\A \to \rho (\A)$ (corresponding to the standard left inverse of $\rho$)
and performing the basic Jones construction one ends up with
the reduced field bundle $\F_{red}\supset\A\supset\rho(\A)$.
As has been noticed by [NR,R] the symmetry algebra associated
with this inclusion can also be described by a weak Hopf algebra.
(Actually, the first proposal for axioms of weak Hopf algebras
by [N1] had been based on this example).

\bigskip
Both examples share the further property of having depth 2.
Motivated by the irreducible case we therefore expect a
Hopf algebra like
symmetry encoded in any reducible finite index and
depth-2 inclusion of von-Neumann algebras. Our claim is
that the appropriate notion of a symmetry which can always be recovered in
such a setting is indeed the concept of a {\em weak $C^*$-Hopf algebra}.

\bigskip
As a motivation we recall
that already the investigation of the inner symmetry
of the Ising model has led to a weaker concept then the
familiar one of a Hopf-symmetry, see [MS]. There the authors
relaxed  the usual
property that the coproduct of the identity should be the tensor product of
the identity
with itself, which is one of the axioms for Hopf algebras.
Instead they only required that the coproduct of the identity is a projection.
In the case of [MS], this projecton is the one onto the ``good
sub-representations" in a tensor product of good representations.

Now in Jones theory we always have a two step periodicity .
Therefore we expect a concept of a symmetry algebra $\A$ such that the dual
object $\hat\A$ is an algebra of the same type. In this way we can develop
the picture of a Jones tower generated by taking alternating crossed products
with $\A$ and $\hat A$, respectively, reflecting the two step periodicity.
Thus, we are naturally lead to depart from the Mack-Schomerus setting in two
essential ways. First, in order that $\hat \A$ be an associative algebra
we require the coproduct on $\A$ to be strictly coassociative (in [MS]
this had been relaxed to quasi-coassociativity  in the sense of [Dr]).
Second,
if the coproduct on $\A$ is non-unital then the counit on $\hat\A$ cannot be
an algebra morphism.
Hence we drop the requirement of counits being algebra maps
altogether (equivalently, the coproduct on $\hat\A$ is allowed to be
non-unital as well).

\bigskip
Before putting these ideas into a precise definition let
us give some more heuristic pictures of the kind structures to
be expected of a ``symmetry algebra" $\A$ appearing in a
reducible depth-2 Jones tower
$\dots\M_{i-1}\subset\M_{i}\subset\M_{i+1}\subset\dots$.

First, by experience from the irreducible depth-2 case the
algebras $\A=\M'_{i-1}\cap\M_{i+1}$ and
$\hat\A=\M_i'\cap\M_{i+2}$ should be dual to each other. Next,
as an intrinsic feature we should always have two distinguished
commuting subalgebras $\A_{L/R}\subset\A$ given by
\bea
\A_L &=& \M_{i-1}'\cap\M_i
\lb{0.1}\\
\A_R &=& \M_{i}'\cap\M_{i+1}
\lb{0.2}
\eea
and similarly for $\hat\A$.
Also, proceeding alternatingly up the tower we would expect
intrinsic identifications $\hat\A_L\cong\A_R$ and
$\hat\A_R\cong\A_L$, whereas $\A_L$ and $\A_R$ should naturally
be anti-isomorphic. If we are dealing with factors $\M_i$, then
we should have
$\A_L\cap\A_R=\CC=\hat\A_L\cap\hat\A_R$
and otherwise we would expect
\bea
\A_L\cap\A_R &=& C(\M_i)
\lb{0.3}\\
\hat\A_L\cap\hat\A_R &=& C(\M_{i+1})
\lb{0.4}
\eea
where $C(\M_i)$ denotes the center of $\M_i$.
Next, observe that the depth-2 condition gives
$\M_{i+1}=\M_i\vee\A$ and therefore there should exist a natural isomorphism
\beq\lb{0.5}
\hat\A_L\cap\hat\A_R \cong C(\M_{i+1})\cong\A_R\cap C(\A)\ .
\eeq
Applying the modular conjugation of $\M_i$ (which intertwines
$\A_L$ and $\A_R$) we further conclude
\beq\lb{0.5a}
\C(M_{i+1})\cong\C(\M_{i-1})=\A_L\cap\C(\A)
\eeq
Finally, we would expect
\beq\lb{0.6}
C(\M_i)\cap C(\M_{i+1})\cong\A_L\cap\A_R\cap C(\A)\cong
\hat\A_L\cap\hat\A_R\cap C(\hat\A)
\eeq
since in any Jones tower this abelian algebra is globally
fixed, i.e. we have for all $n\in\NN$
\beq\lb{0.7}
C(\M_i)\cap C(\M_{i+1}) = C(\M_{i+n})\cap C(\M_{i+n+1})
\eeq
as a strict identity .

It turns out that the existence of the subalgebras
$\A_{L/R}\subset\A$ and $\hat\A_{L/R}\subset\hat\A$ satisfying
all these properties will indeed follow from our axioms for
(finite dimensional) weak Hopf algebras $\A$ and $\hat\A$ below,
provided $\A$ acts regularly (in a suitable sense) on $\M_i$
such that $\M_{i-1}$ is the fixed point algebra under this action
and $\M_{i+1}=\M_i\cros\A$ is the crossed product equipped with
its canonical dual $\hat\A$-action.
Note that in the particular example of a partly inner group action on a
factor $\M$ one has
\bea
\lb{0.7a}
\A&:=&\N'\cap(\M\cros G)=\span\{u(h)g\mid h\in H,\ g\in G\}
\\\lb{0.7b}
\A_L&:=&\N'\cap\M = \span\{u(h)\mid h\in H\}
\\\lb{0.7c}
\A_R&:=&\M'\cap(\M\cros G)=\span\{u(h)h^{-1}\mid h\in H\}
\eea
The weak Hopf algebra structure associated with this example
will be described in Section 2.5 and, more generally, in
Appendix B.

\subsection{The basic setting}

We now give our axioms for weak $C^*$-Hopf algebras, where we restrict
ourselves
to the finite dimensional case.
A first proposal of such a structure has been made by [N1] and
independently by [BSz,Sz].
Meanwhile the theory is well established [N2,BNS], and in
particular we have unified our axioms (see [N2] for a detailed
discussion of the relations between various sets of axioms).
We remark, that the {\em quantum groupoids} of Ocneanu [Oc2] and the
{\em face algebras} of Hayashi [Hay] are special types of weak Hopf
algebras in our sense,
where the above mentioned subalgebras
$\A_{L/R}\subset\A$ are restricted to be abelian.
Moreover, as has been shown in [N2], weak $C^*$-Hopf algebras with an involutive antipode precisely yield the generalized Kac algebras of Yamanouchi [Ya1]. 

\Definition{1.1}{
({\em Axioms for finite dimensional weak $C^*$-Hopf algebras})
\\
A finite dimensional weak $C^*$-Hopf-Algebra $\A$ is a finite
dimensional $C^*$-algebra with linear structure maps
$$
\begin{array}{rcll}
\Delta &:& \A\to \A\otimes \A \quad & \mbox{coproduct}\\
\varepsilon &:& \A\to \CC & \mbox{counit}\\
S &:& \A\to \A & \mbox{antipode}
\end{array}
$$
satisfying the following axioms for $x,y,z\in\A:$
\beanon
I.&a)& \Delta(xy) = \Delta (x)\Delta(y)\\
&b)& \Delta(x^*) = \Delta(x)^*\\
&c)& (\Delta \otimes id) \circ \Delta = (id\otimes \Delta)\circ\Delta
\hbox to 7cm{}\\
&d)& (\Delta(\one)\otimes\one)\,(\one\otimes\Delta(\one)) =
\one_{(1)}\otimes\one_{(2)}\otimes\one_{(3)}\\
~\\
II.&a)& (\epsilon \otimes id ) \circ \Delta =id =(id \otimes
\varepsilon)\circ \Delta\\
&b)& \varepsilon (xyz)=\varepsilon(xy_{(1)})\varepsilon(y_{(2)}z)
\\
\, \\
III.
  &a)& S(x_{(1)}) x_{(2)} = \one_{(1)} \varepsilon (x\one_{(2)})\\
  &b)& x_{(1)} S(x_{(2)}) = \varepsilon (\one_{(1)}x)\one_{(2)} \\
  &c)& S(x_{(1)})x_{(2)}S(x_{(3)}) = S(x)
\eeanon
}
Here we have used the Sweedler notation $\Delta(x) =x_{(1)} \otimes
x_{(2)},\ (\Delta \otimes id) \circ \Delta(x) = (id\otimes
\Delta)\circ\Delta(x)=x_{(1)} \otimes x_{(2)}\otimes x_{(3)}$,
etc., implying as usual a summation convention in the
corresponding tensor product spaces.
Note that $IIIa)$ or $IIIb)$ also imply
\beq\lb{1.1}
 x_{(1)} S(x_{(2)})x_{(3)} = x\ .
\eeq
Moreover, due to the compatibility requirement $Ib)$,
the axioms $Id)$ and $IIb)$ are equivalent to,
respectively,
\beanon
I&d')&
(\one\otimes\Delta(\one))\,(\Delta(\one)\otimes\one) =
\one_{(1)}\otimes\one_{(2)}\otimes\one_{(3)}
\\
II&b')&
\varepsilon (xyz)=\varepsilon(xy_{(2)})\varepsilon(y_{(1)}z)
\eeanon
A weak Hopf algebra is an ordinary Hopf algebra, if
the counit $\varepsilon$ is multiplicative, or if
$\Delta(\one) = \one\o\one$, or if instead of $III$ the
antipode $S$ satisfies the usual axioms
$S(x_{(1)}) x_{(2)} =
\varepsilon(x){\bf 1}= x_{(1)} S(x_{(2)})$
(in which case axiom $IIIc)$ would be a consequence).
As for ordinary
Hopf algebras, under the
axioms $I- III$ the antipode is uniquely determined
provided it exists [N2,BNS]. Moreover, it is invertible and satisfies
\bea
 S(xy) &=& S(y) S(x)\lb{1.2}\\
 \Delta \circ S &=& (S\otimes S)\circ \Delta^{op}\lb{1.3}\\
 S(x^*)^* &=& S^{-1}(x)\lb{1.4}
\eea
We remark that it will often be easier to check (\ref{1.2}) and (\ref{1.3})
explicitely, in which case the axioms $IIIc$ and $IIb),b')$
become redundant and $Id),d')$ may be relaxed by only requiring [N2]
\beq\lb{1.11'}
[(\one\otimes\Delta(\one))\,,\,(\Delta(\one)\otimes\one)] = 0
\eeq
Note that $Id)$ and $Id')$ together always imply \no{1.11'}.

Under the present setting the following axioms of [BSz,Sz]
are a consequence [N2,BNS]
\bea
 S(x_{(1)}) x_{(2)}\otimes x_{(3)} &=& \one_{(1)} \o x\one_{(2)}\lb{1.5}\\
 x_{(1)}\o x_{(2)} S(x_{(3)}) &=& \one_{(1)}x\o \one_{(2)}\lb{1.6} \\
 S^{-1}(x_{(3)}) x_{(2)}\o x_{(1)} &=& \one_{(2)} \o x\one_{(1)}\lb{1.7}\\
 x_{(3)}\o x_{(2)} S^{-1}(x_{(1)}) &=&\one_{(2)}x\o \one_{(1)}\lb{1.8}\ .
\eea

\bigskip
Given a finite dimensional weak $C^*$-Hopf algebra $\A$ as above
we denote $\Rep\A$ the category of unital $*$-representations
$\pi_\alpha :\A\to\End V_\alpha$ on finite dimensional Hilbert
spaces $V_\alpha$. Then $\Rep\A$ becomes a rigid $C^*$-tensor category with
monoidal product
$\pi_{\al\x\be}:=(\pi_\al\o\pi_\be)\circ\Delta$ and
$V_{\al\x\be}:=\pi_{\al\x\be}(\one)\,(V_\al\o V_\be)$.
It is shown in [BNS] that the counit $\e$ provides a positive
(unnormalized) state on $\A$. Let $\pi_\varepsilon$
denote the associated GNS-representation of $\A$
with cyclic vector $\Omega_\e\in V_\e$.
Then $\pi_\e$ is the unit object in $\Rep\A$ with unitary
intertwiners $U_\al^r :V_\al\to V_{\al\x\e}$ and
$U_\al^\ell :V_\al\to V_{\e\x\al}$ given on $v\in V_\al$ by
[BSz,Sz,N2]
\bea
U_\al^r\,v &:=& \pi_{\al\x\e}(\one)\,(v\o\Omega_\e)\lb{1.9}\\
U_\al^\ell\,v &:=& \pi_{\e\x\al}(\one)\,(\Omega_\e\o v)\ .\lb{1.10}
\eea
As has been pointed out in [N2], the intertwining relations
and the triangle identities for these intertwiners are
precisely guaranteed by the axioms $IIb)$ and $IIb')$.
Following [BSz,Sz] we call $\A$ {\em pure} if $\e$ is pure and
therefore $\pi_\e$ is {\em irreducible}.

As for ordinary Hopf algebras, our antipode axioms also provide
a rigidity structure on $\Rep\A$, see [N2] for a general
discussion and [BSz,Sz] for the $C^*$-formulation.

\bigskip
Associated with any weak Hopf algebra
$(\A,\one,\varepsilon,\Delta,S)$
we can immediately define the dual object $(\hat \A,\hat\one,\hat\varepsilon,
\hat\Delta,\hat S)$ by dualising all structure maps: Let $\hat\A$ be the dual
space of $\A$ and denote
$$\bra\cdot\mid\cdot\ket :\hat\A \otimes \A \to \CC$$
the canonical pairing. Then for $\psi,\varphi\in\hat\A$ and $x\in\A$
\beq\lb{1.11}
\langle\psi\varphi\mid x\rangle :=
\langle\varphi\otimes \psi\mid\Delta(x)\rangle
\eeq
defines an associative product on $\hat A$. Similarly we get
the coproduct, the unit, the counit, the antipode and the
$*$-structure by putting for $\psi \in\hat\A$ and $x,y\in\A$
\bea
\langle \hat\Delta\psi\mid x\otimes y\rangle &:=&
\langle \psi\mid xy\rangle\lb{1.12}\\
\langle \hat\one\mid x\rangle &:=& \varepsilon(x)\lb{1.13}\\
\hat\varepsilon(\psi) &:=& \langle \psi\mid\one\rangle\lb{1.14}\\
\langle \hat S(\psi)\mid x\rangle &:=& \langle \psi\mid S(x)\rangle\lb{1.15}\\
\langle\psi^*\mid x\rangle &:=& \overline{\langle \psi\mid x_*\rangle}\lb{1.16}
\eea
where in the last line we have introduced the multiplicative
antilinear involution
\beq\lb{1.17}
 x_* := S(x)^*\ .
\eeq
It turns out that the quintuple $(\hat\A,\hat\one,\hat\varepsilon,\hat\Delta,
\hat S)$ again defines a weak $C^*$-Hopf algebra [BNS].
Note in particular, that due to \no{1.13} and \no{1.14} the
axioms $Id),d')$ and $IIb),b')$, respectively, become dual to
each other.

\bigskip
The theory of integrals generalizes from ordinary Hopf algebras as
follows. An element $l\in\A$ is called a {\em left integral}, if one -
and hence all - of the following three equivalent conditions hold
\bea
(\one\o a)\,\Delta(l)
& = & (S(a)\o\one)\,\Delta(l),\quad\forall a\in\A\lb{1.18}\\
a\,l & = & a_{(1)}S(a_{(2)})\,l\equiv \e(\one\1 a)\one\2 l,
\quad\forall a\in\A\lb{1.19}\\
a\,l & = & a_{(2)}S^{-1}(a_{(1)})\,l\equiv \e(\one\2 a)\one\1 l,
\quad\forall a\in\A\lb{1.20}\ .
\eea
The equivalence of these conditions follows from $IIa$ and
(\ref{1.5})-(\ref{1.8}) [BNS]. Similarly,  $r\in\A$ is called a
{\em right integral}, if one -
and hence all - of the following three equivalent conditions hold
\bea
\Delta(r)(a\o\one) &=& \Delta(r)(\one\o S(a)),\quad\forall a\in\A\lb{1.21}\\
r\,a &=& r\,S(a_{(1)})a_{(2)}\equiv r\one\1\e(a\one\2),
\quad\forall a\in\A\lb{1.22}\\
r\,a &=& r\,S^{-1}(a_{(2)})a_{(1)}\equiv r\one\2\e(a\one\1),
\quad\forall a\in\A\lb{1.23}\ .
\eea
Note that $l\in\A$ is a left integral if and only if $S(l)$ (or $l^*$)
is a right integral. If $l$ is a left and a right integral, then it is
called a {\em two sided integral}.
A left (right) integral is called {\em positive} and/or {\em
nondegenerate}, if it is positive and/or nondegenerate as a
linear functional on $\hat\A$.
As opposed to ordinary finite
dimensional Hopf algebras, in our case the space of left (right)
integrals is more than one dimensional. In fact,
as a linear space it will be naturally isomorphic to the subspaces
$\A_{L,R}\subset\A$ mentioned in Section 1.1, ([BNS], see also
\lem{2.22} below).
Following [BSz,Sz] we define the {\em normalized Haar integral}
$h\in\A$ to be the
unique two sided integral satisfying the normalization condition
\beq\lb{1.24}
S(h_{(1)})h_{(2)} = h_{(1)}S(h_{(2)}) = \one\ .
\eeq
The uniqueness of $h$ is obvious from the definitions which in
particular imply
\beq\lb{1.25}
h=h^2=h^*=S(h)
\eeq
The existence of $h$ is proven in [BNS], where it is also shown
that as a functional on $\hA$ it is positive and nondegenerate.

A more detailed investigation of weak $C^*$-Hopf algebras, including a
theory of Fourier transformations and the modular theory for
the normalized Haar integral, is given in [BNS].
From now on, by a weak Hopf algebra $\A$ we will always mean a finite
dimensional weak $C^*$-Hopf algebra.

\bigskip
Next, we need the notion of an $\A$-module algebra $\M$ allowing for a natural
definition of a `fixed point' subalgebra $\N\equiv \M^{\A} \subset \M$ and
a crossed product extension $\M\subset \M\cros \A$.
\Definition{1.2}{{\em ($\A$-module algebras)}
\\
A $*$-algebra $\M$ together with a left action $\lef:\A \otimes \M \to \M ,\ (a\otimes m)\mapsto a\lef m,$ is
called an $\A$-module algebra, if the following properties hold for all $a,b \in \A$ and $ m,n \in \M$
\bea
(ab)\lef m & = & a\lef (b\lef m) \lb{1.29}\\
\one_\A\lef m &=& m \lb{1.30}\\
a\lef(mn)& = &(a_{(1)} \lef m)(a_{(2)}\lef n) \lb{1.31}\\
(a\lef m)^*&=& a_*\lef m^* \lb{1.32}\\
a \lef \one_{\M} &=& (a_{(1)}S(a_{(2)}))\lef\one_{\M}
\equiv\e(\one\1a)\one\2\lef\one_\M \lb{1.33}
\eea
If $\M$ is a $C^*$- or a von-Neumann algebra, then we also
require $a\lef$ to be norm or weakly continuous, respectively,
for all $a\in\A$.
}
We will see that the axiom \no{1.33} may equivalently be
replaced by
\beq\lb{1.33'}
a\lef\one_\M =( a\2S^{-1}(a\1))\lef\one_\M
\equiv\e(\one\2a)\one\1\lef\one_\M
\eeq
An example of an $\A$-module action is provided by the natural left action
$\A\o\hA\ni a\otimes \phi \mapsto (a\arr \phi )\in\hA$
 of $\A$ on its dual $\hat{\A}$, which together with the
opposite right $\A$-action
$\hA\o\A\ni\phi\o a\mapsto (\phi\arl a)\in\hA$ is defined as
for ordinary Hopf algebras, i.e.
\bea\lb{1.34}
\bra a\arr \phi \mid b\ket &:=& \bra\phi \mid ba\ket
\lb{left}
\\
\lb{right}
\bra\phi\arl a\mid b\ket &:=&\bra\phi\mid ab\ket
\eea
where $\phi \in \hat{\A}$ and $a,b \in \A$.
For the left action \no{1.34} the identities \no{1.33} and
\no{1.33'} follow from the axioms $IIb$ and
$IIIb$. As a natural cyclic $\A$-submodule we denote
$$
\hat\A_R := \A\arr\one_{\hat\A}\subset\hat\A
$$
which as an object in $\Rep\A$ is in fact equivalent to the
``trivial'' representation $\pi_\e$ of $\A$. Moreover,
(\ref{1.33}) implies (see Section 2.2)
that as a cyclic $\A$-submodule
\beq\lb{1.35}
\M_R := \A\lef\one_{\M} \subset \M \lb{M_R}
\eeq
is homomorphic to $\hA_R$  via the map
\beq\lb{1.36}
\tau_\lef:\hat\A_R\ni(a\arr\one_{\hat\A})\,\mapsto\,(a\lef \one_{\M})\in\M_R
\eeq
We will see (\lem{2.5}) that if
$\A$ is pure, (\ref{1.36}) even defines an
$\A$-module isomorphism.
Moreover, we will show in \prop{2.6}
that for any $\A$-module algebra
$\M$ the submodule $\M_R\subset\M$
also is a $*$-subalgebra of $\M$ and that the map (\ref{1.36}) also
provides a homomorphism of $*$-algebras.
Interchanging the r\^ole of $\A$ and $\hat\A$ in (\ref{left})
we may put similarly
\bea
\A_R &:=& \{ (id\otimes \phi)(\Delta(\one_\A)) \mid \phi\in\hat\A\}
\subset\A \lb{A_R}\lb{1.37}\\
\A_L &:=& \{(\phi \otimes id)(\Delta(\one_{\A})) \mid \phi \in \hat{\A}\}
\equiv S(\A_R)\lb{A_L}\ .\lb{1.38}
\eea
Consequently, $\A_{R/L}\subset\A$ are also $*$-subalgebras.
Moreover, $\A_{L/R}\cong\hA_{R/L}$ as $*$-algebras, the
isomorphism being given by [BSz,Sz,N2,BNS]
\bea\lb{1.39}
\A_L\ni a \mapsto (a\arr\one_{\hat\A}) \in \hat\A_R
\\
\A_R\ni a \mapsto (\one_{\hat\A}\arl a) \in \hat\A_L
\eea
with inverses given by the dual versions, respectively.
Composing (\ref{1.39}) with (\ref{1.36}) we get a
$*$-algebra epimorphism
\beq
\mu_{\lef} :\, \A_L\ni a \mapsto (a\lef \one_{\M}) \in \M_R
\lb{1.40}
\eeq
which becomes an isomorphism if $\A$ is pure.
Since by \no{1.11'}
$\A_L$ and $\A_R$ always commute and since one also always has
\beq\lb{1.40'}
\mu_\lef(\A_L\cap\A_R)\subset C(\M)
\eeq
(see \cor{2.17}(i)), we already arrive at the first
ingrediences expected to appear in reducible Jones towers as
discussed in our motivation in Section 1.1 (see \Eq{0.3}).

\bigskip
We now turn to what by analogy we would like to call the
fixed point subalgebra
$\N\equiv\M^\A\subset\M$, for which the following definitions
will be appropriate.

\Definition{1.4}
{{\em (Fixed point algebra)}\\
Given an $\A$-module algebra $\M$ we define the ``fixed point
algebra" $\N\equiv \M^{\A} \subset \M$ to be the unital
$*$-subalgebra given by the elements $n\in\M$ satisfying either of the
following conditions
$$
\begin{array}{rrcll}
i)\qquad & a\lef (mn) & = & (a\lef m)\,n &,\quad\forall a\in\A,\ m\in\M \\
ii)\qquad & a\lef n & = & a_{(1)}S(a_{(2)})\lef n &,\quad\forall a\in\A\\
iii)\qquad & a\lef (nm) & = & n\,(a\lef m) &,\quad\forall a\in\A,\ m\in\M \\
iv)\qquad & a\lef n & = & a_{(2)}S^{-1}(a_{(1)})\lef n &,\quad\forall
a\in\A\ .
\end{array}
$$
}
We  will show in \prop{2.11} that the conditions $(i-iv)$
are in fact all equivalent. Putting $m=\one_\M$ in condition $i)$ and
$iii)$ then implies
\beq\lb{1.41}
\M_R\subset\N'\cap\M.
\eeq
and we will later fix appropriate conditions guaranteeing equality in
\no{1.41}.
Since eventually we want to identify $A_L$ with the relative
commutant $\N'\cap\M$, we now propose the following

\Definition{1.3}
{An $\A$-module algebra $\M$ is called
{\em standard}, if $\mu_\lef$ in \no{1.40} provides an
isomorphism $\A_L\cong\M_R$.
}
Note that for ordinary Hopf algebras one always has
$ \A_L = \A_R=\M_R =\CC$ and therefore in this case
standardness trivially holds.
Moreover, condition (ii) of \defi{1.4} together with
\no{1.33} and \no{1.41} also imply
\beq\lb{1.41'}
\mu_\lef(\A_L\cap C(\A))\subset C(\N)
\eeq
and under suitable regularity conditions we will later also
have equality in \no{1.41'}, thus reproducing \no{0.1} and
\no{0.5a}.
%
%

Finally, we define the crossed product $\M \cros\A$ as an ``amalgamated"
tensor product over $\A_L$.

\Definition{1.6}
{ {\em (Crossed products)\\}
Let $\M$ be an $\A$-module $*$-algebra. The {\em
crossed product} $\M\cros\A$ is defined to be the linear space
\beq\lb{1.42}
\M\cros \A =\M\otimes_{\A_L} \A
\eeq
where $\A_L$ acts on $\A$ by left multiplication and on $\M$ by
right multiplication via its image under $\mu_{\lef}$.
The $*$-algebra structure on $\M\cros \A$ is defined in
the same way as for ordinary crossed products, i.e.
for $m,m'\in\M$ and $a,a'\in\A$
\bea
&(m\cros a)(m'\cros a') = (m(a_{(1)} \rhd m')\cros a_{(2)} a')&\lb{1.43}\\
&(m\cros a)^* =(\one_\M\cros a^*)(m^*\cros\one_\A)
\equiv ((a_{(1)}^* \rhd m^*)\cros a_{(2)}^*)&
\lb{1.44} .
\eea
}
We will show in \thm{3.1} that with
these definitions $\M\cros\A$ indeed becomes a well defined
$*$-algebra extending $\M\equiv(\M\,\cros\, \one_\A)$.
Moreover, as for crossed products by ordinary Hopf algebras we
will have
\beq\lb{1.44b}
((a\lef m)\cros \one_\A)=(\one_\M\cros a\1)\,(m\cros \one_\A)\,(\one_\M\cros S(a\2))\ .
\eeq
Confirming with the scenario in Section 1.1 we will also get
\bea
(\one_\M\,\cros\, \A) &\subset& \N'\cap\M\cros\A
\lb{1.44d}\\
(\one_\M\,\cros\, \A_R) &\subset& \M'\cap\M\cros\A
\lb{1.44c}\\
\one_\M\,\cros\, (\A_R\cap C(\A)) &\subset& C(\M\cros\A)
\lb{1.44e}\\
\one_\M\,\cros\, (\A_L\cap\A_R\cap C(\A)) &\subset& \M\cap C(\M\cros\A)
\equiv\N\cap C(\M\cros\A)
\lb{1.44f}
\eea
where equality in all these inclusions will again
be proven under suitable regularity conditions in Section 4.
This concludes the presentation of our basic setting and we now
proceed to describe our main results.

\subsection{Summary of results}

Many consequences of our axioms described so far are already
contained in [BSz,Sz,N2,BNS] and are presented here only to the
extend making this work sufficiently selfunderstood. Our main
focus here is on applications to the theory of $\A$-
(co)module algebras $\M$ and their crossed product extensions,
which have not been treated in
general in the above papers.

Section 2.1 starts with generalizing the well kown duality relation between
actions and coactions from ordinary Hopf algebras to our
setting.
In Section 2.2 we study the submodule $\A\lef\one_\M$ and prove that
it is a subalgebra of $\M$ homomorphic to $\A_L\cong\hat\A_R$.
In Section 2.3 we verify that our axioms for the fixed point
algebra in \defi{1.4} are in fact all equivalent.
In Section 2.4 we review various relations between $\A_{L/R},\
\A_L\cap\A_R,\ A_{L/R}\cap C(A)$ and their dual counterparts and
show how they interplay with $\M_R$ and $C(\M)$ in the spirit
of our discussion in Section 1.1.
In Section 2.5 we apply our formalism to the example of a
partly inner group action on a factor $\M$.

In   Section 2.6 we review and apply the theory of left integrals as
developed in [BNS]. We show that positive and normalized
left integrals $l\in\A$ give rise to conditional expectations
$E_l:\M\to\N\equiv\M^\A$ via $E_l(m):=l\lef m$. Under the
assumption $\A_L\cong\M_R=\N'\cap\M$ the correspondence
$l\leftrightarrow E_l$ will be one-to-one and $E_l$ will be
faithful if and only if $l$ is nondegenerate.
Considering the special case $\M=\hat\A$
we have $\N=\hat\A_L$ and we denote $\Ind l\in C(\hat\A)$ the
index of $E_l:\hat\A\to\hat\A_L$.
We will see that in fact
\beq\lb{1.45}
\Ind l\in C(\hat\A)\cap\hat\A_R
\eeq
for all nondegenerate left integrals $l\in\A$.
Note that under the scenario of Section 1.1 we would have
$\tau_\lef(C(\hat\A)\cap\hat\A_R)=\mu_\lef(\A_L\cap\A_R)=C(\M)$
suggesting that in the general case the index of $E_l:\M\to\N$
might be given by
\beq
\Ind E_l = \tau_\lef(\Ind l)\in C(\M)\ .
\eeq
Under suitable regularity conditions this will indeed be a
result in Section 4.

\bigskip
In Section 3 we turn to the study of crossed products. In
Section 3.1 we prove that $\M\cros\A$ is a well defined
$*$-algebra extending $\M\equiv(\M\,\cros\,\one_\A)$ and
satisfying \no{1.44b} - \no{1.44f}.
In Section 3.2 we show that as for ordinary Hopf algebras
there is a natural
$\hat\A$-module algebra structure on $\M\cros\A$.
The fixed point algebra under
this action is given by $\M\equiv(\M\,\cros\,\one_\A)$ and we have
$(\M\cros\A)_R=(\one_\M\cros \A_R)$.
In Section 3.3 we improve the duality theory for nondegenerate
left integrals (i.e. the theory of Fourier transformations) of
[BNS] to develop a new notion of {\em p-duality} for {\em positive}
nondegenerate left integrals, which is patterned after the
corresponding notion of  Haagerup duality for finite index conditional
expectations in Jones theory. Associated with any normalized
positive and nondegenerate left integral $l\in\A$ we define a
projection $e_l\in\A$ and a positive and nondegenerate left
integral $\l_l\in\hat\A$ - the p-dual of $l$ - such that
\beanon
&e_l m e_l = e_l E_l(m) = E_l(m) e_l&,\quad\forall m\in\M\
\\
&\hat E_{\l_l}(e_l) =\one_\M&
\eeanon
as identities in $\M\cros\A$.
Proceeding by alternating crossed products we also provide a
generalized Temperley-Lieb-Jones algebra.

Starting from Section 3.4 we focus our attention to the case
where $\M$ is a von-Neumann algebra. We show that in this case
also $\M\cros\A$ becomes a von-Neumann algebra by identifying
it with a (non-unital) $*$-subalgebra of $\M\o\End\A$.
In particular, associated with any
(faithful, normal) $*$-representation $\pi$ of $\M$ on a
Hilbert space $\H$, we
obtain a so-called {\em regular} (faithful, normal)
representation $\pi_{cros}$ of
$\M\cros\A$ on (a subspace of) $\H\o L^2(\A,\hat h)$,  thus
generalizing standard results for crossed products by ordinary
Kac algebras [Pe, St, ES].

\bigskip
In Section 4 we apply our general results to Jones theory,
guided by the special example of a partly inner group action.
In Section 4.1 we
show that any GNS-representation $\pi_\om$ of $\M$ associated
with a faithful normal $\A$-invariant state $\om$ on $\M$ extends to a
representation -- still denoted $\pi_\om$ -- of $\M\cros\A$ (which
may actually be identified with a subrepresentation of the
regular representation $\pi_{cros}$ mentioned above).
We then prove that the basic Jones construction
for $\M_{-1}\equiv\pi_\om(\M^\A)\subset\pi_\om(\M)\equiv\M_0$
is precisely given by $\M_1:=\pi_\om(\M\cros\A)$, which in
general need not be isomorphic to $\M\cros\A$. Rather, we have
$\M_1\cong\M e_l \M$ which might be a nontrivial ideal in
$\M\cros\A$. We also show that $\M_{-1}\subset\M_0\subset\M_1$
is always of finite index and depth 2 and that
$$
\Ind E_l \le \tau_\lef(\Ind l)
$$
for any positive nondegenerated left integral $l\in\A$, where
equality holds if and only if $\M_1$ is a {\em faithful} image of
$\M\cros\A$.
We also point out that this in particular holds, if $\M$ is
itself a crossed product, $\M=\N\cros\hat\A$, with canonical
$\A$-action.

In Section 4.2 we introduce an appropriate outerness condition
for $\A$-actions on $\M$ and show that $\A$ acts outerly iff
$\M'\cap(M\cros\A)=(C(\M)\cros\A_R)$. For outer actions of a pure weak
Hopf algebra $\A$ on a factor $\M$ this will imply that
$\M\cros\A$ is a factor and therefore
$\M\cros\A\cong\M_1$ and
\beq\lb{Jo}
\N\equiv\M^\A\subset\M\subset\M\cros\A
\eeq
is a Jones triple of factors obeying $\M'\cap(\M\cros\A)=\A_R$.
Moreover, in this case also $\hat\A$ is pure and acts outerly
on $\M\cros\A$, such that our construction iterates.
This generalizes well known facts for outer actions by groups
[NaTa] or Kac algebras [St, ES].

In Section 4.3 we envoke methods of Tomita-Takesaki theory to
prove that under the above setting $\N'\cap\M=\A_L$ and
$\N'\cap\M\cros\A=\A$, thus establishing the heuristic picture
developped in Section 1.1.

Finally, in Section 4.4 we generalize our results to
non-factors $\M$ and non-pure weak Hopf algebras $\A$ by
requiring as a regularity condition standardness and outerness
of the $\A$-action together with
$C(\M)=(\A_L\cap\A_R)\lef\one_\M$.
\footnote{In view of \no{1.40'} this means that the center of
$\M$ is required to be as small as possible.}
Again, in this case $\hat\A$ also acts regularly on $\M\cros\A$
and \no{Jo} still provides a Jones triple, where the relative
commutants are again given by $\A_L,\ \A_R$ and $\A$, respectively,
as above. Also, in this case the lower bounds in \no{1.41'},
\no{1.44e} and \no{1.44f} are saturated, i.e
\beanon
C(\N) &=& \A_L\cap C(\A)
\\
C(\M\cros\A) &=& \A_R\cap C(\A)
\\
C(\N)\cap C(\M\cros\A) &\equiv& C(\M)\cap C(\M\cros\A)
= A_L\cap\A_R\cap C(\A)
\eeanon
as anticipated in Section 1.1.
Here we have identified $\A\equiv(\one_\M\cros\A)$ and
$\A_L\equiv\A_L\lef\one_\M\equiv(\one_\M\cros\A_L)\equiv\M_R$.

\bigskip
In Appendix A we generalise the notion of Galois actions [CS]
to our setting and show that it is equivalent to
$\M^\A\subset\M\subset\M\cros\A$ being a Jones triple. In
Appendix B we analyze the appearance of cocycles for partly
inner group actions within our framework.

\bigskip
In part II of this project [NSW] we will go the opposite way,
i.e. we will show that for any finite index and depth-2 Jones
tower
$$
\N\subset\M\subset\M_1\subset\M_2\subset\dots
$$
of von-Neumann algebras with finite dimensional centers
the relative commutants $\A:=\N'\cap\M_1$ and
$\hat\A:=\M'\cap\M_2$ give a dual pair of weak Hopf algebras
acting regularly on $\M$ and $\M_1$, respectively, such that
$\N=\M^\A,\ M_1=\M\cros\A$ and $\M_2=\M_1\cros\hat\A$.
Moreover, in this case $\M$ is a factor iff $\hat\A$ is pure,
$\N$ is a factor iff $\A$ is pure and the inclusions are
irreducible iff $\A$ and $\hat\A$ are ordinary Hopf algebras.

\bsn
{\bf Note added:}
After writing up this work we were informed by L. Vainermann about a preprint by M. Enock and J.-M. Vallin [EV] treating a similar approach.



\sec {$\A$-Module Algebras}

Throughout  this section let $\A$ and $\hat\A$ be a dual pair
of finite dimensional
weak $C^*$-Hopf algebras and let $\M$ denote an $\A$-module algebra as
described in \defi{1.2}. The units in $\A,\ \hat A$ and $\M$
will always be
denoted by $\one \in\A,\ \hat\one \in\hat\A$ and $\one_\M\in\M$.
Elements of $\A$ will be
denoted by Roman letters $a,b,c...$ and elements of $\hat\A$
by Greek letters $\phi,\psi,\xi,...$ .

We will freely identify $\A=\hat{\hat\A}$ and denote the
dual pairing by
$\langle a|\phi\rangle\equiv \langle\phi|a\rangle ,\
a\in\A,\ \phi\in\hat\A$.
Let $\A_{op}$ be the weak Hopf algebra $\A$ with opposite multiplication and
$\A^{cop}$ the weak Hopf algebra $\A$ with opposite comultiplication. Then
$\A_{op},\ \A^{cop}$ and $\A^{cop}_{op}$ are again weak $C^*$-Hopf algebras,
where the antipode of $\A_{op}$ and $\A^{cop}$ is given by
$S^{-1}$ and the antipode
of $\A_{op}^{cop}$ by $S$ [N2,BNS].

A general theory for weak Hopf algebras is developped in
[BSz,Sz,N2,BNS], where in
particular equs. (\ref{1.1})-(\ref{1.8}) are proven. The
aim of this section is to
investigate general properties of $\A$-module algebras $\M$ and their
$\A$-invariant subalgebras $\N\equiv \M^\A\subset\M$ and
to establish the connection between conditional expectations
$E:\M\to\N$ and left integrals $l\in\A$.

\subsection {Coactions}

We start with reformulating \defi{1.2} in terms of an
$\hat\A$-coaction on $\M$.
\Def{2.1}{{\em ($\A$-Comodule Algebras)}\\
A (right) $\A$-comodule algebra $(\M,\rho)$ is a $*$-algebra $\M$ together with a (in general
non-unital) $*$-algebra homomorphism $\rho:\M\to \M\otimes \A$ satisfying
\bea
(\rho\otimes id_\A)\circ\rho &=& (id_\M\otimes \Delta)\circ \rho
\lb{2.1}\\
(id_\M\otimes \e )\circ\rho &=& id_\M \lb{2.2}\\
(\one_\M\otimes\Delta (\one)) (\rho(\one_\M)\otimes\one) &=&
((id_\M\otimes \Delta)\circ\rho)(\one_\M) \lb{2.3}
\eea
where $\Delta$ and $\e$ denote the coproduct and counit,
respectively, on $\A$.
If $\M$ is a von-Neumann algebra, we also require $\rho$ to be
weakly contiuous.
}
For simplicity, we call such a $\rho$ an $\A$-coaction on $\M$.
Note that since $\rho$ is assumed $*$-preserving, the axiom (\ref{2.3})
could equivalently be
replaced by
\bea\lb{2.4}
(\rho(\one_\M)\otimes\one)(\one_\M\otimes\Delta(\one)) &=&
((id_\M\otimes\Delta)
\circ \rho)(\one_\M)
\eea
Without a $*$-structure we would require both axioms, \no{2.3}
and \no{2.4}.
An immediate example of an $\A$-coaction is given by $\M=\A$
and $\rho=\Delta$.
$\hat\A$-coactions are defined correspondingly.
We now have
\Prop{2.2}
{{\em (Left Action $\leftrightarrow$ Right Coaction)}\\
There is a one-to-one correspondence between left
$\A$-module algebra actions $\lef :\A\otimes\M\to\M$ and right $\hat\A$-comodule algebra coactions $\hat\rho:\M\to\M\otimes\hat\A$
given by
\beq\lb{2.5}
a\lef m=(id_\M\otimes a)(\hat\rho(m)).
\eeq
}
{\bf Proof:} Since $\A$ is finite dimensional, (\ref{2.5})
defines a one-to-one
correspondence between $\lef$ and $\hat\rho$ as linear maps. Equ. (\ref{2.1})
is then equivalent to (\ref{1.29}),  equ. (\ref{2.2}) is
equivalent to (\ref{1.30}),
the homomorphism property of $\hat\rho$ is equivalent to (\ref{1.31}) and
equ. (\ref{1.32}) is equivalent to $\hat\rho(m^*)=\hat\rho(m)^*,\ \forall m\in\M$.
We are left to relate (\ref{1.33}) with (\ref{2.3}) (and
terefore (\ref{2.4})). By applying
$id_\M\otimes a\otimes b$ to both sides equ. (\ref{2.3}) is equivalent to
\beq \lb{2.6}
ab\lef\one_\M =\e(a\1 b)a\2\lef\one_\M,\ \forall a,b\in\A
\eeq
and equ. (\ref{2.4}) is equivalent to
\beq \lb{2.7}
ab\lef\one_\M=(a\1\lef \one_\M)\e(a\2b)~~\forall a,b\in\A
\eeq
Putting $a=\one$ equ. (\ref{2.6}) yields
\beq \lb{2.8}
b\lef\one_\M=b\1 S(b\2)\lef\one_\M
\eeq
by axiom IIIb) of \defi{1.1}. Hence (\ref{2.3}) implies
(\ref{1.33}). Conversely,
suppose (\ref{1.33}) holds together with (\ref{1.29}-{\ref{1.32}). Then
\beanon
ab\lef \one_\M &=&ab\1S(b\2) \lef \one_\M\\
&=& \e(a\1 b)a\2 \lef\one_\M
\eeanon
where we have used equ. (\ref{2.9}) of \lem{2.3} below. Thus the axioms
(\ref{1.29}) - (\ref{1.33}) imply (\ref{2.6}) and therefore,
(\ref{2.3}). \qed\\
\\
We are left to prove
\Lem{2.3}
{For $a,b\in\A$ the following identities hold
\bea
ab\1 S(b\2) &=& \e(a\1b)a\2 \lb{2.9}\\
S(b\1)b\2 a &=& a\1\e(ba\2) \lb{2.10}\\
ab\2S^{-1}(b\1) &=& \e(a\2b)a\1 \lb{2.11}\\
S^{-1}(b\2)b\1 a &=& a\2\e (ba\1) \lb{2.12}
\eea
}
{\bf Proof:} Using $\Delta(b)=\Delta(\one )\Delta(b)$ and equ.
(\ref{1.2}) we get
\beanon
ab\1 S(b\2) &=& a\one\1 b\1S(b\2)S(\one\2)\\
&=& a\1b\1 S(b\2)S(a\2)a\3\\
&=&\e(\one\1 a\1b)\one\2a\2\\
&=& \e(a\1b)a\2
\eeanon
where we have used \no{1.7} in the second line and the antipode axiom IIIb of
\defi{1.1} in the third line. Thus we have proven (\ref{2.9}).
Equs. (\ref{2.10})-(\ref{2.12}) reduce to (\ref{2.9}) in $\A_{op}^{cop},
\A^{cop}$ and $\A_{op}$, respectively. \qed
\Cor{2.4}
{Let $\lef :\A\otimes\M\to\M$ obey
(\ref{1.29})-(\ref{1.32}).
Then the axiom (\ref{1.33}) is equivalent to}
\beq\lb{2.13}
a\lef\one_\M =a\2 S^{-1}(a\1) \lef\one_\M,\quad\forall a\in\A
\eeq
{\bf Proof:} Replacing (\ref{1.33}) by (\ref{2.13}) amounts to saying that
$\lef :\A^{cop} \otimes\M_{op}\to\M_{op}$
provides a left $\A$-module algebra action. Now
$\widehat{\A^{cop}} =(\hat\A)_{op}$ and according to \defi{2.1}
$\hat\rho:\M\to\M\otimes\hat \A$ is a coaction iff
$\hat\r:\M_{op}\to M_{op}\otimes (\hat\A)_{op}$
also is a coaction. \qed

\bsn
As a warning we remark that the adjoint action of $\A$ on
itself given by $a\lef b:=a\1 bS(a\2)$ in general is {\em not} an 
$\A$-module algebra action in our sense, since it fails the axioms \no{1.30}
and \no{1.31}. More specifically, in place of \no{1.31} we have
$(a\1\lef b)(a\2\lef c) =a\lef ((\one\lef b)(\one\lef c))$,
for all $a,b,c\in\A$. However one can show that acting
by $\one$ on $\A$ via the adjoint action defines a conditional expectation
$\one\lef : \A\to\A_R'\cap\A$ and
that the adjoint action restricts to a $\A$-module algebra action on
$\A_R'\cap\A$.

\subsection {The submodules $\A\lef\one_\M$}

We now study the cyclic submodule
\beq \lb{2.14}
\M_R:=\A\lef\one_\M \subset\M
\eeq
of a given $\A$-module algebra $\M$. First we consider
$\M = \hat\A$ with
$\A$-left action (\ref{1.34}) corresponding to the $\hat\A$-coaction
$\hat\rho=\hat\Delta:\hat\A\to\hat\A\otimes\hat\A$. We note that
the action of $\A$
on the vector space $\hat\A_R\equiv(\A\arr\hat\one)
\subset\hat\A$ is equivalent to
$\pi_\e$. In fact, let $(\pi_\e,\H_\e,\Omega_\e)$ be the
GNS-triple associated
with the counit $\e$ and define $T:\hat\A_R\to\H_\e$ by
\beq\lb{2.15}
T(a\arr\hat\one) :=\pi_\e(a)\Omega_\e,\ a\in\A
\eeq
Then $T$ is immediately seen to be a well defined $\A$-linear bijection.
Next, we define $\tau_\lef:\hat\A_R\to \M_R$ by putting for $a\in\A$
\beq\lb{2.16}
\tau_\lef(a\arr\hat\one):= a\lef \one_\M\equiv a\1S(a\2)\lef\one_\M
\eeq
\Lem{2.5}
{{\rm ($\M_R\subset\M$ as an $\A$-submodule)\\}
For any $\A$-module algebra $\M$ the map $\tau_\lef$
provides a well defined $\A$-module epimorphism, which becomes
an isomorphism if $\A$ is pure or if the $\A$-action $\lef$ on
$\M$ is faithful.
}
{\bf Proof:}
$\tau_\lef$ is well defined since $a\arr\hat\one=0$ implies $\e(ba)=0$ for all
$b\in\A$ and, therefore, $a\1S(a\2)\equiv \e(\one\1a)\one\2=0$.
Clearly, $\tau_\lef$ is surjective and intertwines $\hat\A_R$ and $\M_R$ as left
$\A$-modules. If $\A$ is pure then $\hat\A_R\cong\H_\e$ is
$\A$-irreducible implying the injectivity of $\tau_\lef$. By
\prop{2.9} and \prop{2.10} below $\tau_\lef$ is also injective if
the action $\lef$ is faithful. \qed

\bsn
We now show that $\M_R\subset \M$ (and in particular
$\hat\A_R\subset\hat\A$) are in fact $*$-subalgebras
and that $\tau_\lef:\A_R\to\M_R$ is also a $*$-algebra homomorphism.

\Prop{2.6}
{{\rm ($\M_R\subset\M$ as $*$-subalgebra) \\}
Let $\M$ be an $\A$-module $*$-algebra and let
$\M_R=\A\lef \one_\M$. Then $\M_R\subset \M$ is a unital
*-subalgebra and the
$\A$-module map $\tau_\lef:\hat\A_R\to\M_R$ (\ref{2.16}) is also a
unital *-algebra homomorphism.
}
{\bf Proof:} To prove that $\M_R\subset \M$
is a unital *-subalgebra we first
note $\one_\M =\one \lef\one_\M \in \M_R$
and $(a\lef\one_\M)^* =(a_*\lef\one_\M)
\in \M_R$, for all $a\in\A$. Using
\beq \lb{2.27}
n m = (\one\1\lef n)(\one\2\lef m) ,\ \forall m,n\in \M
\eeq
we now compute for all $m\in\M$
\bea \lb{2.28}
(a\lef\one_\M)m &=& (a\1 S(a\2)\lef \one_\M)m\nonumber\\
&=& (\one\1 a\1 S(a\2)\lef \one_\M)(\one\2 \lef m)\nonumber\\
&=& (a\1 S(a\4)\lef \one_\M)(a\2S(a\3)\lef m)\nonumber\\
&=& a\1 \lef ((S(a\3)\lef \one_\M)(S(a\2)\lef m))\nonumber\\
&=& a\1\lef (S(a\2)\lef (\one_\M m))\nonumber\\
&=& a\1 S(a\2) \lef m
\eea
where in the third line we have used (\ref{1.6}), in the fourth
line (\ref{1.31}),
in the fifth line (\ref{1.3}) and (\ref{1.31}) and in the last
line (\ref{1.29}).
Putting $m =b\lef\one_\M$ we conclude
\beq\lb{2.29}
(a\lef\one_\M)(b\lef\one_\M) = a\1 S(a\2)b\lef\one_\M.
\eeq
Thus $\M_R\subset M$ is a unital $*$-subalgebra.
Since \no{2.29} holds for all $\A$-module algebras $\M$
and in particular for $\M=\hat\A$, the map
$\tau_\lef:\hat A_R\to\M_R$
(\ref{2.16}) is in fact a $*$-algebra homomorphism.
\qed

\bsn
Interchanging $\A$ with $\hat\A$ we conclude that $\A_R\subset
\A$ also becomes a unital $*$-subalgebra. Similarly
$\hat\A_L=(\hat\A^{cop})_R \subset\hat\A$ and
$\A_L=(\A^{cop})_R \subset\A$ are
unital $*$-subalgebras.
Moreover, \Eq{1.11'} and its dual version imply
\beq
\lb{2.11'}
[\A_L\,,\,\A_R]=0 \quad , \quad [\hA_L\,,\,\hA_R]=0\ .
\eeq
Let us next
introduce, for $\phi\in\hat\A$, the linear maps $\phi_L:\A\to\hat\A$ and
$\phi_R:\A\to\hat\A$ given by
\bea
\bra\phi_L(a)|b\ket &:=&\bra\phi|ab\ket \lb{2.17}\\
\bra\phi_R(a)|b\ket &:=& \bra\phi|ba\ket\lb{2.18}
\eea
Note that for all $\phi\in\hat\A$ the maps $\phi_L$ and
$\phi_R$ are transposes of each other. Also note that choosing
$\phi=\e$ this now allows to rewrite elements of $\hat\A_{R/L}$
as $\e_R(a)\equiv a\arr\hat\one$ and
$\e_L(a)\equiv\hat\one\arl a,\ a\in\A$, respectively.
Introducing similarly $\hat\e_{L/R}:\hat\A\to\A_{L/R}$,
where $\hat\e\equiv\one_\A\in\A$ is the counit on $\hat\A$, the
axioms $III$a,b for the antipode $S$ can then be rewritten as
\bea
S(a\1)a\2 &=& \hat\e_R\e_L(a) \lb{2.23}\\
a\1S(a\2) &=& \hat\e_L\e_R(a) \lb{2.24}\\
S^{-1}(a\2)a\1 &=& \hat\e_L\e_L(a) \lb{2.25}\\
a\2S^{-1} (a\1) &=& \hat\e_R\e_R(a) \lb{2.26}
\eea
where (\ref{2.25})-(\ref{2.26}) follow from (\ref{1.7})-(\ref{1.8}).
Note that these identities in particular imply
\beq
S\hat\e_R\e_R = \hat\e_L\e_R\qquad,\qquad
S\hat\e_L\e_L = \hat\e_R\e_L
\lb{2.25b}
\eeq
Next we mention the important identities [N2]
\bea
\e_\sigma\hat\e_{\sigma'}\e_\sigma &=&\e_\sigma \lb{2.19}\\
\hat\e_\sigma\e_{\sigma'}\hat\e_\sigma&=&\hat\e_\sigma \lb{2.20}
\eea
valid for all choices of $\sigma,\sigma'\in\{L,R\}$.
From these one immediately concludes that
for all $\sigma,\sigma'\in\{L,R\}$ the restrictions
$
\e_\sigma:\A_{\sigma'}\to\hat\A_\sigma
$
and
$
\hat\e_\sigma:\hat\A_{\sigma'} \to \A_\sigma
$
are bijections of vector spaces.
}
In fact, we even have
\Prop{2.9}
{{\rm [BSz,Sz,N2,BNS]
(The isomorphism $\A_{L/R}\cong\hat\A_{R/L}$)\\}
The restrictions $\mu_R\equiv\e_R|\A_L$ and
$\mu_L\equiv\e_L|\A_R$ provide $*$-algebra isomorphisms
$\mu_R:\A_L\to\hat\A_R$ and $\mu_L:\A_R\to\hat\A_L$ obeying
$\mu_R=S\circ\e_L|\A_L,\ \mu_L=S\circ\e_R|\A_R$ and
$\mu_{R/L}^{-1}=\hat\mu_{L/R}$.
}
%
\prop{2.9} immediately generalizes to arbitrary 
$\A$-module algebras $\M$ as follows

\Prop{2.10}
{For any $\A$-module $*$-algebra $\M$
the map
\beq\lb{2.31c}
\mu_\lef :\A_L \ni a\mapsto (a\lef\one_\M)\in\M_R
\eeq
provides a $*$-algebra epimorphism,
which is an isomorphism if $\A$ is pure or if
the $\A$-action $\lef$ on $\M$ is faithful.
More generally, there exists a central projection
$z_\re=z_\re^*=z_\re^2\in\A_L\cap C(\A)$ such that
$\Ker\mu_\re=z_\re\A_L$ and $\Ker\tau_\re=z_\re\arr\hA_R$.
}
\proof
The first statement follows by combining \prop{2.6} and
\prop{2.9}, since $\mu_\lef=\tau_\lef\circ\mu_R$.
In particular, $\mu_\lef$ is an isomorphism iff $\tau_\lef$ is an
isomorphism, the later being true if $\A$ is pure by \lem{2.5}.
Moreover, we will show in \cor{2.17}iii) below that
\beq\lb{2.31a}
\Ker \mu_\lef =\A_L\cap \K_\lef
\eeq
where $\K_\lef\subset\A$ is the ideal annihilated by the action
$\lef$.
\beq\lb{2.31b}
\K_\lef :=\{a\in\A\mid a\lef\M=0\}\ .
\eeq
Hence, $\mu_\lef$ (and therefore $\tau_\lef$)
is also an isomorphism if the $\A$-action
$\lef$ is faithful.
Next, identify $\hA_R\cong\H_\e$ via \no{2.15} and put
$P_\re:\hA_R\to\Ker\tau_\re\subset\hA_R$ the orthogonal
projection. By \lem{2.5} $P_\re$ is $\A$-linear and by
[Sz, Eq.(3.3)] (see also [N2, Prop.4.6]) there exists a
central projection $z_\re\in\A_L\cap C(\A)$ such that
$P_\re\phi=z_\re\arr\phi,\ \phi\in\hA_R$. \prop{2.9} then
implies $\Ker\mu_\re=z_\re\A_L$.
\qed

\bsn
Note that this (together with \cor{2.17}iii) below) also
proves the last statement of \lem{2.5}.
In most parts of this paper we will not care about whether
$\mu_\lef$ (equivalently $\tau_\lef$) is injective or not.
However, for the main results in connection with Jones theory
in Sections 4.3 and 4.4  we need this assumption, see also the
discussion in our motivation in Section 1.1.
Repeating our \defi{1.3} we
say that an $\A$-module algebra $\M$
is {\em standard} if $\mu_\lef$
(equivalently $\tau_\lef$) are isomorphisms.

\subsection{The fixed point subalgebras $\M^\A \subset\M$}

Next, we look at the ``fixed point" algebra $\N\equiv \M^\A$ and show that the
conditions i)-iv) of \defi{1.4} are in fact all equivalent.

\Prop{2.11}
{{\rm (The fixed point algebra)\\}
For $n\in\M$ the following conditions are equivalent
$$
\begin{array}{rrcll}
i)\qquad & a\lef (mn) &=& (a\lef m)n, &\quad \forall a\in\A,\ m\in\M\\
ii)\qquad & a\lef n &=& a\1S(a\2)\lef n\equiv
\e(\one\1 a)\one\2\lef n, &\quad\forall a\in\A\\
iii)\qquad & a\lef (nm) &=& n(a\lef m), &\quad\forall a\in\A,\ m\in\M\\
iv)\qquad & a\lef n &=& a\2 S^{-1} (a\1)\lef n\equiv
\e(\one\2 a)\one\1\lef n, &\quad\forall a\in\A
\end{array}
$$
}
{\bf Proof:} i) implies ii) by putting $m=\one_\M$ and using (\ref{1.33}).
Converseley, suppose ii) holds. Then for all $m\in\M,\ a\in\A$
\beanon
a\lef (mn) &=& (a\1\lef m)(a\2\lef n)\\
&=& (a\1 \lef m)(a\2 S(a\3)\lef n)\\
&=& (\one\1 a\lef m)(\one\2\lef n)\\
&=& (a\lef m)n
\eeanon
where in the third line we have used (\ref{1.5}) and in the
last line (\ref{2.27}).
Hence i) $\Leftrightarrow $ii). The same argument for
$\lef:\A^{cop} \otimes \M_{op}\to \M_{op}$ proves iii) $\Leftrightarrow $iv).
Now putting $m=\one_\M$
also give iii) $\Rightarrow $ii) by (\ref{1.33}) and i)
$\Rightarrow $iv) by (\ref{2.13}).
\qed

\bsn
We note that in terms of the coaction $\hat\rho:\M\to
\M\otimes\hat\A$ associated
with $\lef$ the condition i) of \prop{2.11} is equivalent to
\beq \lb{2.32}
\hat\rho(mn) = \hat\rho(m)(n\otimes\hat\one),\quad\forall m\in\M
\eeq
and the condition iii) is equivalent to
\beq \lb{2.33}
\hat\rho(nm)=(n\otimes\hat\one) \hat\rho(m),\quad\forall m\in\M\ .
\eeq
Clearly, it is enough to require these identities for
$m=\one_\M$.
We now define the fixed point algebra $\N\equiv\M^\A\subset\M$ to be the unital
$*$-subalgebra given by the elements $n\in\M$ satisfying one (and hence all) of
the conditions of \prop{2.11} or, equivalently,
(\ref{2.32}) or (\ref{2.33}). Putting $m=\one$ in \prop{2.11}
then shows that $\N$ commutes with $\M_R$,
\beq\lb{2.33a}
\M_R\subset\N'\cap\M
\eeq
Moreover, we have
\Lem{2.12}
{Let $\hat\rho:\M\to\M\otimes\hat\A$ be the coaction
corresponding to
$\lef :\A\otimes\M\to\M$. Then $\N\equiv\M^\A$ coincides with
\beq \lb{2.34}
\N=\hat\rho^{-1}(\M\otimes\hat\A_L)
\eeq
}
{\bf Proof:} Using \prop{2.11} ii), equ. (\ref{2.24}) and
the fact that the transpose of $\hat\e_L\e_R$ is given by
$\e_L\hat\e_R$ we get
$$
n\in\N\Leftrightarrow\hat\rho(n)=(id_\M\otimes\e_L\hat\e_R)(\hat\rho(n))
$$
and the r.h.s. implies $\hat\rho(\N)\subset\M\otimes\hat\A_L$. Conversely, if
$\hat\r(n)\in\M\otimes
\hat\A_L$ then $\hat\r(n)=(id\otimes\e_L\hat\e_R)(\hat\r(n))$ by (\ref{2.20})
and therefore $n\in\N$. \qed

\subsection
{The left and right subalgebras $\A_{L/R}\subset\A$}

In this subsection we summarize further results on the subalgebras
$\A_{L/R}\subset\A$.
The importance of these algebras in our setting has been
motivated in Section 1.1.
In particular, we prove the inclusions \no{1.40'},
$\mu_\lef(\A_L\cap\A_R)\subset C(\M)$, and \no{1.41'},
$\mu_\lef(\A_L\cap C(\A))\subset C(\N)$, see \cor{2.17}.
By \prop{2.19} this will be the reason
why in Section 4.2 we will have to restrict
ourselves to $\A$ and $\hat\A$ being {\em pure} in order to produce a
Jones tower of {\em factors} by alternating crossed products with $\A$
and $\hat\A$.
First we note
\Lem{2.15}
{{\rm [BSz,Sz,N2,BNS]}
For $a\in\A$ the following equivalences hold
$$
\begin{array}{rrcccl}
i)\qquad & a\in\A_L &\Leftrightarrow&
\Delta(a)=a\one\1\otimes\one\2 & \Leftrightarrow &
\Delta(a)=\one\1 a\otimes\one\2 \\
ii)\qquad & a\in\A_R &\Leftrightarrow &
\Delta(a)=\one\1\otimes a\one\2 &\Leftrightarrow&
\Delta (a)=\one\1\otimes\one\2 a
\end{array}
$$
}
Analogous identities hold for $\hA_L$ and $\hA_R$, which
therefore imply
\Cor{2.13}
{The fixed point algebra $\N\subset\hat\A$ under the canonical
$\A$-left action is given by $\N=\hat\A_L$.
}
Also note, that \lem{2.15} and the axiom $Id)$ of \defi{1.1}
imply
\bleq{2.37'}
\Del(\one)\in\A_R\o\A_L.
\eeq
As an application we now show that the left actions of $\A_{L/R}$
on any $\A$-module algebra $\M$ are always inner in the
following sense.
\Lem{2.16}
{For all $m\in\M$ we have\\
i) $a\in\A_L\Rightarrow a\lef m=(a\lef\one_\M)m=(S^{-1}(a)\lef\one_\M)m$
\\
ii) $a\in\A_R \Rightarrow a\lef m=m(a\lef\one_\M)=m(S(a)\lef\one_\M)$
\\
iii) $a\in\A_L\cap\A_R \Rightarrow (a\lef\one_\M)\in C(\M)$
\\
iv) $a\in\C(\A)\Rightarrow
(a\lef\one_\M)\in\M_R\cap\M^\A\subset\C(\M^\A)$
}
{\bf Proof:} If $a\in\A_L$ then
$$
a\lef m = (a\1 \lef\one_\M)(a\2\lef m)
= (\one\1a\lef\one_\M)(\one\2\lef m)
= (a\lef\one_\M)m
$$
where in the second equation we have used \lem{2.15}i) and in
the third one (\ref{2.27}). Moreover,
if $a\in\A_L$ then $a=\hat\e_L \e_L(a)=S^{-1}(a\2)a\1$ by
\no{2.20}
and (\ref{2.25}).
Therefore, using (\ref{1.33}) and (\ref{1.3})
$$
S^{-1} (a)\lef\one_\M=S^{-1} (a\2)a\1\lef\one_\M=a\lef\one_\M\ .
$$
Part ii) follows since $\A_R=(\A^{cop})_L$ and since
$\lef:\A^{cop}\otimes \M_{op}\to\M_{op}$ is also a left $\A$-module algebra
action. Part iii) is obvious from i) and ii).
Finally, part (iv) follows from \no{2.33a} and the fact, that
$\C(\A)\lef\one_\M\subset\M^\A$ by \prop{2.11}ii) and \Eq{1.33}.
\qed
\Cor{2.17}
{Let $\mu_\lef$ be as in \prop{2.10} and $\K_\lef$ as in
equ. \no{2.31b}. Then
\\
i) $\mu_\lef(\A_L\cap\A_R)\subset C(\M)$
\\
ii) $\mu_\lef(\A_{L}\cap C(\A))=C(\A)\lef\one_\M
\subset C(\M^\A)$
\\
iii) $\Ker\mu_\lef = \A_L\cap\K_\lef$
}
\proof
It remains to proof the first identity of part (ii). Since
\no{1.33} implies $a\lef\one_\M=\hat\e_L\e_R(a)\lef\one_\M$,
this follows from the dual of \Eq{2.35} below.
\qed

\bsn
\cor{2.17}iii) concludes the proofs of the statements in
\lem{2.5} and \prop{2.10} referring to the faithfulness of the
action $\lef$.

Next, we look at the intersection $\A_L\cap\A_R$
and note that it is isomorphic to
$C(\hat\A)\cap\hat\A_{L/R}$ and
that it is trivial if and only if $\hat\A$ is pure.
\Lem{2.18}
{{\rm [Sz,N2]}
For $\sigma\in\{\mbox{L,R}\}$ the $*$-isomorphisms
$\e_L|\A_{R}$ and $\e_R|\A_{L}$ of \prop{2.9} restrict to $*$-isomorphisms
$\e_\sigma: \A_L\cap\A_R \to C(\hat\A)\cap\hat\A_\sigma$
with inverse given by
$\hat\e_{L}|(C(\hat\A)\cap\hat\A_{\sigma})
\equiv\hat\e_{R}|(C(\hat\A)\cap\hat\A_{\sigma})$.
More generally, $\hat\e_L|(C(\hat\A))\equiv\hat\e_R|(C(\hat\A))$ is a
$*$-algebra homomorphism invariant under $S$ and
\bleq{2.35}
\hat\e_{\sigma}(C(\hat\A))=\hat\e_{\sigma}(C(\hat\A)\cap\hat\A_{L/R})
=\A_L\cap\A_R
\eeq
}
%
Using the identity $\pi_\e(\A)'=\pi_\e(\A_L\cap C(\A))$
([Sz, Eq.(3.3)], see also [N2, Prop.4.6]), where the prime
denotes the commutant in $\End\H_\e$, \lem{2.18} leads to

\Prop{2.19}
{{\rm [Sz] (Pureness of $\A$)\\}
For a weak $C^*$-Hopf algebra $\A$ the following conditions
are equivalent
\\
i) $\A$ is pure
\\
ii) $\hat\A_L\cap\hat\A_R = \CC\cdot\hat\one$
\\
iii) $\A_\sigma\cap C(\A) = \CC\cdot\one,\quad\sigma=L
\ \mbox{or}\  \sigma=R\ .$
}
%
Note that \cor{2.17} and \prop{2.19} imply that for
standard weak Hopf actions $\A$ is pure iff $\M^\A\cap\M_R=\CC$
and $\hat\A$ is pure iff $C(\M)\cap\M_R=\CC$.
The first condition in particular holds if $C(\M^\A)=\CC$ (by
\no{2.33a}) and the second if $C(\M)=\CC$.
In Section 4.2 this will be the reason why in applications to
depth-2 inclusions of factors we will have to restrict
ourselves to pure weak Hopf algebras $\A$ and $\hat\A$.

We finally mention that \lem{2.18} implies [Sz]
\beq\lb{Hypcen1}
\A_L\cap\A_R\cap C(\A)\cong\hat\A_L\cap\hat\A_R\cap C(\hat\A)
\eeq
and \cor{2.17} implies
\beq\lb{Hypcen2}
\mu_\lef(\A_L\cap\A_R\cap C(\A))\subset\N\cap C(\M)\ .
\eeq
The abelian algebra $\Z:=\A_L\cap\A_R\cap C(\A)$ has been called the
``hyper-center" of $\A$ in [Sz], since it behaves like scalars
in many respects.
More presicely one has a direct sum decomposition of weak Hopf
algebras [Sz,BNS]
\beq\lb{Hypcen3}
\A=\oplus_p \A_p\qquad , \qquad \hat\A=\oplus_p\hat\A_p
\eeq
where $\A_p\equiv p\A$ and $\hat\A_p\equiv \hat p\hat\A$ are
dual to each other and where
$p$ runs through the minimal projections of $\Z\cong\hat\Z$.
Correspondingly, we could also cut $\N$ and $\M$ with the
central projections $p\lef\one_\M\in\N\cap C(\M)$ to obtain
$\N_p\equiv(p\lef\one_\M)\N$ as the fixed point algebra of
$\M_p\equiv(p\lef\one_\M)\M$ under the action of $\A_p$. Under
the regularity conditions of Section 4.4 this would also imply
$\N_p\cap C(\M_p)=\CC$. Thus, throughout one might without much
loss assume trivial hypercenter.

\subsection{A partly inner group action}

We now analyze the example of a partly inner group action
within our framework. Let $G$ be a finite group,
$H\subset G$ a normal subgroup and put
\bleq{2.5.1}
\A:=\CC H\crosAd G
\eeq
where $G$ acts on $H$ by the adjoint action. Equivalently, $\A$
is the group algebra of the semi-direct product $H\crosAd G$.
Denote $\{(h,g)\mid h\in H,\,g\in G\}$ the natural basis in
$\A$ and define a weak Hopf algebra structure on $\A$ by (see
[N2, Ex. 4])
\bea
\lb{2.5.2}
\Delta (h,g) &:=&
{1\over |H|}\sum_{\tilde h\in H}
(h \,\tilde h^{-1},\ \tilde h g) \o (\tilde h,g)
\\
\lb{2.5.3}
\e (h,g) &:=& |H| \delta (h)
\\
\lb{2.5.4}
S(h,g) &:=& (g^{-1} h\,g,\, g^{-1} h^{-1})
\eea
where $\delta(\one)=1$ and $\delta(h)=0$ else.
One easily verifies the axioms of \defi{1.1}.
As to the antipode axioms we remark that
for $x=(h,g)\in\A$ one computes
\bea
\lb{2.5.5} \e(\one\1x)\one\2 &=& (h,\one) = x\1 S(x\2)\\
\lb{2.5.6} \one\1\e (x\one\2)
&=& (g^{-1} h\,g,\, g^{-1} h^{-1} g) = S(x\1)x\2
\eea
The left and right subalgebras are given by
\bea\lb{2.5.7}
\A_L&=& \span\{(h,\one_G)\mid h\in H\} \cong \CC H
\\
\lb{2.5.8}
\A_R&=& \span\{(h,h^{-1})\mid h\in H\} \cong (\CC H)_{op}
\eea
Clearly, this implies
$
\A_L\cap \A_R =\CC,
$
i.e. the dual algebra $\hat\A$ is pure.
As a linear space $\hA$ is naturally given by the functions
$\phi:\H\x G\to\CC$, with associative multiplication law
\bleq{2.2.5a}
(\phi * \psi)(h,\,g)={1\over |H|}\sum_{\tilde h\in H}
\phi(h\tilde h^{-1},\,\tilde hg)\psi(\tilde h,\,g)\ .
\eeq
Hence, as
an algebra $\hat\A$ may be identified with
$\CC H\>cros \hat\G$, where $\hat \G$
is the abelian algebra of functions on $G$
(i.e. the dual of the group algebra $\G\equiv \CC G)$, and
where the left action
of $\hat \G$ on $\CC H$ is dual to the right coaction
$
\Delta_\G:\CC H \to \CC H \o \G
$
given on the basis $h\in H$ by the standard formula $\Delta_\G(h)=h\o h$.
The coproduct on $\hA$ is given by
\bleq{2.5.2b}
(\hat\Del\phi)((h_1,\,g_1),\,(h_2,\,g_2))=\phi(h_1g_1h_2g_1^{-1},\,g_1g_2).
\eeq
Computing $\hat\Del\e$ we conclude that $\hA_L$ is given by the
space of functions of the form $\phi(h,\,g)=\varphi(h)$ and
$\hA_R$ is given by functions of the form
$\phi(h,\,g)=\varphi(g^{-1}hg)$.
Hence, $\hA_L\cap\hA_R$ is spanned by the $\Ad_G$ -invariant functions
on $H$.
By \no{2.5.2} the left action of $\hA$ on $\A$ becomes
\bleq{2.5.2c}
\phi\arr(h,\,g)={1\over|H|}\sum_{\tilde h\in H}
\phi(\tilde h,\,g)\,(h\tilde h^{-1},\,\tilde hg).
\eeq
Using $(\hA_L\cap\hA_R)\arr\one_\A=C(\A)\cap\A_R$ by \lem{2.18} we
conclude
\bleq{2.5.2d}
C(\A)\cap\A_R=\span_\chi\{\sum_{h\in H}\chi(h)\,(h,h^{-1})\}
\eeq
where $\chi$ runs through the $\Ad_G$ -invariant functions on $H$.
Hence, by \prop{2.19},  $\A$ is pure if and only if $H$ is trivial.
We leave the remaining details to the reader.

\bsn
Let now $\M$ be a von-Neumann factor, $u:H\to \M$ a unitary
representation and $\al:G\to\Aut\M$ an action satisfying
\bea\lb{2.5.11}
\al_h &=&\Ad u(h),\quad h\in H
\\\lb{2.5.12}
\al_g\circ u  &=& u\circ\Ad g.
\eea
Then $(\al,u)$ provides an $\A$-module algebra action $\lef:\A\o\M\to\M$ by
putting for $h\in H,\ g\in G$ and $m\in\M$
\bleq{2.5.13}
(h,g)\lef m:= u(h)\al_g(m)
\eeq
which the reader is invited to check.
Since $(h,g)\lef \one_M =u(h)$ the action is standard, iff
\bleq{2.5.14}
\M_R\equiv\span\{u(h)\mid h\in H\}\cong\CC H
\eeq
i.e. iff the unitaries $u(h)$ are linearly independent in $\M$.
By \no{2.5.5} the fixed point subalgebra $\M^\A$ coincides with
the fixed point algebra of the action $\al$.
\bleq{2.5.15}
\N:=\M^\A=\{m\in\M\mid\al_g(m)=m,\ \forall g \in G\}\equiv\M^G.
\eeq
We remark that for general non-outer group actions $\al$ the
inner part $u:H\to\M$ might only be a projective
representation and \no{2.5.12} might also only be valid up to a
phase. This case is treated in Appendix B.

\subsection{Left integrals and conditional expectations}

In this subsection we study the relation between
left integrals $ l \in\A$ and conditional expectations
$E:\M\to\M^\A$.
To this end
we first review from
[BNS] some central results on the theory of integrals and
Fourier transformations
on weak Hopf algebras. Let us denote $\L(\A)$ and $\R(\A)$ the space of left
and right integrals in $\A$, respectively. Then we have
\Lem{2.20}
{{\rm [BNS] (Left and right integrals)\\}
i) $ l \in \L(\A) \Leftrightarrow  l\arr \phi \in\hat\A_L$ for all
$\phi \in\hat\A$.\\
ii) $l\in\L(\A)\Rightarrow al=S(a)l$ for all $a\in\A_R$\\
i') $r\in \R(\A)\Leftrightarrow \phi\arl\in\hat\A_R$ for all
$\phi\in\hat\A$.\\
ii') $r\in \R(\A)\Rightarrow ra=rS(a)$ for all $a\in\A_L$
}
Note that \lem{2.20}i)+i') in particular implies that
$\e_R( l)\equiv l\arr\hat\one$
and $\e_L(r)\equiv\hat\one\arl r$ are both in the abelian
$*$-subalgebra $\hat\A_L\cap\hat\A_R\subset\hat\A$.
We call a left integral $ l\in \L(\A)$ (right integral $r \in \R(\A))$
{\em normalizable} if $\e_R( l)\ (\e_L(r))$ is invertible and
we call it {\em normalized} if $\e_R( l)=\hat\one$
($\e_L(r)=\hat\one$).
Normalizable left (right) integrals can always be normalized as follows.
For $l\in\L(\A)$ or $r\in\R(\A)$ and $\sigma=L,R$ we define
their ``$\sigma$-normalization''
\beq\lb{n_s1}
n_\sigma(l):=\hat\e_\sigma\e_R(l)\quad ,\quad
n_\sigma(r):=\hat\e_\sigma\e_L(r)
\eeq
Then $n_\s(l),\,n_\s(r)\in\A_\s\cap C(\A)$ by \lem{2.18} and
one may show that the two
definitions in (\ref{n_s1}) coincide on $\L(\A)\cap\R(\A)$.
Moreover
\beq\lb{n_s2}
n_L = S^{\pm 1}\circ n_R\
\eeq
and one has
\Lem{2.21}
{{\rm [BNS] (Normalizable integrals)\\}
Left integrals $l\in\L(\A)$ (right integrals $r\in\R(\A)$)
are normalizable iff $n_\sigma(l),\,n_\sigma(r)$
are invertible, respectively. In this case we have
$n_L(l)^{-1}l=n_R(l)^{-1}l$
and $r\,n_L(r)^{-1}=r\,n_R(r)^{-1}$, which provide normalized left
and right integrals, respectively.
}
Note that by \no{1.19},\no{1.20}
\beq\lb{l^2}
l^2 = n_\s(l)l,\ \s=L,R\ .
\eeq
Applying $\e_R$ to both sides we conclude that $l$ is idempotent iff
$l\arr\hat\one$ is idempotent.
By \prop{2.19} it follows, that in pure Hopf algebras nonzero left
and right integrals are normalized (normalizable), if and only if they are
idempotent (non-nilpotent), respectively.

Next, we have a theory of ``Radon-Nikodym derivatives" for left
(right) integrals.
\Lem{2.22}
{{\rm [BNS] (``Radon-Nikodym derivatives")\\}
Let $\A_\s,\ \s=L,R$, act on $\L(\A)$ by right
multiplication and on $\R(\A)$ by left multiplication.
Then under these actions the normalized Haar integral
$h$ is cyclic and separating.
More precisely, defining
for $ l\in \L(\A)$ and $r\in\R(\A)$ ``Radon-Nikodym
derivatives" with respect to $h$ by
\beq\lb{2.36}
d_\sigma( l):=\hat\e_\sigma\e_L(l)\quad,\quad
d_\sigma(r):=\hat\e_\sigma\e_R(r)
\eeq
we have
$l = hd_\sigma( l),\ r = d_\sigma(r)h$, and
$d_\s(ha)=d_\s(ah)=a,\ \forall a\in\A_\s$. For two-sided
integrals the two definitions of $d_\s$ in (\ref{2.36}) coincide and
$\L(\A)\cap\R(\A) = (\A_\s\cap C(\A))\,h$. In particular, two-sided
integrals are invariant under the antipode.
}
\lem{2.22} follows easily by noting that $h$ is a left-unit in the right ideal $\L(\A)$ and a right-unit in the left-ideal $\R(\A)$.
We also note that under a rescaling of left integrals $l$
by central elements $c\in C(\A)$
one has
\bea\lb{n_s3}
n_\s(cl)&=&(\hat\e_\s\e_{\s'})(c)n_\s(l)
\\
\lb{d_s}
d_\s(cl)&=&(\hat\e_\s\e_{\s'})(c)d_\s(l)\ ,
\eea
and an analogous formula for right
integrals.

Next, we report from [BNS] an important duality
concept for nondegenerate left integrals, which in fact should
be viewed as an appropriate generalization
of the notion of Fourier transformation to weak Hopf algebras.
We define for a left integral
$l \in \L(\A)$ the maps $ l_\sigma:\hat \A\to \A,\
\sigma =L,R$, similarly as in \no{2.17}, \no{2.18}, i.e.
by putting for $\psi\in\hat\A$
\beq\lb{2.41'}
\ba{rcl}
 l_L(\psi) &:=&  l\arl\psi\\
 l_R(\psi) &:=& \psi \arr l
\ea
\eeq
Then $ l$ is nondegenerate as a functional on $\hat\A$ iff
$ l_L$ (equivalently
$ l_R$) are invertible.
It turns out that - as for ordinary Hopf algebras [LS] - a
left integral $ l\in\A$ is nondegenerate if and only if
$\one\in\A$ is in the
image of $ l_R$.
\Prop{2.23}
{{\rm [BNS] (Duality for left integrals)\\}
Let $ l \in \L(\A)$ be a left integral and assume
$\lambda\in\hat A$ to satisfy $ l_R(\lambda)=\one$. Then $\lambda$ is a
uniquely determined left integral in $\hat \A$ and satisfies
$\lambda_R( l)= \hat\one$. Moreover, $\lambda$ and $ l$ are
both nondegenerate with
\bea
\lb{2.37}
 l_R^{-1} &=& \lambda_L \circ S^{-1}\\
\lb{2.38}
\lambda_R^{-1} &=&  l_L\circ\hat S^{-1}
\eea
}
According to \prop{2.23} nondegenerate
left integrals always come in dual
pairs. Similar statements with different left-right conventions are obtained
by passing to
$\A_{op}$, $\A^{op}$ and $\A^{cop}_{op}$, respectively. We will
see later that this kind of
duality is closely related to the notion of Haagerup duality for
conditional expectations in Jones theory.

In weak Hopf algebras the modular theory for the
normalized Haar integral turns out to be non-trivial even in
the finite dimensional case and is in fact very similar to
Woronowicz's results for the Haar state on non-finite compact quantum
groups.

\Thm{2.24}
{{\rm [BNS] (Modular theory for the Haar integral)\\}
Let $h\in\A$ and $\hat h\in\hat\A$ be the normalized Haar
integrals. Then \\
i) $h$ and $\hat h$ are nondegenerate positive
functionals on $\hat\A$ and $\A$, respectively, and
$(\hat h\arr h)\in\A_L,\ (h\arl\hat h)\in\A_R$ are positive
and invertible elements. \\
ii) Let $g_L:=(\hat h\arr h)^{1/2},
\  g_R:=(h\leftarrow\hat h)^{1/2}$ and $g:=g_Lg_R^{-1}$, and let $\hat
g_{L/R},\,\hat g$ be defined analogously in $\hat\A$. Then $g_\s\in\A_\s,\
\hat g_\s\in\hat\A_\s$ and
\bea
\hat g_\s &=& \e_\s(g_{\s'})\quad \forall \s,\s'\in\{L,R\}
\lb{2.44a}\\
g_R &=& S^{\pm 1}(g_L)
\lb{2.44c}\\
S^2(a) &=& gag^{-1},\quad \forall a\in\A
\lb{2.44d}\\
\Delta(g) &=& (g\o g)\Delta(\one) = \Delta(\one)(g\o g)
\lb{2.44e}\\
\Delta_{op}(h) &=& (\one\o g)\Delta(h)(\one\o g)
\lb{2.44f}
\eea
and analogous statements with $\A$ and $\hat\A$ interchanged.
}
Note that \thm{2.24} in particular implies that the left integral
$\lambda$ dual to $h$ is given by
\beq\lb{2.44b}
\lambda = \hat hg_L^{-2} = \hat hg_R^{-2}
\eeq
In the example $\A=\CC H\crosAd G$ of Sect. 2.5 the
normalized Haar integral is given by
\beq\lb{Haar}
e_{Haar} = \frac{1}{|G|}\sum_{g\in G}(\one_H,g)
\eeq
and therefore a basis in the space of left integrals is given
by the elements%
\footnote
{We trust that the reader will not be confused by the fact that in {\em
this} example $h\in H$ denotes the elements of the normal
subgroup $H\subset G$, whereas in {\em general} $h\in\A$ denotes the
normalized Haar integral.}
\bleq{2.6.1}
l_h:=e_{Haar}h =\frac{1}{|G|}\sum_{g\in G}(ghg^{-1},g),\quad h\in H.
\eeq
Thus, a general left integral is of the form
\bleq{2.6.2}
l=\sum_{h\in H}c(h)\,l_h
\eeq
and using \no{2.5.3} $l$ is normalized iff the
coefficients $c(h)\in\CC$ satisfy $c(\one_H)=1$ and
$\sum_{g\in G}c(ghg^{-1})=0,\ \forall\one_H\neq h\in H$.
The left integrals $\l\in\L(\hA)$ are of the form
\bleq{2.6.3}
\l(h,\,g)=\delta(hg)\hat c(h),
\eeq
and $\l\in\L(\hA)$ is dual to $l\in\A$ iff
$$
{1\over|G||H|}\sum_{\tilde h\in H}c(\tilde hh)\hat c(\tilde h)=\delta(h).
$$
In this example $e_{Haar}$ is a trace on $\hA$, i.e.
$\Del(e_{Haar})=\Del_{op}(e_{Haar})$, and the normalized Haar
integral $\l_{Haar}\in\hA$ is given by
\bleq{2.6.4}
\l_{Haar} (h,g) :=|H|\delta(h)\delta(g),
\quad h\in H,\,g\in G,
\eeq
which is also a trace on $\A$. In particular, in this case the
structural elements $g_{L/R}\in\A_{L/R}$ of \thm{2.24}
are just multiples of
the identity as for ordinary groups,
$g_{L/R}=|G|^{-1/2}\one_\A$, which fits the identity $S^2=\id$ in this
example.

\bsn
Coming back to the general theory we now
prepare the relation between left integrals and conditional
expectations by showing that
the nondegeneracy and/or positiviy of a left integral
$l\in\L(\A)$ shows up in its ``$\A_R$-Radon-Nikodym derivative"
$d_R(l)$ as follows.

\Prop{2.25}
{{\rm (Positivity and non-degeneracy for left integrals)\\}
i) A left integral $l\in\L(\A)$ is non-degenerate and/or positive
if and only if $d_R(l)$ is invertible and/or positive, respectively, as an
element in $\A_R$.\\
ii) If $l$ is nondegenerate, then the dual left integral
$\lambda\in\L(\hat\A)$ satisfies
\beq\lb{2.45}
d_R(\l) = \e_R(g_L d_L(l) g_L)^{-1}
\eeq
where $g_L:=(\hat h\arr h)^{1/2}$ as in \thm{2.24}.
In particular, $\lambda$ is positive if and only if $d_L(l)>0$ as an
element in $\A_L$.\\
iii) A nondegenerate left integral $l$ is normalizable iff $l^2$ is
nondegenerate and it is normalized iff $l^2=l$.
\\
iv) If $l\in\L(\A)$ is nondegenerate and positive, then
$n_\s(l)\in\A_\s\cap C(\A)$ is positive and invertible and
$\dot{l}:=n_\s(l)^{-1}l$ is a normalized nondegenerate and
positive left integral
}
\proof
i) Using $l=hd_R(l)$ and \lem{2.16}ii) we have for $\phi\in\hat\A$
$$
\bra l\mid\phi\ket = \bra h\mid\phi\e_R(d_R(l))\ket
$$
Since $h$ is non-degenerate, we conclude that $l$ is
non-degenerate if and only if $\e_R(d_R(l))$ is invertible,
which by \prop{2.9} is equivalent to $d_R(l)$ being
invertible. Now assume $d_R(l)\ge 0$. Then using
$S^{-1}(\A_R)=\A_L$ and \lem{2.20}ii')
we conclude
$l = hd_R(l)=  hS^{-1}(d_R(l)^{1/2})d_R(l)^{1/2}$
and therefore, again by \lem{2.16},
\bea
\bra l\mid\phi\ket &=&
\bra h\mid\e_R\left(S^{-1}(d_R(l)^{1/2})\right)\,
\phi\,\e_R(d_R(l)^{1/2})\ket\nonumber\\
&=&\bra h\mid\xi_l^*\,\phi\,\xi_l\ket\lb{2.40}
\eea
where we have introduced
$\xi_l:=d_R(l)^{1/2}\arr\hat\one\equiv \e_R(d_R(l)^{1/2})$
implying
$$
\xi_l^*=(d_R(l)^{1/2})_*\arr\hat\one=\e_R(S^{-1}(d_R(l)^{1/2}))
$$
by \no{1.16}, \no{1.17} and the self-adjointness of
$d_R(l)^{1/2}$. Hence, since $h$ is positive as a functional on
$\hat\A$ by \thm{2.24}, \no{2.40} implies the same for $l$.
We are left to show that converseley, if $l$ is a positive
weight on $\hat\A$ then $d_R(l)\ge 0$. To this end let
$B=\sum_i x_i\o y_i\in\A\o\A$ be given by $B:=S(l\1)\o l\2$.
Identifying $\A\o\A$ with the space of sesquilinear forms on
$\hat\A$ via
$$
\bra\phi\,,\,\psi\ket_B := \sum_i\overline{\bra\phi\mid
x_i^*\ket}\bra\psi\mid y_i\ket,\quad\phi,\psi\in\hat\A,
$$
we conclude from \no{1.16} and \no{1.17} that $l$ is a positive
weight on $\hat\A$ if and only if $B=S(l\1)\o l\2\in\A\o\A$
defines a positive (possibly degenerate) sesquilinear form
$\bra\cdot,\cdot\ket_B$ on $\hat\A$. By standard
diagonalization procedures this implies $B$ to be of the form
$$
B=\sum_i a_i^*\o a_i
$$
for suitable elements $a_i\in\A$. Hence
$$
d_R(l)\equiv S(l\1)l\2 = \sum_i a_i^*a_i\ge 0
$$
To prove part ii) let $l$ be nondegenerate, then $d_R(l)$ is
invertible by i). Putting $\lambda=\lambda_h\hat
S^{-1}(\e_L(d_R(l)^{-1}))$, where $\lambda_h\in\L(\hat\A)$ is
the left integral dual to $h$, we get by \lem{2.16}ii) and
\no{2.20}
$$
\lambda\arr l = \lambda_h\arr(ld_R(l)^{-1}) = \lambda_h\arr h =\one\ .
$$
Hence $\lambda$ is the left integral dual to $l$.
Now we use the identities $d_R(l)=S(d_L(l)),\ \hat
S^{-1}\circ\e_L = \e_R\circ S$ and the fact that
\no{2.44d} implies $S^2(d_L(l))=g_Ld_L(l)g_L^{-1}$ (since $g_R\in\A_R$
commutes with $\A_L$)to conclude
\beanon
\lambda=\lambda_h\e_R(g_Ld_L(l)^{-1}g_L^{-1}) &=&
\lambda_h\e_R(g_L^2)\e_R(g_L d_L(l) g_L)^{-1}\\
&=&\hat h\e_R(g_Ld_L(l)g_L)^{-1}
\eeanon
where we have used \no{2.44b}. This proves \no{2.45} and
since by \prop{2.9}
$\e_R: \A_L\to\hat\A_R$ is a $*$-algebra isomorphism we
conclude that $d_R(\l)$ is positive if and only if $d_L(l)$ is
positive.
To prove part iii) let $d_R(l)$ be invertible. By \no{l^2} we have
$d_\s(l^2)=d_\s(l)n_\s(l)$ and hence $ d_R(l^2)$ is invertible iff
$n_R(l)$ is invertible and $ l^2=l\Leftrightarrow
d_R(l^2)=d_R(l)\Leftrightarrow n_R(l) =1$.
Finally, to prove part iv) let $l\in\L(\A)$ provide a positive
faithful weight on $\hat\A$. Then $l\arr\hat\one\equiv(id\o
l)(\hat\Delta(\hat\one))\in\hat\A_L\cap\hat\A_R$ is positive.
If $l\arr\hat\one$ was not invertible there would exist a projection
$\chi=\chi^2=\chi^*\in\hat\A_L\cap\hat\A_R$ such that
$(l\arr\hat\one)\chi=0$. Consequently, $(id\o
l)(\hat\Delta(\chi))\equiv(l\arr\hat\one)\chi=0$ and therefore $\chi=0$
by the faithfulnes of $l$ and the injectivity of $\Delta$.
Thus  $l\arr\hat\one$ and by Propositon 2.9 also
$n_\s(l)\equiv\hat\e_\s(l\arr\hat\one)$ are positive and invertible.
Moreover, $\dot l = h d_R(l)n_R(l)^{-1}$ implies $d_R(\dot
l)=d_R(l)n_R(l)^{-1}$ which is invertible and positive since $n_R(l)$
is central. Hence, by part i) $\dot l$ is nondegenerate and positive.
\qed

\bsn
Using $d_R(l)=S(d_L(l))$ by \no{2.25b}, \prop{2.25}i) implies
that the left integrals \no{2.6.2} of our example are
nondegenerate iff
$
d_L(l)\equiv\sum_{h\in H}c(h)\,(h,\,\one_G)\in\CC H
$
is invertible, and they are positive iff
$
d_R(l)\equiv\sum_{h\in H}c(h)\,(h,\,h^{-1})\in\CC H\crosAd G
$
is a positive element.

A straightforward generalization
of the above methods to arbitrary $\A$-module algebra actions now gives

\Thm{2.26}
{{\rm (Left integrals and conditional expectations)\\}
Let $(M,\lef)$ be an $\A$-module von-Neumann algebra
and denote $\N\equiv \M^\A$ the $\A$-invariant subalgebra.
For $ l\in \L(\A)$
a left integral denote $E_l(m):= l\lef m$. Then \\
i) $E_l(\M)\subset\N$ and $E_l$ is a
weakly continuous $\N$-$\N$ bimodule map.\\
ii) For $d_{L/R}\in\A_{L/R}$ we have
$E_{ld_L}(m)=E_l((d_L\lef\one_\M)m)$ and
$E_{ld_R}(m)=E_l(m(d_R\lef\one_\M))$\\
iii) $E_l(\one_\M)=n_\s(l)\lef\one_\M \in C(\N)\cap\M_R$.
In particular, if $ l$ is normalizable then $E_l(\one_\M)$ is
invertible and if $l$ is
normalized then $E_l(\one_\M)=\one_\M$ and therefore $E_l(\M)=\N$.\\
iv) If $ l$ is a positive functional on $\hat\A$ then $E_l$ is positive
and
\beq\lb{2.41}
E_l(m)=E_h(z_l^*mz_l),
\eeq
where $z_l=d_R(l)^{1/2}\lef\one_\M$.
If in addition $ l$ is nondegenerate then $E_l$ is faithful.\\
v) If the $\A$-action on $\M$ is standard and $\N'\cap\M=\M_R$
then $ l \mapsto E_l$ provides a linear
bijection from $\L(\A)$ to $Hom ({_\N}\M_\N \to {_\N} \N_\N)$.
In this case $E_l$ is normalized (nondegenerate, positive) if and only if
$ l$ is normalized (nondegenerate, positive).
}
{\bf Proof:}
$E_l= l \lef$ is weakly continuous by definition. The
$\N$-$\N$ bimodule property follows from
\prop{2.11}. Also $a\lef E_l(m)=
a\1 S(a\2)\lef E_l (m)$ for all $a\in\A$ and $\M\in m$ by \no{1.19}. Hence
$E_l(\M)\subset\N$ by \prop{2.11}ii), which proves part i).
Part (ii) follows from \lem{2.16}i)+ii).
Part iii) follows from \no{1.33}, \no{2.13} and the identities
\no{2.24} and \no{2.26}. Moreover, $E_l(\one_\M)\in C(\N)$ by i),
and by \prop{2.10} $E_l(\one_\M)$ is invertible if
$n_L(l)$ is invertible.
Next, if $ l$ provides a positive (and faithful) weight on $\hat\A$ then
$E_l=(id\otimes  l)\circ\rho$ is positive (and faithful), since by \no{2.2}
$\rho$ is an injective $*$-algebra map. Moreover, \no{2.41}
follows similarly as \no{2.40} (note that $z_l=\tau_\lef(\xi_l)$).
This proves part iv).
To prove v) we use that iii), i) and the fact that
the normalized Haar integral $h\in\A$ is a positive faithful state on
$\hat\A$ [BNS] imply
$E_h:\M\to\N$ to be
a faithful conditional expectation. Hence any weakly continuous
$\N-\N$ bimodule map $E:\M\to\N$ is of the form $E(m)=E^z_h(m):=E_h(zm)$, for
some $z\in \N'\cap \M$%
\footnote
{Using a Pimsner-Popa basis $u_i$ for $E_h$
we have $z=\sum E(u_i)u_i^*$ since $E(m)=\sum E(u_iE_h(u_i^* m))
=\sum E(u_i)E_h(u_i^*m)=\sum
E_h(E(u_i) u_i^* m)$. For $n\in\N$ it follows $nz=\sum E(n u_i)u_i^*=
\sum E_h(z n u_i)u_i^*=zn$ and therefore $z\in\N'\cap\M$.
}
and the correspondence $z\leftrightarrow E_h^z$ is one-to-one.
If $\N'\cap\M=\M_R\cong \A_L$ then
there exists a uniquely determined $d\in\A_L$ such that
$z=d\lef\one_\M$ and therefore, by \lem{2.16}i) $E_h^z =E_{hd}.$
By \lem{2.22} the map $\A_L\ni d\mapsto hd\in \L(\A)$
is one-to-one and therefore
$\L(\A)\ni l\mapsto E_l \in Hom({_\N}\M_\N\to {_\N}\N_\N)$ is one-to-one.
Now $E_h^z$ is nondegenerate if and only if $z$ and therefore
$d$ are invertible, which by \prop{2.25} is equivalent to
$ l=hd$ being a nondegenerate left integral.
Next, if $E\in Hom({_\N}\M_\N\to {_\N}\N_\N)$ is positive then $E(m)=
\sum_i E_h(z_imz_i^*)$ for some elements $z_i\in\N'\cap\M$
and all $m\in\M$ [Lo2].
By \no{1.32} and \lem{2.16} this implies $E=E_l$ for
$ l=\sum_i ha_ia_{i*}$, where $z_i=a_i\lef\one_\M,\ a_i\in\A_L$
(note that $\A_{L*}= S(\A_L)^*=\A_R)$.
Using \lem{2.15} this implies further
$$
d_R( l)\equiv S( l\1) l\2 =\sum_i S(a_i)S(h\1)h\2 a_{i*}
=\sum_i S(a_i) S(a_i)^* \ge 0
$$
and therefore $ l$ is a positive weight on $\hat\A$ by
\prop{2.25}.
Finally, $E_l(\one_\M)=n_L(l)\lef\one_\M$
and therefore, using $\M_R\cong \A_L,\ E_l$ is normalized if and only if
$ l$ is normalized. \qed

\bsn
The above results will become useful in Sect. 4, where
the conditions of standardness and  $\N'\cap\M=\M_R$ in
\thm{2.26}v) will hold if $\A$ acts regularly on $\M$, see
\defi{4.14} and \thm{4.16}. In particular this holds, if
$\A$ is pure, $\M$ is a factor and $\A$ acts
outerly on $\M$, which are precisely the
conditions of \thm{4.10} guaranteeing that $\N\subset\M\subset\M\cros\A$
provides a Jones triple of factors.

\bsn
Note that for the left integrals \no{2.6.2} of our example we
get
\bleq{2.42}
E_l(m)={1\over |G|}\sum_{g\in G,\,h\in H}c(h)\,\al_g(u(h)m).
\eeq
We conclude this Section with studying the {\em index} of $E_l$.
Picking up the terminology of Watatani [Wa]
we call $ E_l$  of
{\em index-finite type} if there exists a {\em quasi-basis} $u_i,v_i\in \M$
such that for all $m\in\M$
$$\sum_i u_i E_l(v_i m)=m=\sum_i E_l(mu_i)v_i$$
In this case the (Watatani) index of $ E_l$ is defined by
\footnote{This is independent of the
choice of quasi-basis as in [Wa], even
if $ E_l$ is not normalized.}
\bea\lb{3.9}
\Ind E_l =\sum _i u_iv_i \in C(\M)
\eea

The following Lemma says that if $E_l$ is of index-finite
type for some
nondegenerate left integral $l\in\L(\A)$ then it is of
index-finite type for all nondegenerate left integrals in $\A$.

\Lem{2.27}
{Let $l\in\L(\A)$ and let $u_i,v_i\in\M$ be a quasi-basis for
$E_l$. For $d_\s\in\A_\s$ invertible and $\s=L,R$ let
$x_\s:=d_\s^{-1}\lef\one_M$.
Then\\
i) $(u_ix_R,\,v_i)$ and $(u_i,\,v_ix_L)$ are quasi-basises
for $E_{ld_R}$ and $E_{ld_L}$, respectively.
\\
ii) If $d\in\A_L\cap\A_R$ is invertible then
\beq\lb{index}
\Index{E_{ld}} = (d^{-1}\lef\one_\M)\,\Index{E_l}
\eeq
}
\proof
Since $x_\s\in\M_R$ commutes with $\M^\A$, part i) follows
immediately from \lem{2.16}i)+ii) and part ii) follows from
\lem{2.16}iii).
\qed

\bsn
Considering the special case $\M=\hA$ with its canonical left
$\A$-action we define
\Def{2.28}
{For any nondegenerate left integral $l\in\L(\A)$ we define
its index $\Ind l$ to be the Watatani index of
$E_l:\hat\A\to\hat\A_L$
}
To see that this is well defined we have to check that
$E_l:\hA\to\hA_L$ is indeed of index-finite type.

\Proposition{2.29}
{{\rm (Index, quasi-basis and dual left integral)\\}
Let $ l\in \L(\A)$ and $\lambda\in \L(\hat\A)$ be a dual
pair of nondegenerate left integrals and put
$n_R(l)\equiv\hat\e_R\e_R(l)\equiv l\2S^{-1}(l\1)$ as in (\ref{n_s1}).
Denote $E_l:\hA\ni\phi\mapsto (l\arr\phi)\in\hA_L$ and similarly
$\hat E_\l:\A\ni a\mapsto (\l\arr a)\in\A_L$.
Then\\
i) $E_l$ is of index finite type with quasi-basis
$\sum_i u_i\otimes v_i=\l\2\otimes \hat S^{-1}(\l\1)$
\\
ii) $\hat E_\l$ is of index finite type with quasi-basis
$\sum_i u_i\otimes v_i=l\2\otimes S^{-1}( l\1)$
\\
iii)
\bea\lb{3.9c}
\Ind l&=\l\2\hat S^{-1}(\l\1) \equiv& n_R(\l)\in\hat\A_R\cap C(\hat\A)
\\
\Ind\lambda &= l\2 S^{-1}(l\1)\equiv& n_R(l)\in\A_R\cap C(\A).\lb{3.9c'}
\eea
iv) Assume $l$ normalized and $\Ind l$
invertible. If $\Ind l$ is hypercentral, i.e.
$\Ind l\in\hat \Z:=\hat\A_L\cap\hat\A_R\cap C(\hat\A)$, then the normalized
left integral $\dot\l=n_R(\l)^{-1}\l$
satisfies
\bleq{3.9e'}
\hat E_\l = (\Ind \dot\l)\hat E_{\dot\l}
\eeq
In this cases, $\Ind\dot\l\in \Z:=\A_L\cap\A_R\cap C(\A)$
and under the natural isomorphism $\Z\cong\hat\Z$ these indices
coincide, i.e. for $\s=L,R$
\bleq{3.9e}
\Ind\dot\l =\hat\e_\s(\Ind l) \quad ,\quad
\Ind l = \e_\s(\Ind\dot\l)
\eeq
}
{\bf Proof:}
Part (i) is the dual of part (ii). To prove (ii) we compute
$$
\sum_i (\l\arr (au_i))v_i = (\lambda\arr(a l\2))S^{-1}( l\1)
= a\1 S(a\2)a\3 = a
$$
where we have used the identity
$$
\ba{rcl}
l\1\o(\lambda\arr al\2) &=& l\1\o a\1l\2\bra\lambda\mid a\2l\3\ket\\
&=&(\one\o a\1)\Delta(l_R\circ\lambda_L(a\2))\\
&=& S(a\3)\o a\1 S(a\2)
\ea
$$
following from $l_R\circ\lambda_L=S$, see \prop{2.23}.
Similarly,
\beanon
\sum_i u_i (\l\arr (v_i a))&=&
l\3 S^{-1} ( l\2)a\1\bra\lambda\mid S^{-1} ( l\1) a\2\ket
\\
&=&l\3 S^{-1}( l\2)a\1\bra \hat S^{-1} (\lambda_R(a\2))\mid  l\1\ket
\\
&=&  a\3 S^{-1} (a\2)a\1= a
\eeanon
where in the third line we have used $ l_L\circ \hat S^{-1} \circ
\lambda_R=id_\A$ by \no{2.38}. Thus $u_i\otimes v_i$ provides a quasi
basis for $\hat E_\l$. Eqs. \no{3.9} and \no{2.26} then imply
$$
\Ind\l= l\2 S^{-1}(l\1)=\hat\e_R\e_R(l)= n_R(l)
$$
by \no{n_s2}.
Together with its dual version this proves part (iii). To prove
part (iv) let now $l$ be normalized implying
$\Ind\l=\one$ by \no{3.9c}. If $\Ind l\equiv n_R(\l)$ is
invertible then $\l$ is normalizable and
$\dot\l:=n_R(\l)^{-1}\l$ is normalized. If $\Ind l\equiv
n_R(\l)\in\hat\A_L\cap\hat\A_R\cap C(\hat\A)$
then \lem{2.27}ii) and $\Ind\l=\one$ imply
$$
\Ind\dot\l=n_R(\l)\arr\one=\hat\e_\s(\Ind l)\in\A_L\cap\A_R\cap C(\A)
$$
where we have used \lem{2.18}. \Eq{2.19} then also implies
$
\Ind l = (\e_\s\hat\e_{\s'})(\Ind l) =\e_\s(\Ind\dot\l)
$
proving \no{3.9e}.
Finally, \no{3.9e'} follows from
\thm{2.26}ii).
\qed

\bsn
In \thm{3.6} we will generalize \prop{2.29} to the case $\M=\N\cros\hat\A$
with its canonical left $\A$-action, and in \thm{4.3} we will
show that for all $\A$-module von-Neumann algebras
$\M$ and for all
positive normalized and nondegenerate left integrals
$l\in\L(\A)$ the index of $E_l:\M\to\N$ is always bounded by
$\Ind E_l\le\tau_\lef(\Ind l)$.
Also note that in a pure weak Hopf algebra $\A$ the condition
$\Ind l\in\hat \Z\equiv\hat\A_L\cap\hat\A_R\cap C(\A)$
means $\Ind l\in\CC$ by \prop{2.19}ii). On the other hand,
if $\hat\A$ is pure
then $\Ind l\in\CC$ always holds by \prop{2.19}iii).



\sec{Crossed Products}

In this section $\A,\ \hat\A$ and $\M$ will have the same meaning as in the
previous section.
We also continue to denote the units by $\one_\M\in\M,\
\one\in\A$ and $\hat\one\equiv\e\in\hat\A$, respectively.
Using the epimorphism
\beq \lb{3.1}
\mu_\lef :\A_L\ni a\mapsto a\lef\one_\M\in \M_R
\eeq
given in \prop{2.10} we now define the crossed product $\M\cros \A$ to
be the $\CC$-vector space
$$
\M\cros\A=\M\otimes_{\A_L}\A
$$
where $\A_L$ acts on $\A$ by
left multiplication and on $\M$ by right multiplication
via its image $\M_R$ under $\mu_\lef$. Thus $\M\cros\A$ is the linear span of
elements $(m\,\cros\, a),\ m\in\M,\ a\in\A$, modulo the relation
\beq \lb{3.2}
(m\,\cros\, ba)=(m(b\lef\one_\M)\,\cros\, a),\ \forall b\in\A_L
\eeq
In Section 3.1 we prove that equs.
(\ref{1.43})-(\ref{1.44}) in \defi{1.6} indeed provide a $*$-algebra
structure on $\M\cros\A$.
In Section 3.2 we construct the dual $\hat\A$-action on
$\M\cros\A$, having $\M$ as its $\hat\A$-invariant subalgebra.
We then apply our results of Section 2.5 to show that positive
normalized left integrals $\l\in\hat\A$ give rise to
conditional expectations $\hat E_\l:\M\cros\A\to\M$ which are
always of index finite type.
In fact, their quasi-basis can always be given explicitely in
terms of their dual left integrals $l\in\A$.

In Section 3.3 we develop a Plancherel--duality theory for {\em
positive} left integrals fitting with the well
known notion of Haagerup--duality for conditional expectations in Jones
theory. This improves the duality concept of [BNS],
 which did not respect positivity.

In Section 3.4 we pass to a von-Neumann algebraic setting by
associating a so-called regular faithful Hilbert space
representation $\pi_{cros}$ of
$\M\cros\A$ on $\H\o L^2(\A,\hat h)$ with any
representation $\pi$ of $\M$ on $\H$.
This is done by identifying $\M\cros\A$ with a (non-unital)
$*$-subalgebra of $\M\o\End \A$.
In particular, $\M\cros\A$ becomes a $C^*$- or von-Neumann
algebra whenever $\M$ is endowed with such a structure.
This will be the starting point for the relation with Jones
theory in Section 4.

\subsection{The $*$-algebra $\M\cros\A$}

We now prove that the definitions
(\ref{1.43})-(\ref{1.44}) provide a $*$-algebra
structure on $\M\cros\A$.

\Thm{3.1}
{{\rm (Crossed product)\\}
i) The crossed product $\M\cros\A$ becomes a $*$-algebra with
multiplication and $*$-structure given for $m,m'\in\M$ and $a,a'\in\A$ by
\bea
(m\,\cros\, a)(m'\,\cros\, a')&:=& (m(a\1\lef m')\,\cros\, a\2a')\lb{3.3}\\
(m\,\cros\, a)^* &:=& (\one_\M\,\cros\, a^*)(m^*\,\cros\, \one)\lb{3.4}
\eea
ii) For all $m\in\M$ and $a\in\A$  we have
\beq\lb{3.4b}
((a\lef m)\,\cros\, \one) =
(\one_\M\,\cros\, a\1)\,(m\,\cros\, \one)\,(\one_\M\,\cros\, S(a\2))
\eeq
iii) The unit in $\M\cros\A$ is given by $(\one_\M\,\cros\, \one)$.
Moreover, we have $\M\cong (\M\,\cros\, \one)\equiv
\M\otimes_{\A_L} \A_L$ and
$\A/I_\lef \cong (\one_\M\,\cros\, \A)\equiv\M_R\o_{\A_L}\A$
as unital *-subalgebras of $\M\cros\A$, where $I_\lef\subset\A$
is the two-sided ideal generated by $\Ker \mu_\lef$.
\\
iv) Putting $\N\equiv\M^\A$ the inclusions \no{1.44d} - \no{1.44f}
hold, i.e.
\bea
\one_\M\,\cros\, \A &\subset& \N'\cap(\M\cros\A)
\lb{3.5a}\\
\one_\M\,\cros\, \A_R &\subset& \M'\cap(\M\cros\A)
\lb{3.5b}\\
\one_\M\,\cros\, (\A_R\cap C(\A)) &\subset& C(\M\cros\A)
\lb{3.5c}\\
\one_\M\,\cros\, (\A_L\cap\A_R\cap C(\A)) &\subset& \M\cap C(\M\cros\A)
\lb{3.5d}
\eea
v)\hspace{8em} $\M\cap C(\M\cros\A) = \N\cap C(\M) = \N\cap C(\M\cros\A)$
}
To prepare the proof of \thm{3.1} we need the following 2 Lemmata.

\Lem{3.2}{{\rm [N2]}
For $a\in\A$ and $b\in\A_L$ we have
\bea
ab &=& \hat\e_L\e_R(a\1 b)a\2\lb{3.5} \\
ba &=& a\2\hat\e_L \e_L (ba\1)\lb{3.6}
\eea
}
\Lem{3.3}
{The two-sided ideal $I_\lef$ generated by $\Ker\mu_\lef$ satisfies
$I_\lef\subset\K_\lef$ and
$$
I_\lef =\A\,(\Ker\mu_\lef) =(\Ker\mu_\lef)\,\A=I_\lef^*.
$$
}
{\bf Proof:}
By \cor{2.17}iii) $I_\lef\subset\K_\lef$ and by \prop{2.10}
$\Ker\mu_\re=z_\re\A_L$ for a central projection
$z_\re\in\A_L\cap C(\A)$. Hence $I_\re=z_\re\A$, proving the
assertions.
\qed

\bsn
{\bf Proof of \thm{3.1}}\\
To show that the multiplication (\ref{3.3}) is well defined let
$b\in\A_L$. Then
\beanon
(m\,\cros\, ba)(m'\,\cros\, a') &=& (m(b\1a\1\lef m')\,\cros\, b\2a\2a')\\
&=& (m(ba\1\lef m')\,\cros\, a\2a')\\
&=& (m(b\lef\one_\M)(a\1\lef m')\,\cros\, a\2a')\\
&=& (m(b\lef\one_\M)\,\cros\, a)(m'\,\cros\, a')
\eeanon
where in the second line we have used \lem{2.15}i) and in the third line
\lem{2.16}i). Next,
\beanon
(m\,\cros\, a)(m'\,\cros\, ba') &=& (m(a\1\lef m')\,\cros\, a\2 ba')\\
&=& (m(a\1\lef m')\,\cros\, \hat\e_L\e_R(a\2b) a\3a')\\
&=& (m(a\1\lef m')(a\2b\lef\one_\M)\,\cros\, a\3a')\\
&=& (m\,\cros\, a)(m'(b\lef\one_\M)\,\cros\, a')
\eeanon
where in the second line we have used \lem{3.2} and in the third line
(\ref{3.2}) and \no{1.33}. Thus the multiplication
(\ref{3.3}) is well defined on $\M\cros\A$. The associativity of the product
(\ref{3.3}) is straightforward and is left to the reader. Next, we compute
\beanon
(m\,\cros\, \one)(m'\,\cros\, a) &=& (m(\one\1\lef m')\,\cros\, \one\2a)\\
&=&(m(\one\1\lef m')(\one\2\lef\one_\M)\,\cros\, a)\\
&=& (mm'\,\cros\, a)
\eeanon
since $\Delta(\one)\in\A_R\otimes\A_L$ by \no{2.37'}, and
\beanon
(m\,\cros\, a)(\one_\M\,\cros\, a') &=& (m(a\1\lef \one_\M)\,\cros\, a\2a')\\
&=& (m(a\1S(a\2)\lef\one_\M)\,\cros\, a\3a')\\
&=& (m\,\cros\, a\1S(a\2)a\3a')\\
&=& (m\,\cros\, aa')
\eeanon
where we have used (\ref{1.33}) in the second line and (\ref{1.1}) in
the last line.
Thus $(\one_\M\,\cros\, \one)$
is the unit in $\M\cros \A$ and $m\mapsto (m\,\cros\, \one)$ and
$a+I_\lef \mapsto (\one_\M\,\cros\, a)$
define unital inclusions $\M \subset\M\cros \A$
and $\A/I_\lef\subset\M\cros\A$, respectively.
\footnote{Note that by definition
$(\one_\M\,\cros\, a)=0$ iff $a\in(\Ker \mu_\lef)\A\equiv I_\lef$.}

Next, to show that (\ref{3.4})
provides a well defined $*$-structure let again
$b\in\A_L$. Then
\beanon
(m\,\cros\, ba)^* &=& (\one_\M\,\cros\, a^*b^*)(m^*\,\cros\, \one)\\
&=& (\one_\M\,\cros\, a^*)((b\1^*\lef m^*)\,\cros\, b\2^* )\\
&=& (\one_\M\,\cros\, a^*)((b\1^* \lef m^*)(b\2^*\lef\one_\M)\,\cros\,\one)\\
&=& (\one_\M\,\cros\, a^*)((b^*\lef m^*)\,\cros\, \one)\\
&=& (\one_\M\,\cros\, a^*)((S^{-1}(b)\lef m)^*\,\cros\, \one)\\
&=& ((S^{-1} (b)\lef m)\,\cros\, a)^*\\
&=& (m(b\lef \one_\M)\,\cros\, a)^*
\eeanon
where in the third line we have used $\Delta(b^*)\in\A \otimes \A_L$, in the
fifth line (\ref{1.32}) and in the last line \lem{2.16}ii)
(note that $S^{-1}(b)\in\A_R$).  Thus $(m\,\cros\, a)^*$ is well
defined. Moreover, it is
involutive since
\beanon
(m\,\cros\, a)^{**} &=& [(\one_\M\,\cros\, a^*)(m^*\,\cros\, \one)]^*\\
&=& ((a\1^*\lef m^*)\,\cros\, a\2^*)^*\\
&=& (\one_\M\,\cros\, a\2)((S^{-1}(a\1)\lef m)\,\cros\, \one)\\
&=& ((a\2 S^{-1}(a\1)\lef m)\,\cros\, a\3)\\
&=& ((\one\1\lef m)\,\cros\, \one\2 a)\\
&=& (m\,\cros\, a)
\eeanon
where in the third line we have used (\ref{1.32}) and in the fourth line
(\ref{1.8}). Finally, to prove the antimultiplicativity
of the $*$-operation it is enough
to check for all $a\in\A,\ m\in\M $
\beanon
(m\,\cros\, \one)^*(\one_\M\,\cros\, a)^*
&=& (m^*\,\cros\, \one)(\one_\M\,\cros\, a^*)\\
&=& (m^*\,\cros\, a^*)\\
&=& (m^*\,\cros\, a^*)^{**}\\
&=& [(\one_\M\,\cros\, a)(m\,\cros\, \one)]^*
\eeanon
Thus $\M\cros \A$ is a $*$-algebra containing
$\M\cong (\M\,\cros\, \one)$ and
$\A/I_\lef\cong(\one_\M\,\cros\, \A)$ as *-subalgebras, which proves i)
and iii).
To prove part ii) we use \no{3.3} to calculate
\beanon
(\one_\M\,\cros\, a\1)\,(m\,\cros\, \one_\A)\,(\one_\M\,\cros\, S(a\2))
&=&((a\1\lef m)\,\cros\, a\2S(a\3))\\
&=&((a\1\lef m)\,(a\2S(a\3))\lef\one_\M)\,\cros\, \one_\A)\\
&=&((a\1\lef m)\,(a\2\lef\one_\M)\,\cros\, \one_\A)\\
&=&((a\lef m)\,\cros\, \one_\A)
\eeanon
where we have used \no{1.33} and $a\2 S(a\3)\in\A_L$.
To prove the statements in part iv) first note that \no{3.5a}
follows by putting $m=\one_\M$ in Proposition
2.11iii) and using $((a\1\lef\one_\M)\,\cros\, a\2)=(\one_\M\,\cros\, a)$.
To prove \no{3.5b}
let $a\in\A_R$ and use \lem{2.15}ii) to
compute
$$
(\one_\M\,\cros\, a)(m\,\cros\, \one)=((a\1\lef m)\,\cros\, a\2)
=((\one\1\lef m)\,\cros\, \one\2)(\one_\M\,\cros\, a)
=(m\,\cros\, \one)(\one_\M\,\cros\, a)\ .
$$
The inclusion \no{3.5c} follows, since
$\M'\cap(\one_\M\,\cros\,\A)'\cap\M\cros\A=C(\M\cros\A)$ and
the inclusion \no{3.5d} follows, since
$(\one_\M\,\cros\,\A_L)\subset(\M\,\cros\,\one)$.
Finally, to prove part v) let
$(n\,\cros\,\one)\in\M\cap C(\M\cros\A)$.
Then $n\in C(\M)$ and by \no{3.4b}
\beanon
\left( (a\lef(nm))\,\cros\,\one\right) &=&
(\one_\M\,\cros\,a\1)\,(nm\,\cros\,S(a\2))
\\
&=& (n\,\cros\,\one)\,(\one_\M\,\cros\,a\1)\,(m\,\cros\,S(a\2))
\\
&=& \left(n(a\lef m)\,\cros\,\one\right)
\eeanon
for all $a\in\A$ and all $m\in\M$. Hence $n\in\N$ by
\prop{2.11}iii) implying $\M\cap C(\M\cros\A)\subset\N\cap C(\M)$.
Now by \no{3.5a} $\N\cap C(\M)\subset\N\cap C(\M\cros\A)$ and
the inclusion $\N\cap C(\M\cros\A)\subset\M\cap C(\M\cros\A)$
holds trivially. This proves part v) and therefore concludes
the proof of \thm{3.1}.
\qed

\bsn
Applying the above crossed product construction to the example
of the partly inner group action in Section 2.5 one immediately
verifies that
\beq\lb{3.7a}
\M\>cros \A\ni (m\,\>cros (h,g))
\mapsto (mu(h)\,\>cros_\alpha\, g)\in\M\>cros_\alpha \,G
\eeq
provides an identification $\M\>cros \A=\M\>cros_\alpha G$.
Thus, if the action is standard (i.e. if the implementers
$u(h),\ h\in H$, are linearly independent in $\M$), then Eqs.
\no{3.5a}-\no{3.5d} confirm with \no{0.7a}-\no{0.7c}.

\bsn
In Section 4.2 we will show that for a
factor $\M$ we have $\M'\cap(\M\cros\A)=(\one_\M\,\cros\, \A_R)$ if and
only if the $\A$-action on $\M$ is outer in a suitable sense.
If in addition the action is standard (i.e. $\Ker\mu_\lef=0$)
we will also have
$\N'\cap(\M\cros\A)=(\one_\M\,\cros\, \A)$,\
$\N'\cap\M=\M_R\cong\A_L$ and
$C(\M\cros\A)=\one_M\,\cros\,(\A_R\cap C(\A))$, thus verifying
the scenario of Section 1.1.
A generalization of these identities to non-factorial
von-Neumann algebras will be given in Section 4.4.

We also remark that for non-standard actions one might try to
pass to a quotient weak Hopf algebra structure on
$\B:=\A/\A_\lef,\ \A_\lef:=I_\re+S(I_\re)\subset\K_\re$,
which would act standardly on $\M$.
In fact, one may show that $\B$ and
$\hat\B:=\A_\re^\perp\subset\hA$ again give a dual pair of weak
$*$-bialgebras. But since $\e(z_\re)>0$ (unless $z_\re=0$ and
$\B=\A$) by the faithfulness of $\e|_{\A_L}$, one has
$\one_{\hat\B}\neq\one_{\hA}$ and therefore we are not sure
whether $\B$ and $\hat\B$ still satisfy all weak Hopf axioms.

\subsection{Dual actions and conditional expectations}

In this subsection we show that as for ordinary crossed products the
algebra $\M\cros\A$ allows for a natural dual $\hat\A$-left action given by
\beq\lb{3.8a}
\phi\lef (m\,\cros\, a)=(m\,\cros\, (\phi\arr a))
\eeq
where $\phi \in\hat\A$ and $(m\,\cros\, a)\in\M\cros\A$. To show that this is well
defined, we consider equivalently the associated dual coaction.

\Prop{3.4}
{{\rm (Dual coaction)}\\
Let $\rho:\M\cros\A\to (\M\cros \A)\otimes\A$ be given by
$$ \rho(m\,\cros\, a):= (m\,\cros\, a\1)\otimes a\2$$
Then $\rho$ provides a well defined right coaction of $\A$ on $\M\cros \A$
such that $(\M\cros\A)_R=(\one_\M\,\cros\, \A_R)$.
}
{\bf Proof:}
If $b\in \A_L$ then by \lem{2.15}i)
$\rho(m\,\cros\, ba)=\rho(m(b\lef\one_\M)\,\cros\, a)$
showing that $\rho$ is well defined.
The axioms \no{2.1}-\no{2.3} and the
identity $(\M\cros\A)_R=(\one_\M\,\cros\, \A_R)$ follow immediately.
Using \no{3.3} and \no{3.4} one also immediately checks that $\rho$ provides a
*-algebra homomorphism.
\qed

\Cor{3.5}
{The $\hat\A$-invariant subalgebra of $\M\cros\A$ is given by
$\M\equiv\M\cros\A_L$.
}
{\bf Proof:} By \lem{2.12} $(m\,\cros\, a)\in\M\cros\A$
is $\hat\A$-invariant if and only if
$\Delta(a) \in\A\otimes\A_L$ and therefore $a\equiv\e(a\1)a\2\in\A_L$. \qed

\bsn
Similarly as in \thm{2.26} we now obtain for every left integral
$\lambda\in\hat\A$ a
$\M$-$\M$ bimodule map $\hat E_\lambda:\M\cros\A\to \M$ by putting
\beq\lb{3.8}
\hat E_\lambda(m\,\cros\, a):=(m\,\cros\, (\lambda\arr a))
\equiv (m\mu_\lef(\lambda\arr a)\,\cros\, \one)
\eeq
where the second equality follows from $\lambda\arr a\in \A_L$, which indeed
implies $\hat E_\l(\M\cros \A)\subset\M\cros\A_L\cong \M$.

In our example of the partly inner group action $(\al,\,u)$ of
Sect. 2.5 we have $\M\cros\A=\M\cros_\al\,G$ by \no{3.7a}
and for a left integral $\l\in\L(\hA)$ of the form \no{2.6.3}
we get
\bleq{3.7b}
\hat E_\l(m\cros_\al\,g)=
{1\over|H|}\sum_{h\in H}\hat c(h)\delta(hg)\,mu(h^{-1})
\eeq
In general $\hat E_\l$ is normalized (positive, nondegenerate)
if and only if
its restriction to $(\one_\M\,\cros\, \A)$ is normalized
(positive, nondegenerate).
If $(\one_\M\,\cros\, \A)\cong\A$ (i.e. $\Ker\mu_\lef=0)$, then applying
\thm{2.26}v) to the $\hat\A$-action on $\A$ this is further
equivalent to $\lambda$ itself being normalized (positive,
nondegenerate).
In fact, most of the following results already follow from their validity in
the case $\M\cros \A=\A$ (i.e. $\M=\M_R=\hat \A_R)$.
The following is a straightforward generalization from
\prop{2.29}.

\Thm{3.6}
{{\rm (Index, quasi-basis and dual left integral)\\}
Let $ l\in \L(\A)$ and $\lambda\in \L(\hat\A)$ be a dual
pair of nondegenerate left integrals and put
$n_R(l)\equiv\hat\e_R\e_R(l)\equiv l\2S^{-1}(l\1)$ as in (\ref{n_s1}).
For any crossed product $\M\cros\A$ let $\hat
E_\l:\M\cros\A\to\M$ be given by \no{3.8}.
Then $\hat E_\l$ is of index finite type with quasi-basis
$\sum_i u_i\otimes v_i=
(\one_\M\,\cros\,  l\2)\otimes (\one_\M\,\cros\, S^{-1}( l\1))$
and we have
\bea
\hat E_\l (\one_\M\,\cros\,  l)&=&\one_\M
\lb{3.9a}\\
\hat E_\l(\one_\M\cros\one_\A) &=& \tau_\re(n_R(l))\equiv
 \tau_\lef(\Ind l) \in C(\M)\cap\M_R
\lb{3.9d}\\
\Ind\hat E_\l &=&
(\one_\M\,\cros\,  \Ind\lambda)\in (\M\cros\A)_R\cap C(\M\cros\A)\
\lb{3.9b}
\eea
Moreover, under the conditions of \prop{2.29}iv) we have
$\hat E_\l=\tau_\lef(\Ind l)\, \hat E_{\dot\l}$ and
\bleq{3.9f}
\tau_\lef(\Ind l)\equiv(\tau_\lef(\Ind l)\cros\one_\A)
=(\one_M\cros\Ind\dot\l)\in\N\cap C(\M\cros\A).
\eeq
}
{\bf Proof:}
By definition
$\hat E_\l(\one_\M\,\cros\,  l)=(\one_\M\,\cros\,(\l\arr l))=
(\one_\M\,\cros\, \one)$, which proves \no{3.9a}.
\Eq{3.9d} follows, since by \no{3.9c} and the definition
\no{2.16} of $\tau_\re$
$$
\tau_\lef(\Ind l) =\tau_\lef(\hat\e_R\e_R(\l))=
(\lambda\arr\one)\lef\one_\M = \hat E_\l(\one_\M\,\cros\, \one)
$$
where in the last equation we have identified $\M_R$ with
$(\one_\M\,\cros\,\A_L)\subset\M\cros\A$.
Since by the $\M$-$\M$--bimodule property (or by
\lem{2.16} iii) ) $\hat E_\l(\one_\M\,\cros\, \one)\in C(\M)$, this proves
\no{3.9d}.
The fact that
$u_i\otimes v_i =
(\one_\M\,\cros\,  l\2)\otimes (\one_\M\,\cros\, S^{-1}( l\1))$
 is a quasi-basis for $\hat E_\l$ follows
similarly as in the proof of \prop{2.29}ii). Indeed by the same
arguments
$$
\sum_i \hat E_\l((m\,\cros\, a)u_i)v_i
= (m\,\cros\, \lambda\arr(a l\2))S^{-1}( l\1)=a
$$
and
\beanon
\sum_i u_i \hat E_\l(v_i (m\,\cros\, a))&=&
\sum_i u_i(\lambda\1\lef v_i)(\lambda\2\lef(m\,\cros\, a))
\\
&=& (\one_\M\,\cros\,  l\3 S^{-1} ( l\2))
(m\,\cros\, a\1)\bra\lambda\mid S^{-1} ( l\1) a\2\ket\\
&=& (m\,\cros\,  l\3 S^{-1}( l\2)a\1)\bra \hat S^{-1} (\lambda_R(a\2))
\mid  l\1\ket\\
&=& (m\,\cros\, a\3 S^{-1} (a\2)a\1)\\
&=& (m\,\cros\, a)
\eeanon
where in the third line we have used $ l\3 S^{-1}( l\2)\in\A_R$, which
commutes with $\M$, and in the fourth line $ l_L\circ \hat S^{-1} \circ
\lambda_R=id_\A$ by \no{2.38}. Thus $u_i\otimes v_i$ provides a
quasi-basis and
$
\Ind\hat E_\l=(\one_\M\,\cros\, l\2 S^{-1}(l\1))=
(\one_\M\,\cros\, \Ind\l)\in\one_\M\cros (\A_R\cap C(\A))
\subset(\M\cros\A)_R\cap C(\M\cros\A)
$
by \no{3.9c'}, \no{3.5c} and \prop{3.4},
which proves \no{3.9b}.
Finally, if $l$ is normalized and $\Ind l\in\Z $ is invertibel, then
by \no{3.9e}, \no{2.16} and \cor{2.17}i)+ii)
$\tau_\lef(\Ind l)=\mu_\lef(\Ind\dot\l)\in\N\cap C(\M)\equiv
\N\cap C(\M\cros\A) $
and \lem{2.27} implies
$\hat E_\l=\tau_\lef(\Ind l)\, \hat E_{\dot\l}$.
\qed

\bsn
As an application of these methods we now show that the
relative commutant $\N'\cap\M\cros\A$ saturates the lower bound
\no{3.5a} if and only if $\N'\cap\M$ saturates the lower bound
\no{2.33a}.
Note, by the way, that in our example of Sect 2.5 these lower bounds are
always saturated, see also \no{0.7a} and \no{0.7b}.

\Prop{3.7}
{Let $\N\equiv\M^\A\subset\M$ be the $\A$-invariant
subalgebra. Then $\A\lef(\N'\cap\M)\subset\N'\cap\M$ and
$\N'\cap(\M\cros\A)=(\N'\cap\M)\cros\A$. In particular,
$\N'\cap\M=\M_R$ if and only if
$\N'\cap(\M\cros\A) = (\one_\M\cros\A)$.
}
\proof
Let $n\in\N,\ m\in\N'\cap\M$ and $a\in\A$. Then
$n(a\lef m)=a\lef(nm)=a\lef(mn)=(a\lef m)n$ implying
$\A\lef(\N'\cap\M)\subset\N'\cap\M$ and
by \no{3.5a} $(\N'\cap\M)\cros\A\subset\N'\cap(\M\cros\A)$.
To prove the inverse conclusion observe that
$\hat E_\l(\N'\cap(\M\cros\A))=\N'\cap\M$ for any normalized
nondegenerate left integral $\l\in\L(\hat\A)$.
Now by \thm{3.6}i) and \Eq{3.5a} a quasi-basis $u_i,v_i$
for $\hat E_\l$ may be chosen in $(\one_\M\,\cros\,\A)$, i.e. in
the commutant of $\N$, and
therefore we get for any $x\in\N'\cap(\M\cros\A)$
$$
x=\sum\hat E_\l(xu_i)v_i\in(\N'\cap\M)\cros\A
$$
Hence $(\N'\cap\M)\cros\A=\N'\cap(\M\cros\A)$. If
$(\N'\cap\M)=\M_R$ this implies
$\N'\cap(\M\cros\A)=(\one_\M\,\cros\,\A)$. Conversely, if
$\N'\cap(\M\cros\A)=(\one_\M\,\cros\,\A)$ then
$\N'\cap\M=\hat E_\l(\one_\M\,\cros\,\A)=(\one_\M\,\cros\,\A_L)=\M_R$.
\qed

\bsn
We conclude this subsection
with using similar arguments to characterize
$\hat E_\l$-invariant (and, more generally, $\hat\A$-invariant) ideals
in $\M\cros\A$ in terms of their intersection with
$\M\equiv \hat E_\l(\M\cros\A)$, thus generalizing standard
results from ordinary crossed product theory.
In particular, this will imply that $\hat\A$-covariant
representations of $\M\cros\A$ are faithful whenever their
restrictions to $\M$ are faithful.

\Lem{3.8}
{{\rm ($\hat E_\l$-invariant ideals)}\\
Let $\lambda\in \L(\hat\A)$ be a nondegenerate left integral
and let $\I \subset\M\cros \A$ be an ideal such that
$\hat E_\l (\I)\subset \I$.
Then $\A\lef(\I \cap\M)=\I\cap \M$ and
$\I =(\I\cap \M)\cros\A$.
}
{\bf Proof:} If $m\in \I\cap\M$ then
$$
((a\lef m)\,\cros\, \one)=
(\one_\M\,\cros\, a\1) (m\,\cros\, \one)(\one_\M\,\cros\, S(a\2))\in\I \cap \M
$$
Since $\I$ is an ideal we have $(\I\cap \M)\cros\A\subset \I$.
Conversely, let
$x\in\I$. Then by \thm{3.6}i)
$$
x=\hat E_\l (x(\one_\M\,\cros\,  l\2))(\one_\M\,\cros\, S^{-1} ( l\1))
\in \hat E_\l (\I)(\one_\M\,\cros\, \A)
\subset(\I\cap \M)\cros\A\quad\qed
$$
\Cor{3.9}
{Let $f:\M\cros\A\to\B$ be a homomorphism of algebras and
let $\lambda \in \L(\hat \A)$ be nondegenerate.
Assume there exists a linear map
$F_\l:\B\to B$ such that
\beq\lb{3.10}
F_\l\circ f=f\circ \hat E_\l
\eeq
Then $f$ is injective iff its restriction to $\M$ is injective.
}
{\bf Proof:} Let $\I=\Ker f$ then $\hat E_\l(\I)\subset\I$
by \no{3.10}. Hence
$\I\equiv(\I\cap\M)\cros\A=0$ if $\Ker f\cap \M=0$. \qed



\subsection{Plancherel-Duality}

In this subsection we develop a modified duality concept for {\em
positive} nondegenerate left integrals $l\in\L(\A)$ and
$\l\in\L(\hA)$, such that for sufficiently regular $\A$-module
algebras $\M$ the conditional expectations $E_l:\M\to\M^\A$ and $\hat
E_\l:\M\cros\A\to\M$ become dual in the sense of Haagerup [Ha], 
see also [Di, La, St].
Simultaneously, with any positive nondegenerate and normalized left
integral $l\in\A$ we will associate a ``Jones projection'' $e_l\in\A$ such
that we have the Jones relation
\beq\lb{3.3.1}
e_lme_l = E_l(m)e_l = e_lE_l(m),\quad\forall m\in\M
\eeq
as an identity in $\M\cros\A$.

\bsn
Note that our purely algebraic duality notion for left integrals in
\prop{2.23} did not touch the question of positivity.
In fact, by \prop{2.25}i)+ii)
the dual
$\lambda\in\L(\hA)$ of a given left integral $l\in\L(A)$
 is positive if and only if $d_L(l)>0$, whereas $l$ itself
is positive if and only if $d_R(l)>0$.
We also warn the reader that the positivity of $l$ as a left integral
(i.e. considered as a functional on $\hA$) in general does not imply
$l$ to be positive
(not even selfadjoint) as an element in the
$C^*$-algebra $\A$. In particular, in \no{3.3.1} we cannot choose
$l=e_l$ as in ordinary finite dimensional $C^*$-Hopf algebras.
Instead we now define, for every positive left integral
$l=hd_R(l)\in\L(\A)$,
\bleq{3.3.1'}
e_l:=d_R(l)^{1/2}\,h\,d_R(l)^{1/2}
\eeq

\Prop{3.15'}
{Let $l\in\L(\A)$ be a positive left integral and put $E_l:\M\ni
m\mapsto l\lef m\in\N\equiv\M^\A$.
For $m\in\M$ consider $m\equiv(m\,\cros\, \one)$ and
$e_l\equiv(\one_\M\,\cros\, e_l)$ as elements in $\M\cros\A$.
Then
\\
i) $e_l$ is positive, commutes with $\M^\A$ and satisfies \no{3.3.1}.
\\
ii) $e_l=e_{l'}\Rightarrow l=l'$
\\
iii) $e_l^2=n_\s(l)e_l$. In particular, if $l$ is normalized then $e_l$ is
a projection.
\\
iv) $l$ is nondegenerate if and only if $e_l$ is nondegenerate (as
functionals on $\hA$). In this case $l$ is normalized if and only if $e_l$ is
a projection.
}
\proof
i) $e_l$ is positive since $h=h^*h$ and
 $e_l\in (\one_\M\,\cros\, \A)$ commutes with
$\M^\A$ by \no{3.5a}.
To prove \no{3.3.1}
we use that for any left integral $l\in \L(\A)$ we have
\bea \lb{4.2}
aml &\equiv& (a\1\lef m) a\2 l
=(a\1\lef m)(a\2\lef\one_\M) l\nonumber\\
&=& (a\lef m) l,\ \forall a\in\A,m\in\M
\eea
by \no{1.19}, \no{1.33} and \no{3.2}.
Using that $\A_R\subset\M\cros\A$ commutes with
$\M\subset\M\cros\A$ by \no{3.5b} this gives
\beanon
e_lme_l &=& d_R(l)^{1/2}  lmhd_R(l)^{1/2}
\\
&=& d_R(l)^{1/2}(l\lef m) hd_R(l)^{1/2}\\
&=& E_l(m)e_l = e_lE_l(m)
\eeanon
where the last identity holds since $E_l(m)\equiv\l\lef m\in\M^\A$.
Thus we have shown i). To prove ii)
let $\r_h\equiv S(\l_h)\in\R(\hA)$ be the right integral dual
to $h$,
i.e. the unique solution of $h\arl\r_h=\one$. Since similarly as in
\cor{2.13} $\A_R\subset\A$ is the fixed point algebra under the
canonical right $\hA$-action we conclude $e_l\arl\r_h=d_R(l)$. Since
$l=l' \Leftrightarrow d_R(l)=d_R(l')$, this proves part (ii).
To prove part (iii) we use $lh=n_\s(l)h$ by \no{n_s1}, \no{1.19} and
\no{1.20} to conclude
$$
e_l^2=d_R(l)^{1/2}lhd_R(l)^{1/2}=n_\s(l)e_l.
$$
Hence, by the argument in the proof of (ii), $e_l^2=e_l$ if and only if
$d_R(l)=n_\s(l)d_R(l)$, which in particular holds for $n_\s(l)=\one$,
i.e. if $l$ is normalized. If $l$ is nondegenerate then $d_R(l)$ is
invertible and $e_l^2=e_l\Leftrightarrow n_\s(l)=\one$. Finally, in
this case we also have
$$
e_l=d_R(l)^{1/2}ld_R(l)^{-1/2}
$$
and therefore
$$
e_l\arl\phi =d_R(l)^{1/2}(l\arl\phi)d_R(l)^{-1/2},\quad\forall\phi\in\hA
$$
proving that $e_l$ is nondegenerate. Converseley, if $e_l$ is
nondegenerate, we use the identity
$$
e_l\arl\phi =d_R(l)^{1/2}(h\arl\phi)d_R(l)^{1/2},\quad\forall\phi\in\hA
$$
to conclude from the nondegeneracy of $h$
$$
d_R(l)^{1/2}ad_R(l)^{1/2}=0\Rightarrow a=0,\quad\forall a\in\A.
$$
Clearly, this implies $d_R(l)$ to be invertible and therefore $l$ to
be nondegenerate. This proves part (iv).
\qed

\bsn
Let now $l\in\L(\A)$ be a positive and nondegenerate left integral.
Inspired by the notion of Haagerup duality for conditional
expectations we seek for a left integral $\l\in\L(\hA)$ such that, for
any $\A$-module algebra $\M$, $\hat E_\l:\M\cros\A\to\M$
satisfies
\bleq{3.3.5'}
\hat E_\l(e_l)=\one_\M.
\eeq
By \prop{3.15'}iv) $e_l$ is nondegenerate and therefore $\l\in\hA$ must be
the unique element satisfying
\bleq{3.3.6''}
\l\arr e_l=\one_\A
\eeq
\Prop{3.16'}
{Let $l\in\L(\A)$ be a positive left integral and assume there exists
$\l\in\hA$ satisfying \no{3.3.6''}. Then $l$ is nondegenerate and $\l$
is uniquely determined. Moreover, $\l$ is a positive and nondegenerate
left integral in $\hA$.
}
\proof
By \lem{2.20}ii) and ii') we have
$
e_l=aha^*
$
where $a=S(d_R(l)^{1/2})\in\A_L$. Hence
$$
\one=\l\arr e_l=a(\l\arr h)a^*
$$
and therefore $a$ is invertible. Thus, $d_R(l)$ is invertible, $l$ is
nondegenerate and $\l$ is unique.
Let now $l':=ha^*a\in\L(\A)$. Then
\bleq{3.3a}
\l\arr l'=\l\arr(a^{-1}e_la)=a^{-1}(\l\arr e_l)a=\one.
\eeq
Thus, $\l$ is the left integral dual to $l'$ and therefore
nondegenerate by \prop{2.23}.
By \prop{2.25}ii) $\l$ is also positive, since $d_L(l')=a^*a>0$.
\qed

\bsn
The above results motivate the following
\Def{3.17}
{For any positive and nondegenerate left integral $l\in\L(\A)$ we
 define the {\em p-dual} (p $\equiv$ ``positive'' or ``Plancherel'')
$\lambda_l\in\L(\hat\A)$ to be the unique positve and nondegenerate left
integral satisfying $\l_l\arr e_l=\one$.
}
We now show that \defi{3.17}
indeed provides a sensible
notion of duality, i.e. converseley $l$ is also the
p-dual of $\lambda_l$.
Moreover, upon iterating our constructions we obtain
a generalized  (i.e. in general with algebra-valued index)
Temperley-Lieb-Jones (TLJ) algebra in $(\M\cros\A)\cros\hat\A$.

\Thm{3.18}{{\rm (P-duality and TLJ-algebra)\\}
Let $l\in\L(\A)$ be a positive nondegenerate left integral
and let $\lambda\in\L(\hat\A)$ its p-dual.
Then
\\
i)
$l$ is also the p-dual of $\lambda$.
\\
ii) The indices $\Ind l$ and $\Ind\lambda$
are both positive and invertible and satisfy
\bleq{3.3b}
\Ind l = n_R(\lambda)\in C(\hat\A)\cap\hat\A_R\quad,\quad
\Ind\lambda = n_R(l)\in C(\A)\cap\A_R
\eeq
iii) For any $\A$-module algebra $\M$ the elements
 $e=((\one_\M\,\cros\, e_l)\,\cros\, \hat\one)\in(\M\cros\A)\cros\hat\A$ and
$\hat e=((\one_\M\,\cros\, \one)\,\cros\, e_\lambda)
\in(\M\cros\A)\cros\hat\A$ satisfy the TLJ-relations
\bea
e^2 &=& e\,\Ind\l\lb{3.3.11}\\
\hat e^2 &=&\hat e\,\Ind l\lb{3.3.12}\\
\hat e\, e\, \hat e &=& \hat e
\lb{3.3.13}\\
e\,\hat e\, e &=& e
\lb{3.3.14}
\eea
}
\proof
By the proof of \prop{3.16'} $\l$ is the p-dual of $l$ if and only if
$\l\arr l'=\one$, where $l'=ha^*a,\ a:=S(d_R(l)^{1/2})$.
Using \no{2.45} we conclude
\bleq{3.3c}
d_R(\l)=\e_R(b_l^*b_l)^{-1}
\eeq
where
\bleq{3.3d}
b_l=ag_L\equiv S(d_R(l)^{1/2})g_L.
\eeq
By \lem{3.19} below, $b_l=b_l^*>0$ and therefore the p-dual $\l$ of $l$ is
uniquely determined by
\bleq{3.3e}
d_R(\l)^{-1/2}=\e_R(b_l)\equiv\e_R(S(d_R(l)^{1/2})g_L).
\eeq
Using $\hat S\circ\e_R=\e_L\circ S^{-1}$ and
$g_R=S^{\pm 1}(g_L)=\hat\e_R(\hat g_L)$
this implies by \prop{2.9}
$$
\hat\e_R(\hat S(d_R(\l)^{1/2})\hat g_L) =
\hat\e_R\e_L(d_R(l)^{-1/2}g_R^{-1})g_R
=d_R(l)^{-1/2}.
$$
Since this relation is precisely the dual version of \no{3.3e}, $l$ is
also the p-dual of $\l$.
To prove part (ii) we note $\Ind\l=n_R(l')$ by \no{3.3a} and
\prop{2.29}iii). Hence,
$$
l'\equiv ha^*a=hS(a^*)a=ld_R(l)^{-1/2}S(d_R(l)^{1/2})
$$
by \lem{2.20}ii').
\lem{2.15} now gives
$$
\Delta(l') = \Delta(l)[S(d_R(l)^{1/2})\o d_R(l)^{-1/2}]
$$
and therefore
\bleq{3.3f}
n_R(l')\equiv l'\2S^{-1}(l'\1) = l\2S^{-1}(l\1) \equiv n_R(l).
\eeq
This proves $\Ind\l=n_R(l)$ and by duality also $\Ind l=n_R(\l)$.
By \prop{2.25}iv) these indices are positive and invertible, thus
proving part (ii).
To prove part (iii) first note that \no{3.3.11} and \no{3.3.12}
follow from \prop{3.15'}iii). Next, the analogue of \no{3.3.1} gives
$$
\hat e e \hat e =\hat E_\l(e)\hat e =\hat e
$$
since $\hat E_\l(e)=\one_\M$ by definition of $\l$.
Finally, we prove \no{3.3.14} as an identity in
$\A\cros\hat\A$.
First we note that for any right integral $r\in\R(\A)$ we have
in $\A\cros\hat\A$
\bea
(r\,\cros\, \phi)(a\,\cros\, \psi)
&=& (ra\1\,\cros\, (\phi\arl a\2)\psi)
\nonumber\\
&=&(r\hat\e_R(\e_L(a\1))\,\cros\, (\phi\arl a\2)\psi)
\nonumber\\
&=& (r\,\cros\, (\one \arl a\1)(\phi\arl a\2)\psi)
\nonumber\\
&=&(r\,\cros\, (\phi\arl a)\psi)
\lb{3.3.21}
\eea
where $a\in\A$ and $\phi,\psi\in\hat\A$, and where we have used
that by \no{3.2}
\beq\lb{3.3.22}
(a\,\cros\, \e_L(b)\phi) = (a\hat\e_R(\e_L(b))\,\cros\, \phi)
\eeq
for all $a,b\in\A$ and $\phi\in\hat\A$.
Next, viewed as elements in $\A\cros\hat\A$ we have
\bea\lb{3.3.23b}
\hat e &=& (\one\,\cros\, d_R(\l)^{1/2}\hat hd_R(\l)^{1/2})
\\
\lb{3.3.23a}
e &=& (d_R(l)^{1/2}hd_R(l)^{1/2}\,\cros\, \hat\one) =
(\one\,\cros\, \xi_l)(h\,\cros\, \one)(\one\,\cros\, \xi_l)
\eea
where
\beq\lb{3.3.23}
\xi_l := \e_L(d_R(l)^{1/2})\in\hA_L.
\eeq
Since $(\one\,\cros\,\hat\A_R)$ commutes with $(\A\,\cros\,\hat\one)$ by
\thm{3.1}iv) we conclude
\bea
e\hat e e &=& \left(\one\,\cros\, \xi_ld_R(\l)^{1/2}\right)\,
\left(h\,\cros\, \xi_l\hat h\xi_l\right)\,
\left(h\,\cros\, d_R(\l)^{1/2}\xi_l\right)
\nonumber\\
&=& \left(\one\,\cros\, \xi_ld_R(\l)^{1/2}
\hat S^{-1}(\xi_l)\right)\,
\left(h\,\cros\, (\hat h\arl h)\right)\,
\left(h\,\cros\, \hat S(\xi_l)d_R(\l)^{1/2}\xi_l\right)
\nonumber\\
&=&(\one\,\cros\, \xi_l\phi^*)\,(h\,\cros\, \hat\one)\,
(\one\,\cros\, \phi\xi_l)\ .
\lb{3.3.24}
\eea
Here we have used \lem{2.20}ii)+ii'), $\hat S^{\pm 1}(\xi_l)\in\hA_R$ and
\Eq{3.3.21} in the second line and, using
$\hat h\arl h=\hat g_R^2$ in the third line,
we have introduced $\phi\in\hA_R$ given by
$$
\phi=\hat g_R\hat S(\xi_l)d_R(\l)^{1/2}\ .
$$
Comparing \no{3.3.24} with \no{3.3.23a} we realize that
\no{3.3.14} follows provided $\phi=\hat\one$.
To prove this we use
$\hat g_R=\e_R(g_L),\ \hat S(\xi_l)=\e_R(S^{-1}(d_R(l)^{1/2}))$
and the fact that $g_L$ implements $S^2$ on $\A_L$ to compute
\beanon
\phi &=&\e_R(g_L S^{-1}(d_R(l)^{1/2}))d_R(\l)^{1/2}
\\
&=&\e_R(S(d_R(l)^{1/2})g_L)d_R(\l)^{1/2}
\\
&=& \hat\one
\eeanon
where the last equation follows from \no{3.3e}.
This concludes the proof of part iii) and therefore of
\thm{3.18}.
\qed

\bsn
We are left to show
\Lem{3.19}
{Let $d\in\A_R$ and $b=S(d)g_L\in\A_L$.
Then $b=b^*\Leftrightarrow d=d^*$ and $b>0 \Leftrightarrow d>0$.}
\proof
We use that $g_L\in\A_L$ is positive and implements $S^2$ on
$\A_L$. Hence $b^*=g_LS^{-1}(d^*)=S(d^*)g_L$ implying
$b=b^*\Leftrightarrow d=d^*$. Assume now $d=c^*c,\ c\in\A_R$.
Then $b=S(c)g_LS^{-1}(c^*)=S(c)g_LS(c)^*>0$.
Converseley, if $b=a^*a,\ a\in\A_L$, then $d=S^{-1}(a^*ag_L^{-1})
=S(a)g_R^{-1}S(a)^*>0$, where we have used
$g_R=S^{\pm 1}(g_L)>0$.
\qed

\subsection{The regular representation}

We now show that if $\M$ is a von-Neumann (or $C^*$-) algebra then also
$\M\cros\A$ naturally becomes a von-Neumann (or $C^*$-)
algebra, respectively.
As for crossed products by duals of compact groups $G$, which may be defined as
subalgebras of $\M\otimes \B(L^2(G))$, we first provide what we
call the {\em regular homomorphism}
$\Lambda_{cros} :\M\cros\A\to \M\otimes End~\A$.

To this end we introduce on $\A$ the scalar product
$$(a,b) := \bra \hat h\mid a^* b\ket$$
where $\hat h\in \hat\A$ is the normalized
two-sided Haar integral. Elements of the Hilbert space
$L^2(\A,\hat h)$ are denoted by $\mid a\ket,\ a\in\A$.
On $L^2(\A,\hat h)$ we define
the following operators
\bea
\lb{3.11}  \ell(a) \mid b\ket &:=& \mid ab\ket\\
\lb{3.12}  \tau_R(\phi)\mid b\ket &:=& \mid\phi\arr b\ket\\
\lb{3.13}  \tau_L(\phi)\mid b\ket &:=& \mid b\arl \hat S^{-1}
(\phi)\ket
\eea
where $a,b\in\A$ and $\phi \in\hat\A$. Using
$(\phi^*\arr b)^*=\hat S^{-1}(\phi)\arr b^*,\
(b\arl \phi^*)^*=b^*\arl \hat S^{-1}(\phi)$, and the fact
that as a two-sided integral $\hat h$ satisfies \no{1.18} and
\no{1.21}, one immediately
checks that
$\tau_R$ an $\tau_L$ provide two commuting faithful *-representations
of $\hat A$, whereas $\ell$
of course provides the GNS-representation of $\A$
associated with $\hat h$.  Moreover, we have

\Lem{3.20}
{$\tau_L(\e_R(a))=\ell(a),\quad\forall a\in\A_L.$
}
{\bf Proof:} Let $a=\hat\e_L(\phi),\phi\in\hat\A$. Then by \no{2.23}
$\e_R(a)=\hat S(\phi\1)\phi\2$ and therefore
\beanon
\tau_L(\e_R(a))\mid b\ket
&=& \mid b\2 \ket \bra b\1\mid S^{-1} (\phi\2)\phi\1\ket\\
&=& \mid b\3\ket\bra b\2 S^{-1} (b\1) \mid\phi \ket\\
&=& \mid\e_L(\phi)b\ket
\eeanon
where the last equation follows from \no{1.8}. \qed

\Thm{3.21}{{\rm (The regular homomorphism)\\}
Let $\Lambda_{cros} :\M\cros \A \to \M\otimes End~\A$ be given
by
\beq \lb{3.14}
\Lambda_{cros} (m\,\cros\, a) := (id_\M\otimes \tau_L)
(\rho(m))\,(\one_\M\otimes \ell(a))
\eeq
Then $\Lambda_{cros}$ defines a (in general non-unital) injective $*$-algebra
homomorphism satisfying for all $m\in\M,\ a\in\A$ and
$\phi\in\hat\A$
\beq \lb{3.15}
\Lambda_{cros} (m\,\cros\, (\phi\arr a))=(\one_\M\otimes \tau_R(\phi\1))\,
\Lambda_{cros}(m\,\cros\, a)\, (\one_\M\otimes \tau_R(\hat S(\phi\2))
\eeq
Moreover, if $\M$ is a von-Neumann (or $C^*$-) algebra then $\Lambda_{cros}
(\M\cros \A)$ is a von-Neumann (or $C^*$-) subalgebra of $\M\otimes End~\A$.
}
{\bf Proof:} Since $\rho:\M\to\M\otimes \hat\A$ is an injective
$*$-algebra map, so is
the restriction of $\Lambda_{cros}$ to $\M$.
To show that $\Lambda_{cros}$ is well defined on $\M\cros\A$ we use
$$
\rho(a\lef \one_\M)= (id_\M\otimes a\arr)(\rho (\one_\M))
= \rho(\one_\M)(\one_\M\otimes \e_R(a))
$$
by \no{2.1} and \no{2.4}. Hence, using \lem{3.20} we get for $a\in \A_L$
$$\Lambda_{cros} (m(a\lef \one_\M)\,\cros\, b) = \Lambda_{cros} (m\,\cros\, ab)$$
which shows that $\Lambda_{cros}$ is well defined.
Clearly, the restriction $\Lambda_{cros}|(\one_\M\,\cros\, \A)$ defines a $*$-algebra
map provided the projection $P=\Lambda_{cros}(\one_\M\,\cros\, \one)$ commutes
with $\one_\M\o \ell(\A)$.
Using \no{1.6} one checks
\beq \lb{3.16}
\ell(a)\tau_L(\phi) =\tau_L(\phi\1) \bra\phi\2 \mid a\1\ket \ell(a\2)
\eeq
which together with \no{2.1} and \no{2.4} implies
\beanon
(\one_\M\o \ell(a))P &=& P(\one_\M\o\tau_L(\e_R(a\1))\ell(a\2))\\
&=& P(\one_\M\o\tau_L(\e_R(a\1 S(a\2)))\ell(a\3))\\
&=& P(\one_\M\o \ell(a))
\eeanon
where in the second line we have used \no{2.24} and \no{2.19} and in
the last line \lem{3.20} and the identity \no{1.1}.
Hence $\Lambda_{cros}|(\one_\M\,\cros\, \A)$ provides a $*$-algebra map. More
generally, \no{3.16} and \no{2.1} also imply
$$
 \Lambda_{cros} (\one_\M\,\cros\, a) \Lambda_{cros} (m\,\cros\, \one)
=\Lambda_{cros} ((a\1 \lef m)\,\cros\, a\2)\ .
$$
Hence, $\Lambda_{cros}$ provides an
algebra map, which therefore also respects the $*$-structure
\no{3.4}.
In particular
$\Lambda_{cros}(\M\cros \A)=\Lambda_{cros}(\M) (\one_\M\otimes\ell(\A))$.
Next, we use that $\tau_R(\hat \A)$ commutes with
$\tau_L(\hat\A)$ and therefore
\no{3.15} follows from
\beanon
\tau_R(\phi\1) \ell(a) \tau_R (S(\phi\2))
&=& \ell(\phi\1 \arr a)\tau_R(\phi\2
S(\phi\3))\\
&=& \ell(\hat\one\1 \phi \arr a)\tau_R(\hat\one\2)\\
&=& \ell(\phi \arr a)
\eeanon
where in the second line we have used \no{1.6}.
To prove that $\Lcros$ is injective let now $\l\in\L(\hA)$ be a
nondegenerate left integral and define
$F_\l:\M\o\End\A\to\M\o\End\A$ by
$$
F_\l(X):=[\one_\M\o\tau_R(\l\1)]\,X\,[\one_\M\o\tau_R(\hat S(\l\2))].
$$
Then \no{3.15} implies
\bleq{3.15'}
F_\l\circ\Lcros = \Lcros\circ\hat E_\l
\eeq
and by \cor{3.9} $\Lcros$ is injective.
Finally, if $\M$ is a von-Neumann (or $C^*$-) algebra, then by
definition $\r:\M\to\M\o\A$ is required to be weakly or norm
continuous, respectively, and therefore
$\Lcros(\M)\subset\M\o\End\A$ is weakly (or norm) closed.%
\footnote
{In the von-Neumann algebra case $\M\subset\B(\H)$ we consider
$\Lambda_{cros}(\M\cros \A)$ to act on $\H\otimes L^2(\A,\hat h)$.
}
Let now $l\in\L(\A)$ be the left integral dual to $\l$ and put
$$
\sum_i u_i\o v_i :=
\ell(l\2)\o\ell(S^{-1}(l\1))\in\End\A\o\End\A.
$$
Then by \no{3.15'} and \thm{3.6}i) we get for all
$X\in\Lcros(\M\cros\A)$
$$
X=\sum_i X_\l^i(\one_\M\o v_i),
$$
where $X_\l^i:=F_\l(X(\one_\M\o u_i))\in\Lcros(\M)$.
Clearly, the coefficient maps $X\mapsto X_\l^i$ are weakly (or
norm) continuous, and therefore $\Lcros(\M\cros\A)$ is weakly
(or norm) closed, respectively.
\qed

\bsn
From now on we consider $\M\subset\B(\H)$ to be a von-Neumann
algebra acting on a Hilbert space $\H$.
Denote $\H_\A :=\H \otimes L^2 (\A,\hat h)$ and
$\H_{cros}:=P\,\H_\A$, where
$P:=\Lambda_{cros} (\one_\M\,\cros\, \one)\in\B(\H_\A)$.
Then $\Lcros(\M\cros\A)\H_{cros}\subset\H_{cros}$ and we call
the resulting faithful unital $*$-representation
$$
\pi_{cros} :\M\cros\A\to\B(\H_{cros})
$$
the {\em regular representation} of $\M\cros\A$.
In this way we may identify $\M\cros\A\cong\pi_{cros}(\M\cros\A)$
as a von-Neumann algebra on $\H_{cros}$.
Note that in this way \no{3.15} implies the dual $\hA$-action on
$\M\cros\A$ to be weakly continuous as well.
We now show that if $\H$ is the GNS-Hilbert space obtained from
a normal state $\om$ on $\M$, then $\pi_{cros}$ is the
GNS-representation of $\M\cros\A$ obtained from the normal
state $\om_{cros}=\om\circ\hat E_{\hat h}$.

\Thm{3.22}{{\rm (The regular GNS-representation)\\}
Let $\Omega \in\H$ be cyclic for $\M$ and denote $\om$ the
induced state, $\om(m):=(\Om,\,m\Om),\ m\in\M$.
Also put
$\Omega_\A :=
\Omega\, \otimes \mid\one_\A\ket\in\H_\A,\
\Om_{cros}:=P\,\Om_\A\in\H_{cros}$ and
$\om_{cros}(m\cros a):=
(\Om_{cros},\,\pi_{cros}(m\cros a)\Om_{cros})$ as a state on $\M\cros\A$.
Then
\\
i) $\Om_{cros}\in\H_{cros}$ is cyclic for $\M\cros \A$.
\\
ii) $\om_{cros}=\om\circ\hat E_{\hat h}$
\\
iii)  $\| \Omega_\A\| ^2 =  \e (\one_\A)\| \Omega_{cros}\| ^2$
and $\|\Omega\| =\|\Om_{cros}\|$
\\
iv)
If $\Om$ is separating for $\M$ then $\Om_{cros}$ is separating
for $\M\cros\A$.
}
{\bf Proof:}
To prove i) we show
\bleq{3.20'}
\pi_{cros}(\M\cros\A)\Om_{cros}
\equiv\Lcros(\M\cros\A)\Om_\A\supset P\,(\M\o\ell(\A))\Om_\A
\eeq
which is clearly dense in $\H_{cros}$.
To this end let us introduce the Dirac notation
\beq\lb{3.18}
\mid m\otimes a\ket_{\omega,\hat h} := (m \otimes \ell(a)) \Omega_\A
\in\H_\A
\eeq
Then for $m_1,m_2 \in \M$ and $a,b\in\A$ the definition \no{3.14} gives
\beq \lb{3.19}
\Lambda_{cros}(m_1\,\cros\, a) \mid m_2 \otimes b\ket_{\omega,\hat h}
= \mid (S^{-1} (a\1 b\1)\lef m_1)m_2
\otimes a\2 b\2\ket_{\omega,\hat h}
\eeq
which yields
\bea
\Lambda_{cros} ((a\1 \lef m)\,\cros\, a\2) \Omega_\A
&=& \mid S^{-1} (a\2) a\1 \lef m
\otimes a\3\ket_{\omega,\hat h}\nonumber\\
&=& \mid (S^{-1} (a\2)a\1 \lef\one_\M)m\otimes
a\3\ket_{\omega,\hat h}\nonumber\\
&=&\mid (S^{-1} (a\1)\lef\one_\M) m\otimes a\2\ket_{\omega,\hat h}\nonumber\\
&=& P\mid m\otimes a\ket_{\omega,\hat h}
\lb{3.20}
\eea
where in the second line we have used \lem{2.16}i), in the third line
\no{1.33} and in the last line \no{3.19}.
This proves \no{3.20'} and therefore part i).
To prove part ii) we use \no{3.19} to get
\beanon
(\Omega_\A, \Lambda_{cros}(m\,\cros\, a) \Omega_\A) &=&
\omega(S^{-1}(a\1)\lef m)\bra\hat h\mid a\2\ket\\
&=&\omega (S^{-1}(\hat h\arr a)\lef m)\\
&=& \omega(m((\hat h\arr a)\lef\one_\M))\\
&=& \omega_{cros} (m\,\cros\, a)
\eeanon
where in the third line we have used
$S^{-1} (\hat h\arr a)\in S^{-1}(\A_L)=\A_R$ together with
\lem{2.16}ii).
To prove part iii) we note that by definition
$\|\Omega_\A\|^2/\|\Omega\|^2
= \bra\hat h\mid \one_\A\ket=\e(\hat h\arr
\one_\A)=\e(\one_\A)$.
The identity $\|\Omega\|=\|\Omega_{cros}\|$ follows from
$\om_{cros}(\one_\M\cros\one_\A)=\om(\one_\M)$.
Finally, part (iv) follows, since the conditional expectation
$\hat E_{\hat h}$ is faithful, and therefore $\om_{cros}$ is
faithful if $\om$ is faithful.
\qed

\bsn
Le us check the result of this construction for our example in
Sect. 2.5. Using the identification \no{3.7a} and applying the
formula \no{3.7b} for $\l_{Haar}\in\L(\hA)$ we conclude from
\no{2.6.4} that $\hat E_{\l_{Haar}}:\M\cros_\al\,G\to\M$ is
given by the standard formula
\bleq{3.54}
\hat E_{\l_{Haar}}(m\cros_\al\,g)=m\delta(g)
\eeq
and therefore $\om_{cros}(m\cros_\al\,g)=\om(m)\delta(g)$.
Thus, if $\H=L^2(\M,\om)$ then
$\H_{cros}:=L^2(\M\cros_\al\,G,\,\om)\cong\H\o\CC G$
and we obtain the standard definition [St, Pe, ES] of
$\M\cros_\al\,G$ as a von-Neumann algebra acting on $\H_{cros}$.



\sec{Jones Extensions}

In this section we will always assume the setting of \thm{3.22}, i.e. we
will consider $\M\cong \pi_\omega(\M)$ as a concrete von-Neumann algebra
acting on $\H_\omega=L^2(\M,\omega)$, where $\omega$ is a faithful normal
state on $\M$ and where $\pi_\omega$ denotes the GNS representation
associated with $\omega$.
For $m\in\M$ we will use the notation
\beq\lb{4.1}
\mid m\ket_\omega := \pi_\omega(m)\Omega_\omega\in\H_\omega
\eeq
where $\Omega_\omega\in\H_\omega$ is the cyclic GNS-vector satisfying
$(\Omega_\omega,\pi_\omega(m)\Omega_\omega)=\omega(m)$.
By \thm{3.22}
$\H_{cros}:= P(\H_\omega\otimes L^2 (\A,\hat h))$ is the GNS-Hilbertspace
associated with the induced faithful normal state
$\omega_{cros} =\omega\circ \hat E_{\hat h}$
on $\M\cros \A$, where $P=\Lambda_{cros}(\one_\M\,\cros\, \one)$ and
where $\hat h\in\L(\hat\A)$ is the normalized Haar integral in
$\hat\A$.
Accordingly, for $(m\,\cros\, a)\in\M\cros \A$ we will denote
\beq
\mid  m\,a\ket_{cros}
:= \Lambda_{cros}(m\,\cros\, a)\Omega_\A\equiv \mid S^{-1} (a\1)\lef
m\otimes a\2\ket_{\omega,\hat h} \in\H_{cros}
\eeq
where the second identity follows from \no{3.19}.
Moreover, from now on we will always assume $\omega$ to be
{\em $\A$-invariant}, by which we mean
\beq\lb{4.3}
\omega\circ E_h =\omega
\eeq
Note that this can always be achieved by replacing $\omega$ by its
``$A$-average" $\omega\circ E_h$. Also note that the state
$\omega_{cros}=\om\circ\hat E_{\hat h}$
on $\M\cros\A$ is $\hat \A$-invariant with
respect to the dual $\hat A$-action on $\M\cros\A$. Since
$\H_{cros}=L^2(\M\cros\A,\omega_{cros})$, this
means that by taking alternating crossed products our construction
will proceed ``up the ladder".

The aim of this section is to relate our theory of crossed
products by weak Hopf algebras with the theory of Jones
extensions. To this
end we will show in Section 4.1 that for any faithful normal
$\A$-invariant state $\omega$ on $\M$ the GNS-representation
$\pi_\omega$ of $\M$ extends to a representation of $\M\cros\A$
(still denoted $\pi_\omega$) such that
$\M_1\equiv \pi_\omega (\M\cros\A)$
coincides with the basic Jones construction
for the inclusion
$\pi_\omega(\M^\A)\equiv\M_{-1}\subset \M_0\equiv \pi_\omega(\M)$.
Moreover, $\pi_\omega$ may be identified with a
subrepresenation of $\pi_{cros}$
and $\M_1$ may be identified with the ideal
$\M e_h\M\subset\M\cros\A$ where $e_h=(\one_\M\,\cros\, h)$
is the Jones projection associated with $E_h:\M\to\N$.
We will also show that the Jones triple
$\M_{-1}\subset\M_0\subset\M_1$ has depth 2 and finite index.
In fact, we will have $\Ind E_h\le\tau_\lef(\Ind h)$, where
equality holds if and only if $\M_1$ is a faithful image of
$\M\cros\A$.

In Section 4.2 we introduce an appropriate notion of
outerness for weak Hopf algebra actions and show that an action
is outer if and only if the relative commutant of $\M$ in
$\M\cros\A$ is minimal, i.e. iff
$\M'\cap\M\cros\A=C(\M)(\one_\M\,\cros\, \A_R)$.
In particular this will
imply that if $\A$ is a pure Hopf algebra acting outerly on a
factor $\M$ then also $\N\equiv\M^\A$ and $\M\cros\A$ are factors and
therefore $\M_1\cong\M e_h \M = \M\cros\A$. Moreover, in this
case also $\hat\A$ is pure and acts outerly on
$\M\cros\A$ and the iterated inclusions
\beq\lb{4.0}
\N\subset\M\subset\M\cros\A\subset(\M\cros\A)\cros\hat\A\subset\dots
\eeq
provide a Jones tower of factors.

To determine the relative commutants of this tower we
investigate in Section 4.3 the actions on $\pi_\om (\A_{L/R})$
of the modular group
$\Delta_{\M,\om}^{it}$ and the modular conjugation $J_{\M,\om}$
associated with $(\M,\om)$. This will prove that under the above setting
$\N'\cap(\M\cros\A)=\A,\ \N'\cap\M=\A_L$ and
$\M'\cap(\M\cros\A)=\A_R$.

Finally, in Section 4.4 we generalize these results to
non-pure weak Hopf algebras $\A$ acting on non-factorial
von-Neumann algebras $\M$ by imposing as a regularity condition
$C(\M)=(\A_L\cap\A_R)\lef\one_\M$
in addition to standardness and outerness.
Note that in view of \cor{2.17}i) this means that the
center of $\M$ is required to be as small as possible.
Again, this kind of regularity proceeds up the tower \no{4.0},
which under these conditions
still provides a Jones tower with lowest relative
commutants given as above. Moreover, in such a setting
we will have $C(\M_{2i})\cong\A_L\cap\A_R,\
C(\M_{2i+1})\cong\hat\A_L\cap\hat\A_R$ and
$C(\M_{i})\cap C(\M_{i+1})\cong\A_L\cap\A_R\cap C(\A)
\cong\hat\A_L\cap\hat\A_R\cap C(\hat\A)$, thus producing
precisely the scenario described in Section 1.1.
Finally, in such a ``regular scenario" the derived tower
$\N'\cap\M_i,\ i\ge 0$, is given by
\beq\lb{4.0'}
\A_L\subset\A\subset\A\cros\hat\A\subset
(\A\cros\hat\A)\cros\A\subset\dots\ .
\eeq

In [NSW] we will show that conversley any Jones tower of finite
index and depth 2 with finite dimensional centers
appears this way, where $\A$ and
$\hat\A$ are given by $\A=\N'\cap\M_1$
and $\hat\A=\M'\cap\M_2$.

\bsn
To simplify our notation we will from now on write elements
$(m\,\cros\, a)\in\M\cros\A$
as products of $m\equiv (m\,\cros\, \one)\in\M$ and
$a\equiv (\one_\M\,\cros\, a)\in\A/I_\lef$, whenever there is no
confusion possible. Note that for $a\in\A_L/Ker\mu_\lef$ this leads to
the identification $a\equiv a\lef\one_\M$.

\subsection{The basic construction}
We start with pointing out that if $\omega$ is $\A$-invariant,
then the restriction
of $\pi_{cros}$ to $\M\subset \M\cros\A$
contains $\pi_\omega$ as a subrepresentation.
First we need

\Lem{4.1}{{\rm ($\A$-invariant states)\\}
The following conditions for a state $\omega$ on $\M$ are
equivalent
\\
i) $\omega$ is $\A$-invariant, i.e. $\om=\om\circ E_h$
\\
ii) $\omega(a\lef m)=\omega((S^{-1}(a)\lef \one_M)m),\
\forall a\in\A,\forall m\in\M$
\\
iii) $\omega(a\lef m)=\omega(m(S(a)\lef\one_M)),\
\forall a\in\A,\forall m\in\M$
}
{\bf Proof:}
If $\omega$ is $\A$-invariant then ii) and iii) follow from
(\ref{1.22}) and (\ref{1.23}), \no{1.33} and \no {2.13} and
\lem{2.16}.
Conversely, putting $a=h$, ii) or iii) immediately imply i). \qed

\bsn
Note that in our example of the partly inner group action
$(\al,\,u)$ the conditions of \lem{4.1} just mean
$$
\om(u(h)\al_g(m))=\om(u(g^{-1}hg)m)=\om(mu(g^{-1}hg)),\quad
h\in H,\,g\in G,
$$
which of course is equivalent to $\om$ being $\al_g$-invariant,
$\forall g\in G$.
Next, to show that the restriction $\pi_{cros}\mid \M$
contains $\pi_\omega$ as a
subrepresentation we provide a unit vector $\Phi\in\H_{cros}$
implementing $\omega$.
Making the Ansatz
\beq
\Phi := \pi_{cros} (l) \Omega_{cros} \in\H_{cros} \lb{4.6}
\eeq
for some left integral $l=hd_L,\ d_L\in\A_L$, we get
\beanon
(\Phi,\pi_{cros} (m)\Phi) &=& \omega_{cros}(l^* m l)
= \om_{cros}( d_L^*(h\re m) l)\\
&=&(\omega \circ \hat E_{\hat h})(E_h(m) d_L^* h d_L)\\
&=& \omega (E_h(m) d_L^* (\hat h\arr h) d_L)
\eeanon
where we have used \no{4.2} in the second equation,
\thm{3.22}ii)
and the fact that $h\re m\equiv E_h(m)\in\M^\A$ commutes with $d_L^*\in\A$
in the third equation and the definition \no{3.8}
of $\hat E_{\hat h}$ in the last equation.
If we now choose $d_L :=g_L^{-1}\equiv
(\hat h\arr h)^{-1/2}$  we get
\beq
(\Phi,\pi_{cros} (m)\Phi)=\omega (E_h(m))=\omega (m)\ .\lb{4.4}
\eeq
Thus, with this choice $\Phi\in\H_{cros}$
indeed implements $\omega$ and we arrive at

\Prop{4.2}
{\
\\
i) Let $l_0 :=hg_L^{-1}\in\L(\A)$.
Then $V:\H_\omega \to \H_{cros},\
V\,|m\ket_\omega := |m\,l_0 \ket_{cros}$, defines an
isometry satisfying $\forall m\in\M$ and $\forall a\in\A$
\bea
V^*\,|m\,a\ket_{cros} &=& |m\mu_\lef(ag_L)\ket_\om
\lb{4.5a}\\
VV^*\,|m\,a\ket_{cros} &=& |m\,ag_Ll_0\ket_{cros}
\lb{4.5b}\\
\pi_\omega(m) &=& V^*\pi_{cros} (m) V\ .\lb{4.5}
\eea
\\
ii) $V\H_\omega\subset\H_{cros}$ is invariant under $\pi_{cros}(\M\cros\A)$
and therefore $\pi_\omega$ extends to a representation (still denoted by
$\pi_\omega)$ of $\M\cros \A$ on $\H_\omega$ by putting
$\pi_\omega(ma):=V^*\pi_{cros}(ma)V$, yielding
\beq\lb{4.7}
\pi_\omega (ma)\,|m'\ket_\omega = |m(a\lef m')\ket_\omega
\eeq
iii) $\M h\M\subset\M\cros\A$ is an ideal orthogonal to
$\Ker\pi_\om$, i.e.
\beq\lb{4.7a}
(\M h\M)(\Ker\pi_\om) = (\Ker\pi_\om)(\M h\M) = 0
\eeq
}
{\bf Proof:}
We have $V\pi_\omega(m)\Omega_\omega = \pi_{cros}(m)\Phi$ and since
$\Phi=\pi_{cros}(l_0)\Omega_{cros}$
implements $\omega$, $V$ extends to an isometry intertwining
$\pi_\omega$ and $\pi_{cros}$.
Moreover, for $m_1,m_2\in\M$ and $a\in\A$ we have
\beanon
\bra m_1\mid V^*\mid m_2a\ket_\om &=&
\bra m_1l_0\mid m_2a\ket_{cros}
=\om_{cros}(l_0^*m_1^*m_2a)\\
&=&\om_{cros}\left(l_0^*(S^{-1}(a)\lef(m_1^*m_2))\right)
=\om\left((g_L\re\one_\M)(S^{-1}(a)\lef(m_1^*m_2))\right)
\\
&=&\om(m_1^*m_2\mu_\lef(ag_L))
\eeanon
where we have used the adjoint of \no{4.2} and
$\om_{cros}=\om\circ\hat E_{\hat h}$ in the second line
and \lem{4.1} together with
$\hat h\arr l_0^*=g_L^{-1}(\hat h\arr h)=g_L=S^{-2}(g_L)$
in the last line.
This proves equs. \no{4.5a} and \no{4.5b} and therefore part (i).
Next, by \no{4.2} we have for $m,m'\in\M$ and
$a\in\A$
$$
\pi_{cros} (ma)\mid m'\,l_0\ket_{cros} =\mid m(a\lef m')l_0\ket_{cros}
$$
which also proves ii).
Finally, to prove iii) we note that by \no{4.2}
$\M h \subset \M \cros \A$ is a left
ideal and therefore
$\M h\M\subset\M\cros\A$ is a two-sided ideal. Moreover, for
$\sum m_ia_i\in\Ker\pi_\om$ we have by \no{4.7}
$\sum m_i(a_i\lef m')=0,\ \forall m'\in\M$, and therefore
$\sum m_ia_im'h = 0$ by \no{4.2}. Eq.  \no{4.7a} follows by
also taking the adjoint of this.
\qed

\bsn
Note that \prop{4.2}
 in particular implies that any $\A$-invariant state $\omega$
on $\M$ extends to a state $\bar\omega$ on $\M\cros \A$ by
\beq\lb{4.8}
\bar\omega(ma)
:=(\Omega_\omega,\pi_\omega (ma)\Omega_\omega)\equiv \omega(m(a\lef\one_\M)).
\eeq
In this way \no{4.7}
may also be viewed as the GNS-representation of $\M\cros\A$
associated with
$\bar \omega$.

\bigskip
We are now in the position to identify the Jones extension of
$\M^\A\subset\M$ with an ideal $\M_1\subset\M\cros\A$. To this
end we work in $\B(\H_\omega)$ and put
\bea
\M_{-1} &:=& \pi_\omega (\M^\A)\lb{4.9}\\
\M_0 &:=& \pi_\omega (\M)\lb{4.10a}\\
\M_1 &:=& \pi_\omega (\M\cros\A)\lb{4.10b}
\eea
Also, for a positive normalized and nondegenerate left integral
$l\in\L(\A)$ let $e_l=d_R(l)^{1/2}hd_R(l)^{1/2}$ be the
associated ``Jones projection" \no{3.3.1'}
and let $\l_l\in\L(\hat\A)$ denote the
p-dual of $l$, see \defi{3.17}.

\Thm{4.3}{{\rm (The basic construction)\\}
Let $\omega$ be a normal faithful $\A$-invariant state on
$\M$, and let $u_i\in\M$ be a Pimsner Popa basis for
$E_h:\M_0\to\M_{-1},\ h\in\A$ the normalized Haar integral.
Then
\\
i) The inclusion
$\M_{-1}\subset \M_0\subset \M_1$ provides a Jones triple of
finite index and depth 2.
\\
ii) $p:=\sum_i u_i hu_i^* \in \M\cros \A$ is a central projection satisfying
$\Ker \pi_\omega =(\one-p)(\M\cros \A)$ and therefore
$\M_1\cong p(\M\cros\A)$.
\\
iii)$ p (\M\cros \A)=\M h\M$ and the map $m\mapsto pm$ provides a
unital inclusion $\M\to\M h\M$.
\\
iv) For any positive nondegenerate normalized left integral $l\in\L(\A)$
the Jones projection $\eee_l\in\M_1$ associated with the
conditional expectation $E_l:\M_0\to\M_{-1}$ is given by
$$
\eee_l=\pi_\om(e_l)
$$
and the (unnormalized) conditional
expectation $E'_l:\M_1\to\M_0$ dual to $E_l$ is given by
\beq
E'_l(\pi_\om(m\,a)) = \pi_\om(m)\,\pi_\om(\hat E_{\l_l}(ap))
\eeq
v)
We have $\Index E_l\le \tau_\lef(\Ind l)$ and equality holds if and only
if $p=\one$, i.e. iff $\M_1\cong\M\cros\A$.
}
{\bf Proof:}
i) Putting
 $\eee_h=\pi_\om(e_h)\equiv\pi_\omega(h)$
we have $\eee_h\,|m\ket_\omega = |E_h(m)\ket_\omega$
and hence the basic construction associated with $E_h:\M_0\to
\M_{-1}$ is given by
$\bra \M_0,\eee_h\ket$. Since $\eee_h\in \M_1$ we have
$\bra\M_0,\eee_h\ket \subset\M_1$.
To prove the inverse inclusion it is enough to show $\pi_\omega(\A)\subset
\bra\M_0,\eee_h\ket$, or equivalently
$\bra\M_0,\eee_h\ket'\subset\pi_\omega(\A)'$.
To this end we use that $X\in\bra\M_0,\eee_h\ket'$ implies for all $m\in\M$
$$
X\,|m\ket_\omega = \pi_\omega(m)X\Omega_\omega
= \pi_\omega (m\, h)X\Omega_\omega
$$
since $X\Omega_\omega =X\eee_h\Omega_\omega =\eee_h X\Omega_\omega$.
Hence, we get
for all $a\in\A$ and $m\in\M$
\beanon
X\pi_\omega(a)\,|m\ket_\omega &=& X\,|a \lef m\ket_\omega
= \pi_\omega ((a\lef m)\, h) X\Omega_\omega\\
&=& \pi_\omega (a)\pi_\omega (m\,h) X\Omega_\omega
= \pi_\omega (a) X\,|m\ket_\omega
\eeanon
where in the third equation we have used \no{4.2}. Thus
$\bra\M_0,\eee_h\ket' \subset\pi_\omega(\A)'\cap\M_0' = \M_1'$
and therefore
$\M_1\subset\bra\M,\eee_h\ket$.

Next, we prove that the triple
$\M_{-1}\subset\M_0\subset\M_1$ has finite index and depth-2.
By definition [Oc1] the depth-2 property means that the derived tower
$$
\M_{-1}'\cap\M_0\subset\M_{-1}'\cap\M_1\subset\M_{-1}'\cap\M_2
$$
also is a Jones tower, where $\M_2\supset\M_1$ is the basic
construction for $\M_1\supset\M_0$.
Moreover, finite index holds if the ``unnormalized conditional
expectation" (more precisely the operator valued weight)
$E'_h:\M_1\to\M_0$ dual to $E_h$ is bounded, i.e.
$E'_h(\one)\equiv\Ind E_h < \infty$.
Recall that $E'_h$ is uniquely defined on $\M_0\eee_h\M_0$
by the requirement $E'_h(\eee_h) =\one$ [BDH].
Now, depth 2 and finite index can be guaranteed simultaneously
if $E'_h$
is of index-finite type with a quasi-basis which can
be chosen to lie in $\M'_{-1}\cap\M_1$.
To show that this indeed holds let now $q\in\M\cros\A$ be the
central projection onto $\Ker\pi_\om$, i.e.
$q(\M\cros\A)=\Ker\pi_\om$.
Then for any positive left integral $\l\in\L(\hat\A)$ the map
$\widetilde E_\l : \M_1\to\M_0$ given for $x\in\M\cros\A$ by
\beq
\widetilde E_\l(\pi_\om(x)) := \pi_\om\left(\hat E_\l((\one-q)\,x)\right)
\eeq
is a well defined (unnormalized) conditional expectation. Moreover,
if $\l$ is nondegenerate, then by \thm{3.6}
$\widetilde E_\l$ is of index-finite type with quasi basis
$\sum x_i\o y_i =\pi_\omega(l\2)\o\pi_\omega(S^{-1}(l\1))$,
where $l\in\L(\A)$ is the left integral
dual to $\l$ (i.e. the unique solution of $l\arr \l=\hat\one$).
In particular,
$x_i,y_i\in\pi_\omega(\A)\subset\M_{-1}'\cap\M_1$.
Choosing now $\lambda=\lambda_h$ to be the p-dual of $h$
we have $\lambda_h\arr h\equiv\l_h\arr e_h=\one$
and therefore, using $hq=0$ by
\no{4.7a},
$$
\widetilde E_{\lambda_h}(\eee_h)=
\pi_\omega(\hat E_{\lambda_h}((\one-q)\,h))
=\pi_\omega(\lambda_h\arr h)=\one\ .
$$
Thus $E'_h:=\widetilde E_{\lambda_h}$ indeed provides the
(unnormalized) conditional expectation
dual to $E_h$, which proves the finite index and depth-2 property.

To prove ii) first note that since $E_h$ has finite index the
Pimsner-Popa basis $\{u_i\}$ is finite.
By definition we have for all $m\in\M$
\beq\lb{4.12}
\sum_i u_i E_h(u_i^*m) = m=\sum_i E_h(mu_i) u_i^*
\eeq
implying $\pi_\omega(p) =\one$ and therefore $\one -p\in\Ker\pi_\om$.
Hence $\one -p=q(\one -p)=q$, since $qp=0$ by \no{4.7a}.

To prove iii) we note that by \no{4.12} and \no{4.2} $p$ is the unit
in $\M h\M$ and therefore
$\M h\M=p (\M\cros \A)\cong \M_1$.
Next, for $m\in\M$ assume $pm=0$. Then
$m=(\one-p)m\in\Ker\pi_\omega\cap\M=0$ by assumption, which
proves part iii).

For $l=h$ part iv) has already been shown in part i).
Let now $l=hd_R(l)$ be an arbitrary
positive nondegenerate and normalized left integral.
By \no{2.41} we have $E_l(m)=E_h(z_l^*mz_l),\ \forall m\in\M$,
where $z_l=d_R(l)^{1/2}\lef \one_\M\equiv\mu_\re(S(d_R(l)^{1/2}))$
(\lem{2.16}ii).
Hence, the Jones projection associated with $E_l$ is given on
$\H_\omega\equiv L^2(\M,\omega)$ by
\footnote%
{Note that in general $E_l$ does not define a selfadjoint projection in
$\B(\H_\omega)$, unless $\omega$ is $E_l$-invariant. Otherwise
let $\omega_l:=\omega\circ E_l$, then
$U|m\ket_\omega:=|mz_l^{-1}\ket_{\omega_l}$ defines an
isomorphism $U:L^2(\M,\omega)\to L^2(\M,\omega_l)$ satisfying
$U\eee_lU^{-1}|m\ket_{\omega_l}=|E_l(m)\ket_{\omega_l}$.
}
$$
\eee_l=\pi_\omega(z_lhz_l^*)
=\pi_\omega(S(d_R(l)^{1/2})hS^{-1}(d_R(l)^{1/2}))=\pi_\om(e_l)
$$
where we have used \lem{2.20}ii)+ii').
Next, let $\l_l\in\L(\hat\A)$ be the
p-dual of $l$ and use $e_l\in\M h\M$ to get $pe_l=e_l$ and
therefore $\hat E_{\l_l}(pe_l)=\l_l\arr e_l=\one$ by \defi{3.17}.
Hence, $E'_l:\M_1\to\M_0$ given by
$$
E'_l(\pi_\omega(ma)):=\pi_\omega(\hat E_{\l_l}(pma))
$$
satisfies $E'_l(\eee_l)=\one$, and therefore $E'_l$
provides the (unnormalized) conditional expectation dual to
$E_l$, proving part iv).
Finally, to prove part v) we use  $p\le\one$ to conclude
(identifying $\M\cong\pi_\om(\M)$)
\beq\lb{4.12a}
\Ind E_l \equiv E'_l(\one) =\hat E_{\l_l}(p)\le\hat
E_{\l_l}(\one)=\tau_\re(n_R(\l_l))=\tau_\lef(\Ind l)
\eeq
where we have used \no{3.9d} and \no{3.3b}.%
\footnote
{Note that by \no{3.9d} and \no{3.3b} $\Ind l=\Ind l'$, where
$l'$ is the algebraic dual of $\l_l$ in the sense of
\prop{2.23}, see also \no{3.3f}.
}
By the faithfulness of $\hat E_{\l_l}$ equality in \no{4.12a}
holds  if and only if $p=\one$.
\qed

\bsn
We remark that a similar version of \thm{4.3} holds if in
place of requiring $\M$ to be a von-Neumann algebra we demand
the conditional expectation $E_h:\M\to\M^\A$ to be of
index-finite type with quasi basis $u_i,v_i\in\M$.
Then $\pi_\om$ given by \no{4.7} may still
be considered as a well defined algebra map
$M\cros\A\to\End_r^\N(\M)$ (i.e. onto the endomorphisms of $\M$
commuting with the right multiplication by $\N$), which is in
fact surjective. Moreover, $p=\sum u_ihv_i\in\M\cros\A$
still is a central projection, which is independent of the joice
of quasi-basis and the remainig parts of \thm{4.3} hold
similarly. We leave the details to the reader.

\bsn
In the following we are mostly interested in the case where $p=\one$,
i.e. where $\pi_\omega$ gives a {\em faithful} representation of
$\M\cros\A$, for any $\A$-invariant faithful state $\om$ on $\M$.
If this holds we say that $\lef : \A\o\M\to\M$ provides a
{\em Galois action}.
We note that for ordinary Hopf algebra (co-)actions this terminology
has been introduced somewhat differently by Chase and Sweedler
[CS] prior to the appearance of Jones theory.
In Appendix A we will translate the CS-notion of Galois (co-)actions
more literally to our setting of weak Hopf algebras and show that for
actions on von-Neumann algebras $\M$ it coincides with the above
definition.

Let us now show that the action of $\A$ on $\M$ is in particular
Galois if
$\M=\N\cros\hat\A$ is itself a crossed product with canonical
$\A$-action given in the same way as in \no{3.8a}.
To see this note that in this case a quasi-basis for $E_h:\M\to\N$ is given
according to \thm{3.6}i) by
$\sum_i u_i\o v_i =
(\one_\N\,\cros\, \lambda\2)\o(\one_\N\,\cros\,\hat S^{-1}(\l\1))$,
where $\l\equiv\l_h\in\L(\hat\A)$ is the left integral dual to $h$.
Hence, in $(\N\cros\hat\A)\cros\A$ we have
$p=\l\2\,h\,\hat S^{-1}(\l\1) = \one$ due to the following

\Lem{4.4}
{Let $l\in\L(\A)$ and $\l\in\L(\hat\A)$ be a dual pair (in the
sense of \prop{2.23}) of nondegenerate left integrals.
Then in $\hat\A\cros\A$ we have for all $a\in\A$
\bea
\l\2(a\arr\hat S^{-1}(\l\1)) &=& \bra a\1\mid \hat S^{-1}(\l)\ket a\2S(a\3)
\lb{4.13}\\
\l\2\,l\,\hat S^{-1}(\l\1) &=& \one_{\hat\A\cros\A}
\lb{4.14}
\eea
}
\proof
By \no{2.24} the r.h.s of \no{4.13} is in $\A_L$ and by \no{2.26}
the l.h.s. is in $\hat\A_R$. To prove that they coincide in
$\hat\A\cros\A$ we only have to verify $\mbox{l.h.s} =
\e_R(\mbox{r.h.s})$. To this end we compute
\beanon
\e_R\left( \bra a\1\mid \hat S^{-1}(\l)\ket a\2S(a\3)\right) &=&
(a\arl\hat S^{-1}(\l))\arr\hat\one
= \hat\one\1\bra a\mid\hat S^{-1}(\l)\hat\one\2\ket
\\
&=& \l\2(a\arr\hat S^{-1}(\l\1))
\eeanon
where in the first equation we have used \no{1.33} and in the last
one \no{1.5}. This proves \no{4.13}. Eq.  \no{4.14} follows
since we have
\beanon
\l\2 l \hat S^{-1}(\l\1) &=& \l\2(l\1\arr\hat S^{-1}(\l\1))l\2
\\
&=& \bra l\1\mid S^{-1}(\l)\ket l\2 S(l\3) l\4
\\
&=& l\arl\hat S^{-1}(\l)
\\
&=& \l_R^{-1}(\l) =\one
\eeanon
where we have used \no{4.13} in the second line,
\no{1.1} in the third line and
\no{2.38} in the last line.
\qed

\Cor{4.5}
{Any tower of alternating crossed product extensions of the form
$$
\N\subset\N\cros\hat\A\equiv\M\subset\M\cros\A\subset\M\cros\A\cros\hat\A
\subset\dots
$$
provides a Jones tower of depth 2.
}
Note that putting $\N=\A_R\cong\hat\A_L$ and $\M=\hat\A$
\cor{4.5} in
particular implies that the alternating crossed products
$\hat\A_L\subset\hat\A\subset\hat\A\cros\A\subset\dots$
provide a Jones tower.

Of course, an $\A$-action on $\M$ must also be Galois if we can assure
that $\M\cros\A$ is a factor (in which case $\pi_\omega$ would
necessarily have to be faithful).
For group actions this is well known to hold if $\M$ is a factor and
the action is outer.
In the next subsection we are going to generalize this result to
weak Hopf algebra actions.

\subsection{Outer actions}

We are now looking for conditions guaranteeing
that $\M\cros\A$ (and therefore
$\N\equiv\M^\A)$ is a factor provided $\M$ is a factor. For crossed products
by finite groups this is well known to hold if and only if the group acts
outerly. Thus we are motivated to seek a sensible
notion of outerness also
for weak Hopf actions such that the above appropriately generalizes to
our setting
(It turns out that we will also have to require $\A$ to be pure). We should
mention that in [Ya2]
an action of a finite dimensional Kac algebra $\A$ on a
factor $\M$ has been called outer if $\M'\cap (\M\cros\A)=\CC$.
However, we feel
that outerness should be defined as opposed to innerness,
i.e. in the sense of
the action not being implementable by elements of $\M$.
The connection with the
above triviality of the relative commutant $\M'\cap (\M\cros\A)$ should then
become a theorem rather than a definition.
In this sense our approach will be
closer to the notion of outerness for coactions by Hopf algebras as
given in [BCM].
We also point out that in our setting we always have the lower bound
$\M'\cap (\M\cros\A)\supset(\one_M\,\cros\,\A_R)$
by \no{3.5b}, which is precisely why weak Hopf
algebras fit with reducible inclusions. Hence for us ``triviality" of
the relative
commutant should mean equality in the above inclusion,
which for a factor $\M$
will be precisely the result of \thm{4.9} below.
However, we emphasize that beyond these applications the methods
and results of this subsection apply to
arbitrary weak $\A$-module algebras $\M$.

We start by defining an {\em inner implementer} of a coaction $\hat\rho:\M
\to \M\otimes\hat\A$ to be an element $T\in\M\o\hA$
satisfying for all $m\in\M$
\beq
T(m\otimes\hat\one)=\hat\rho(m)T   \lb{4.19}
\eeq
More generally, if $\widetilde\M\supset\M$ is some algebra extension then an
implementer of
$\hat\rho$ in $\widetilde\M$ is an element $T\in\widetilde\M\otimes\hat\A$
satisfying \no{4.19} for all $m\in\M$. Let $\lef :\A\otimes\M\to\M$ be the
action dual to $\hat\rho$ and identify $\widetilde\M\otimes
\hat\A=\Hom_\CC(\A,\widetilde\M)$.
Then $T$ is an implementer of $\hat\rho$ if and only if for all
$a\in\A$ and all
$m\in\M$
\beq
T(a)m =(a\1\lef m)T(a\2)  \lb{4.20}
\eeq
We denote $\T_\lef(\widetilde M)$ the space of implementers of $\lef$ (or
$\hat\rho$,
equivalently) in $\widetilde \M$.
Obviously, $\T_\lef(\widetilde\M)$ becomes a $C(\M)$-module by putting
$(zT)(a):=zT(a)$, for $z\in C(M)$.

We now show that for any $\A$-module algebra $\M$ the space of left
integrals in $\hat\A$
always induces a nontrivial subspace of $\T_\lef(\M)$.

\Lem{4.7}
{Let $\hat\rho:\M\to\M\otimes\hat\A$ be the coaction dual to $\lef$ and
for a left
integral $\lambda \in\L(\hat\A)$ put
$T_\lambda:=\hat\rho(\one_\M)(\one_\M\otimes
\lambda)\in\M\o\hA$. Then $T_\lambda\in\T_\lef (\M)$.
}
{\bf Proof:} For $a\in\A$ we have $T_\lambda(a) =(\lambda\arr a)\lef\one_\M$.
Using $(\lambda\arr a)\in\A_L$ and \lem{2.16}i) we get for all $m\in\M$
\beanon
T_\lambda(a)m &=& (\lambda\arr a)\lef m
= [(\lambda\arr a)\1 \lef m][(\lambda\arr a)\2 \lef \one_\M]\\
&=& [a\1\lef m] T_\lambda (a\2)
\eeanon
where we have used $\Delta(\lambda\arr a)=a\1 \otimes(\lambda\arr a\2)$.
\qed

\bsn \lem{4.7} motivates to introduce
\beq
\T_\lef^0 (\M):= \hat\rho(\one_M)(\one_\M\otimes \L(\hat\A)) \lb{4.21}
\eeq
as the space of ``trivial" inner implementers of $\lef$. Note that for
$\lambda\in\L(\hat\A)$ we have $\lambda\arr\A=\A_L$ and therefore
$T_\lambda(\A)=
\M_R$. Hence $\T_\lef^0(\M)\not= 0$ and consequently
$\T_\lef(\M)\supset C(\M)\T_\lef^0(\M)\not= 0$.
Thus we are lead to the following

\Def{4.8}
{A $\A$-module algebra action $\lef:\A\otimes\M\to\M$ is called {\em outer}
(equivalently, the dual coaction $\hat\rho$ is called outer) if
the space of inner implementers
is given by $\T_\lef(\M)=C(\M) \T_\lef^0(\M)$.
}
Let us discuss this definition in the light of our example of
Sect. 2.5. In this example, by \no{4.20} and \no{2.5.2}, an
inner implementer would be a map $T:H\x G\to\M$ satisfying for
all $m\in\M$
$$
T(h,g)m={1\over |H|}\sum_{\tilde h\in H}u(h)\al_g(m)u(\tilde
h^{-1})T(\tilde h,g)
$$
Putting $m=\one_\M$ we conclude that $T$ would be of the form
$T(h,g)=u(h)v(g)$ for some function $v:G\to\M$ implementing the
action $\al$, i.e.
$$
v(g)m=\al_g(m)v(g),\quad\forall g\in G,\,\forall m\in\M.
$$
If $T=T_\l\in\T_\re^0(\M)$ for some left integral
$\l\in\L(\hA)$ of the form \no{2.6.3}, then
$$
T(h,g)=(\l\arr(h,g))\re\one_M\equiv\left\{
\ba{rcl}
0 &,& g\not\in H\\
{\hat c(g^{-1})\, |H|^{-1}\,}u(hg) &,& g\in H
\ea
\right.
$$
where we have used \no{2.5.2c}.
Thus, our weak Hopf action $\re$ associated with $(\al,u)$ is
outer if and only if the inner part of the $G$-action $\al$ is
precisely given by $H$, i.e. iff $v(g)m=\al_g(m)v(g)$ for all
$m\in\M$ implies $g\in H$ or $v(g)=0$.

We now establish the general good use of our Definitions by showing
that the space of inner implementers $\T_\lef(\M)$
is always linearly isomorphic to the relative commutant
$\M'\cap (\M\cros\A)$. To this
end let $\lambda_h\in\L(\hat\A)$ be the left integral dual to the normalized
Haar integral $h\in\A$ and let
$\hat E\equiv \hat E_{\lambda_h} :\M\cros\A\to
\M$ be the (unnormalized) conditional expection given by \no{3.8}. For
$x\in\M'\cap (\M\cros \A)$ and $a\in\A$ we then define
\beq
t_x(a):=\hat E(ax) \lb{4.22}
\eeq

\Thm{4.9}{{\rm (Outer actions)\\}
Let $\lef : \A\otimes \M\to\M$ be an $\A$-module algebra and for $x\in\M'
\cap (\M\cros\A)$ and $a\in\A$ let $t_x(a)\in\M$ be given by
\no{4.22}. Then the
assignment $x\mapsto t_x$ provides a $C(\M)$-module isomorphism $\M'\cap
(\M\cros \A)\to \T_\lef(\M)$ satisfying
\beq
t_{(\one_\M\,\cros\, \A_R)} =\T_\lef^0(\M)   \lb{4.23}
\eeq
In particular, the action $\lef$ is outer if and only if
$\M'\cap(\M\cros \A)=(C(\M)\,\cros\,\A_R)$.
}
{\bf Proof:} For $z\in C(\M)$ we clearly have $t_{zx}(a)=z t_x(a)$. To show
that $t_x\in\T_\lef(\M)$ let $m\in\M$. Then
\beanon
t_x(a)m &=& \hat E(amx)\\
&=& \hat E((a\1\lef m)a\2 x)\\
&=& (a\1\lef m) t_x (a\2)
\eeanon
proving $t_x\in\T_\lef(\M)$. Also, if $x=(\one_\M\,\cros\, b)$
for some $b\in\A_R$ then
\bea
t_{(\one_\M\,\cros\, b)} (a)
&=& \left(\one_\M\,\cros\, (\lambda_h \arr (ab))\right)
\nonumber\\
&=& \left(((\lambda'\arr a)\lef \one_\M)\,\cros\, \one\right)\lb{4.24a}
\eea
where $\lambda'=\lambda_h\varepsilon_R(b)$ and where we have used
$(\one_\M\,\cros\, \A_L)=((A_L\lef \one_\M)\,\cros\, \one)$
in $\M\cros\A$. Hence $t_{(\one_\M\,\cros\, \A_R)}\subset
\T_\lef^0(\M)$. To prove that
$\M'\cap (\M\cros \A)\ni x\mapsto t_x\in\T_\lef(\M)$ is
bijective pick $T\in\T_\lef (\M)$ and define for $a\in\A$
\beq \lb{4.24}
T'(a) := S(a\1)T(a\2)\in \M\cros \A
\eeq
Then $T'(a)\in \M'\cap \M\cros\A$ for all $a\in\A$ since for $m\in\M$ we
have
\beanon
T'(a)m &=& S(a\1)(a\2\lef m) T(a\3)\\
&=& S(a\1)a\2 mS(a\3) T(a\4)\\
&=& mS(a\1)a\2 S(a\3)T(a\4)\\
&=& m T'(a)
\eeanon
Here we have used \no {3.4b} in the second line, the fact that $S(a\1)a\2
\in\A_R$ commutes with $\M$ in the third line and Axiom IIIc
of \defi{1.1}
in the last line. In particular, if $T=T_\lambda \in\T_\lef^0(\M)$ then
\bea\lb{4.25}
T'_\lambda(a) &=& (\one_\M\,\cros\, S(a\1)(\lambda\arr a\2))\nonumber\\
&=& (\one_\M\,\cros\, S(b\1)b\2)\in(\one_\M\,\cros\, \A_R)
\eea
where $b=\lambda\arr a$ and where we have again identified
$((\A_L\lef\one_\M)\,\cros\, \one)=
(\one_\M\,\cros\, \A_L)$ in $\M\cros\A$. We now claim that for the normalized Haar
integral $h\in\A$ we get
\bea
\lb{4.26} t'_x(h) &=& x, \quad\forall x\in \M'\cap(\M\cros\A)\\
\lb{4.27} t_{T'(h)} (a) &=& T(a), \quad
\forall a\in\A, \forall T\in \T_\lef (\M)
\eea
which would prove that the inverse of the map $x\mapsto t_x$ is given by
$T\mapsto T'(h)$ and therefore also $t_{(\one_\M\,\cros\, \A_R)} =\T_\lef^0 (\M)$ by
\no{4.25} and \no{4.24a}. To prove
\no{4.26} we use $h=S^{-1}(h)$ and therefore
\beanon
t'_x (h) &=& h\2 t_x(S^{-1} (h\1))\\
&=& h\2 \hat E(S^{-1} (h\1)x)\\
&=& x
\eeanon
 where we have used that $h\2\otimes S^{-1}(h\1)$ provides a quasi-basis for
$\hat E\equiv \hat E_{\lambda_h}$ by \thm{3.6}i). To prove \no{4.27} we
compute for $T\in\T_\lef (\M)$ and $a\in\A$
\beanon
t_{T'(h)}(a) &=& \hat E(ah\2 T(S^{-1}(h\1)))\\
&=& (\one_\M\,\cros\, \lambda_h\arr h\2)(T(S^{-1} (h\1)a)\,\cros\, \one)\\
&=& (\one_\M\,\cros\, \one\2)(T(S^{-1}(\one\1)a)\,\cros\, \one)\\
&=& (\one_\M\,\cros\, a\1 S(a\2))(T(a\3)\,\cros\, \one)\\
&=& ((a\1\lef\one_\M) T(a\2)\,\cros\, \one) = (T(a)\,\cros\, \one)
\eeanon
Here we have used again $h=S^{-1}(h)$ in the first line,
 the definition \no{3.8} of $\hat E$ and the identity
$
h\1\o ah\2=S(a\1)h\1\o h\2
$
\no{1.18} in the
second line, the identity
$\Delta(\one)=\Delta(\lambda_h\arr h)=h\1\o\lambda_h\arr h\2$
in the third line, Eq.  \no{1.7} in the fourth line and Eq.  \no{1.33}
in the fifth line.
Hence we have proven \no{4.26} and \no{4.27} and therefore
\thm{4.9}.
\qed

\Corollary{4.2.4}
{Let $\M$ have trivial center and let $\A$ act outerly on $\M$.
Then
$$
C(\M\cros\A)=\one_\M\cros(C(\A)\cap\A_R).
$$
In particular, if $\A$ is pure then $\M\cros\A$ has trivial
center, and if $\A$ acts standardly then $\A$ is pure if and
only if $C(\M\cros\A)=\CC$.
}
\proof
$C(\M)=\CC$ and outerness imply
$
C(\M\cros\A)\equiv\M'\cap\A'\cap(\M\cros\A)=\one_\M\cros(C(\A)\cap\A_R)
$
and by \prop{2.19}iii)
$\A$ is pure if and only if $C(\A)\cap \A_R=\CC$.
\cor{4.2.4} follows, since for standard actions
$\one_\M\cros\A\cong\A$.
\qed

\bsn
Applying these results to our example of a partly inner group
action $(\al,u)$ on a factor $\M$ we recover the identity
\no{0.7c}.
Indeed, by the remark following \defi{4.8}, in this case our
outerness condition is equivalent to $H\subset G$ being the
kernel of $\pr_{\Out\M}\circ\al: G\to\Out\M$ and therefore
\thm{4.9} gives
$
\M'\cap\M\cros_\al\,G=\one_\M\cros\A_R
\equiv\span\{u(h)h^{-1}\mid h\in H\}
$
by \no{2.5.8}.
Moreover, from \cor{4.2.4} and \no{2.5.2d} we conclude
\bleq{4.29}
C(\M\cros_\al\,G)=\span_\chi\{\sum_{h\in H}\chi(h)\,u(h)h^{-1}\},
\eeq
where $\chi$ runs through the $\Ad_G$-invariant characters of
$H$. In particular, we recover the well known result that
$\M\cros_\al\, G$ is a factor, if and only if $H$ is trivial,
i.e. iff the $G$-action $\al$ on $\M$ is outer in the
conventional sense.
For pure weak Hopf algebras $\A$
this statement generalizes as follows.

\Theorem{4.10}
{An outer action of a pure weak Hopf algebra $\A$ on a factor
$\M$ is Galois, i.e.
\beq\lb{4.28}
\M^\A\equiv \N\subset \M\subset\M\>cros
\A\subset(\M\>cros\A)\cros\hat\A \subset\dots
\eeq
provides a Jones tower of factors.
Moreover, under these conditions the dual
weak Hopf algebra $\hat\A$ is also pure and
its canonical action on $\M\>cros\A$ is
again outer.
}
{\bf Proof:}
By \cor{4.2.4}, if
$\A$ is pure then $\M\cros\A$ is a factor.
Hence, by \thm{4.3} $\M^\A\subset\M\subset\M\cros\A$
provides a Jones triple and therefore $\N\equiv\M^\A$ must also be a
factor. By \cor{4.5} the sequence \no{4.28}
provides a Jones tower of factors.
By \prop{2.10} $\A$ acts standardly and by
\cor{2.17}i) and \prop{2.19}ii) also $\hat\A$  is pure. The fact
that $\hat\A$ acts again outerly on $\M\>cros\A$ will be proven as part of
\thm{4.16}ii) in Section 4.4.
\qed

\subsection{Minimal actions}

\bigskip
In \thm{4.9} we have seen that for outer actions the relative commutant of $\M$ in
$\M\>cros\A$ is minimal, i.e.
$\M'\cap(\M\>cros\A)=(C(\M)\>cros\A_R)$.
We now look at the lower relative
commutant $\N'\cap\M,\ \N\equiv \M^\A$, and recall
that it is always bigger then $C(\M)\M_R$ by \no{2.33a}.
Similar as for group
actions we say that the action of $\A$ on $\M$ is {\em minimal} if
$\N'\cap\M$
is as small as possible.
\Definition{1.5}
{An $\A$-module algebra action on $\M$ is called {\em minimal} if
  $\N'\cap\M=C(\M)\M_R$
}
As an immediate consequence of this Definition we recall from
\prop{3.7}

\Corollary{4.3.2}
{$\A$ acts minimally on $\M$ if and only if
$
\A\re(C(\M)\M_R)\subset C(\M)\M_R
$
and
$
\N'\cap(\M\cros\A)=C(\M)\M_R\cros\A\equiv C(\M)\cros\A
$.
}
\proof
This follows from \prop{3.7} and the identity
$(\M\cros\one_\A)\cap(U\cros\A)=(U\M_R)\cros\one_\A$
for any linear subspace $U\subset\M$.
\qed

\bsn
We are now aiming to prove that a Galois action is minimal if and only if
it is outer.
In particular, under the setting of \thm{4.10} this will imply
$\N'\cap \M= \A_L,\ \M'\cap(\M\cros\A)=\A_R$ and
$\N'\cap(\M\cros\A)=\A$.
To this end we recall from Jones theory that if
$J_\M$ denotes the modular conjugation
associated with the GNS-representation
$\pi_\omega$ of $\M$ on $\H_\omega\equiv
L^2(\M,\omega)$ then the basic construction for  $\N\subset\M$
is also given by
$$\M_1 := J_\M\, \N\, J_\M$$
implying
$$\N'\cap \M=J_\M(\M'\cap\M_1) J_\M$$
where we have dropped the symbol $\pi_\omega$. Thus, if
$\M_1\cong \M\>cros\A$
and $\A$ acts outerly we can determine $\N'\cap\M$ provided we know how the
modular conjugation $J_M$ acts on $\pi_\omega(\A_R)$.
To this end we introduce
on $\A$ the antilinear involution
\beq\lb{4.2.0}
 a\mapsto \bar a:=g^{1/2}a_*g^{-1/2}\equiv
g^{1/2}S(a)^*g^{-1/2}.
\eeq
where $a_*:=S(a)^*$ and $g=g_Lg_R^{-1} \in\A$
is the positive group-like element implementing
the square of the antipode, see \thm{2.24}.
Note that $a_{**}=a$, $(\A_L)_*=\A_R$ and $g_*=S^{\pm 1}(g)=g^{-1}$, implying
indeed $\bar{\bar a}=a$ as well as the identities
\beanon
\overline{ab} =\bar a\bar b\quad &,& \quad (\bar a)^*=\overline{a^*}
\\
a\in\A_L &\Leftrightarrow & \bar a\in\A_R\ .
\eeanon

\Proposition{4.11}
{Under the setting of \thm{4.3} let $J_\M\equiv J_{\M,\om}$ and
$\Delta_\M\equiv\Delta_{\M,\om}$ denote the modular
conjugation and the modular operator, respectively, associated
with $(\M,\om)$.
Then for all $a\in\A_L\A_R$
\bea
\lb{4.2.1}
\Delta_\M^{it} \pi_\omega (a)\Delta_\M^{-it}
&=& \pi_\omega (g^{it} ag^{-it}),
\quad\forall t\in\RR \\
\lb{4.2.2}
J_\M \pi_\omega(a) J_\M &=& \pi_\omega (\bar a)
\eea
}
{\bf Proof:}
Let us first consider $a\in\A_L$, in which case it suffices to
prove \no{4.2.1} with $g$ replaced by $g_L$.
\lem{4.1} implies
\bleq{4.2.2a}
\om((a\lef\one_\M)m) = \om(m(S^2(a)\lef\one_\M)),
\quad\forall a\in\A,\ m\in\M
\eeq
and for $a\in\A_L$ \lem{2.16}i) gives
\bleq{4.2.2b}
\pi_\om(a)=\pi_\om(a\lef\one_\M),\quad\forall a\in\A_L.
\eeq
Putting $a=g_L=S^2(g_L)$ \no{4.2.2a} and \no{4.2.2b} imply
by the Pedersen-Takesaki Theorem [PeTa, St]
\bleq{4.2.2c}
\Delta_\M^{it} \pi_\omega (g_L)\Delta_\M^{-it} =\pi_\om(g_L),
\quad\forall t\in\RR .
\eeq
More generally, define a new faithful normal state $\om'$ on
$\M$ by
$$
\om'(m):=\om((g_L^{-1/2}\lef\one_\M)m(g_L^{-1/2}\lef\one_\M))
\equiv\om((g_L^{-1}\lef\one_\M)m).
$$
Since $g_L$ implements $S^2$ on $\A_L$, \Eq{4.2.2a} implies
$$
\om'((a\lef\one_\M)m) = \om'(m(a\lef\one_\M))
\quad\forall a\in\A_L,\ m\in\M .
$$
Again by [PeTa] we conclude that $\pi_\om(\A_L)$ is invariant 
under the modular group of $\om'$, which due to
\no{4.2.2c} is implemented by
$$
\Delta_{\M,\om'}^{it}=\Delta_{\M,\om}^{it}\pi_\om(g_L^{-it}).
$$
This proves \no{4.2.1} for all $a\in\A_L$.
To prove \no{4.2.2} for $a\in\A_L$ let
$S_\M=J_\M\Delta_\M^{1/2} =\Delta_\M^{-{1/2}} J_\M$ be the
Tomita operator for $\M$ on $L^2 (\M,\omega)$.
Then for all $m\in\M$ and $a\in\A$ we get from \no{1.32}
$$
S_\M\pi_\omega (a_*) S_\M| m\ket =
|(a_*\re m^*)^*\ket=\pi_\omega (a)|m\ket.
$$
Using $(\bar a)_*= g^{-1/2}a g^{1/2}$ we conclude specifically
for $a\in\A_L$
\beanon
\pi_\om(\bar a)  \mid m\ket &=&
S_\M\pi_\omega (g^{-{1\over 2}} a g^{1\over2})S_\M \mid m\ket
\\
&=& J_\M\Del_\M^{1\over 2}\pi_\omega (g^{-{1\over 2}} a g^{1\over2})
\Del_\M^{-{1\over 2}}J_\M\mid m\ket
\\
&=&J_\M\pi_\om(a)J_\M\mid m\ket
\eeanon
for all $m\in\M$, where the last equation follows from
\no{4.2.1} by analytic continuation.
This proves \no{4.2.2} for $a\in\A_L$, and since $J_\M^2=\id$
also for $a\in\A_R$, whence for all $a\in\A_L\A_R$.
Finally, let $a\in\A_R$, then $\bar a\in\A_L$ and using
$g\in\A_L\A_R$ and $\bar g=g^{-1}$
\beanon
\Delta_\M^{it} \pi_\omega (a)\Delta_\M^{-it} &=&
\Delta_\M^{it} J_\M\pi_\omega (\bar a)J_\M\Delta_\M^{-it}
\\
 &=& J_\M\Delta_\M^{it} \pi_\omega (\bar a)\Delta_\M^{-it}J_\M
\\
&=& J_\M\pi_\omega (g^{it}\bar a g^{-it})J_\M
\\
&=&\pi_\omega (g^{it}a g^{-it}).
\eeanon
This proves \no{4.2.1} also for $a\in\A_R$ and therefore for
all $a\in\A_L\A_R$.
\qed

\bsn
We remark without proof that if $\M=\N\cros\hA$ with canonical
$\A$-action, then \prop{4.11} holds for all $a\in\A$.


\Cor{4.12}
{An outer action of a weak Hopf algebra $\A$ on a von-Neumann
algebra $\M$ is minimal.
A Galois action on $\M$ is minimal if and only if it is outer.
}
{\bf Proof:}
If $\A$ acts outerly, then using the notation \no{4.9} -
\no{4.10b} we have $\M_0'\cap\M_1=\pi_\om(C(\M)\A_R)$ by
\thm{4.9} and  \prop{4.11} gives
\beanon
\pi_\om(\N'\cap\M) &\equiv& J_\M(\M_0'\cap\M_1)J_\M
= J_\M\,\pi_\om(C(\M)\A_R)\,J_\M
\\
&=& \pi_\om(C(\M)\A_L)
= \pi_\om(C(\M)\M_R)
\eeanon
implying $\N'\cap\M=C(\M)\M_R$.
Converseley, if $\A$ acts minimally we get by the same arguments
$\pi_\om(\M'\cap\M\cros\A)=\pi_\om(C(\M)\A_R)$ implying
outerness provided the action is also Galois (i.e. provided
$\pi_\om:\M\cros\A\to\M_1$ is faithful).
\qed

\bsn
Note that for the partly inner group action $(\al,u)$ on a
factor $\M$ \cor{4.12} proves the assertion \no{0.7b}, i.e.
$\N'\cap\M = \span\{u(h)\mid h\in H\}$, and together with
\cor{4.3.2} also the assertion \no{0.7a}, i.e.
$\N'\cap(\M\cros G)=\span\{u(h)g\mid h\in H,\ g\in G\}$, where
$\N\equiv\M^G$.
For general weak Hopf algebras $\A$
acting outerly on a factor $\M$ we conclude

\Corollary{4.13}
{Let $\A$ be a weak Hopf algebra acting
outerly on a factor $\M$ and let $\N\equiv \M^\A$. Then
\bea
\lb{4.2.3}   \N'\cap \M=\M_R &\equiv& \one_\M\cros\A_L\\
\lb{4.2.4}   M'\cap \M\>cros \A &=&\one_\M\cros \A_R\\
\lb{4.2.5}   \N'\cap \M\>cros \A &=& \one_\M\cros\A
\eea
}
{\bf Proof:} Eq. \no{4.2.3} follows from minimality, \no{4.2.4}
follows from \thm{4.9} and
\no{4.2.5} follows from \cor{4.3.2}. \qed

\subsection{Regular actions}

We are finally going to generalize the results of \thm{4.10}
and \cor{4.13} to weak Hopf actions on non-factorial
von-Neumann algebras $\M$, provided the center of $\M$ fits
``regularly" with the $\A$-action as follows.

\Def{4.14}
{An $\A$-module algebra $\M$
is called {\em regular} if $\A$ acts standardly and outerly on $\M$ and if
the center of $\M$ is given by
\beq\lb{4.4.1}
C(\M) = (\A_L\cap \A_R) \re\one_\M.
\eeq
}
We recall that the action $\re$ is called standard if $\mu_\re :\A_L\ni a
\mapsto (a\re\one_\M)\in\M_R$ is an isomorphism.
Also note that by \cor{2.17}i)
we always have $(\A_L\cap\A_R)\re\one_\M\subset C(\M)$ and therefore for
regular actions the center of $\M$ is ``as small as possible".
In particular, under the conditions of \thm{4.10}, i.e.
if $\A$ is pure and acts outerly on a factor $\M$,
then it acts regularly (since for pure $\A$
standardness follows from \prop{2.10}).

The following Theorem substantiates the
full scenario described
in our introductory motivation in Section 1.1.

\Theorem{4.16}{{\rm (Jones towers by regular crossed products)}\\
Under the conditions of \thm{4.3}
let $\A$ act regularly on $\M$ and put $\N\equiv
\M^\A$. Then
\\
i) The $\A$-action is Galois, i.e.
$\N\subset \M \subset\M\>cros \A$ is a Jones triple.
\\
ii) The dual action of $\hat\A$ on $\M\>cros\A$ is also regular.
\\
iii) The relative commutants satisfy
\bea
\lb{4.3.1}  \N'\cap\M=\M_R &\equiv& \A_L\\
\lb{4.3.2}  \N'\cap(\M\>cros\A) &=&\A\\
\lb{4.3.3}  \M'\cap (\M\>cros\A) &=&\A_R
\eea
iv) The centers of $\N,\M$ and $\M\>cros \A$ satisfy in $\M\>cros\A$
\bea
\lb{4.3.4}
C(\N) &=&\A_L\cap C(\A)
\\
\lb{4.3.5}
C(\M) &=& \A_L\cap \A_R
\\
\lb{4.3.6}
C(\M\>cros \A) &=& \A_R\cap C(\A)
\\
\lb{4.3.7}
C(\N)\cap C(\M) &=& C(\M)\cap C(\M\>cros \A)
= \A_L\cap \A_R\cap C(\A)
\eea
In particular, $\M$ is a factor iff $\hat\A$ is pure and $\N$
(equivalently $\M\cros\A$) are factors iff $\A$ is pure.
\\
v) The derived tower $\N'\cap \M_i,\ i=0,1,2,\dots$ is given by
$$
\A_L\subset\A\subset\A\>cros \hat\A\subset
(\A\>cros\hat\A)\>cros \A \subset\dots
$$
}
{\bf Proof:}
Recall from \thm{3.1} that for standard actions the embedding
$\A\to \M\>cros \A$ is injective and we may in particular identify
$\M_R\equiv \mu_\re (\A_L)\subset \M $
with $\A_L\equiv \M\cap \A \subset \M\>cros \A$.
Hence, if $\A$ acts regularly, then $C(\M)= \A_L\cap \A_R \subset\A_R$
and therefore
$\M'\cap (\M\>cros \A)=\A_R$ by \thm{4.9}, proving \no{4.3.3}.
By \cor{4.12} $\A$ acts minimally, yielding \no{4.3.1} and by
\cor{4.3.2} also \no{4.3.2}. This proves part (iii).
Hence, $C(\M\>cros\A)\equiv \M'\cap\A'\cap\M\>cros\A=\A_R\cap\A'$,
proving also \no{4.3.6}.
To prove part (i)
let now $p\equiv \Sigma u_i hu_i^* \in\M\>cros \A$ be given
as in \thm{4.3}. Then $p\in C(\M\>cros \A)=\A_R \cap C(\A)$ and therefore
we get in $\H_\omega\equiv L^2(\M,\omega)$
\bleq{4.51}
|\one_\M\ket_\omega
=\pi_\omega (p)|\one_\M\ket_\omega =| p\re\one_\M\ket_\omega =
|S(p)\re\one_\M\ket_\omega
\eeq
where we have used \no{4.7} and \lem{2.16}ii). Since $\mu_\re$ is injective
and $S(p)\in\A_L$ this proves $p=\one$ and therefore part i).
To prove part ii) we note that
\\
1.) If $\A$
acts standardly on $\M$ then so does $\hat\A$ on $\M\>cros\A$
(since $\hat\A_L\cong \A_R\equiv(\M\cros\A)_R$).
\\
2.) $\hA$ acts minimally on $\M\>cros\A$, since the
fixed point subalgebra under this action is given by $\M$
and since \prop{3.4} and Eq. \no{4.3.3} imply
$\M'\cap\M\>cros \A=\A_R=(\M\>cros\A)_R.$
\\
3.) $\hA$ acts outerly on $\M\>cros\A$ by \cor{4.12}, since
the $\hat\A$-action on $\M\>cros \A$ always is Galois by
\cor{4.5}.
\\
4.) The dual of \lem{2.18} together with Eq. \no{4.3.6} gives
$$
C(\M\cros \A) =C(\A)\cap \A_R=(\hat\A_L\cap\hat\A_R)\re\one_{\M\>cros \A}.
$$
Hence, by 1.), 3.) and 4.)
$\hat\A$ acts regularly on $\M\>cros\A$,
proving part ii).
To prove the remaining identities in part iv), Eq.
\no{4.3.5} holds by definition of regularity
and Eq. \no{4.3.4} follows from \no{4.3.6} and \no{4.2.2},
since in any Jones
triple $\N\subset\M\subset\M_1$ we have $C(\N)=J_\M C(\M_1) J_\M$
and since the
conjugation $a\mapsto \bar a$ maps $\A_R\cap C(\A)$ onto $\A_L\cap C(\A)$.
Eq.  \no{4.3.7} follows trivially from \no{4.3.4} - \no{4.3.6} and the
statements
about pureness of $\A$ and $\hat\A$ follow from \prop{2.19}.
Finally, part v)
follows by induction, using the arguments in the proof of \prop{3.7}
and the fact that the depth 2 property implies
$\N'\cap \M_{i+1} = (\N'\cap\M_i)\vee(\M'_{i-1}\cap\M_{i+1})$.
\qed

\bsn
We remark, that in the above proof
standardness of the $\A$-action was essentially
only needed to conclude $p=\one$ (i.e. the Galois property)
from \no{4.51}.
More generally, for this conclusion
it would have been enough to just require
$\Ker\mu_\re\cap C(\A)=0$, which however, due to \prop{2.10},
already implies standardness.

Finally, we show that for the above scenario our regularity
conditions of \defi{4.14} are in fact also necessary.
More precisely, we have

\Proposition{4.17}
{For an $\A$-module von-Neumann algebra
$\M$ the following conditions are equivalent
\\
i) The $\A$-action is Galois and
$\M'\cap(\M\cros\A)=\one_\M\cros\A_R\cong\A_R$.
\\
ii) $\A$ acts standardly and
$\M'\cap(\M\cros\A)=\one_\M\cros\A_R$
\\
iii) $\A$ acts regularly.
}
\proof
{iii)$\Rightarrow$ i):} Follows from \thm{4.16}.
{i)$\Rightarrow$ ii):}
If the action is Galois then $\pi_\om$ in \thm{4.3}
represents $\M\cros\A$ faithfully implying
$
\one_\M\cros\A_R\cong\pi_\om(\A_R)\cong\pi_\om(\A_L)\cong\M_R
$
by \no{4.2.2}. Hence, $\M_R\cong\A_{L/R}$ and therefore
standardness follows from $\one_\M\cros\A_R\cong\A_R$.
{ii)$\Rightarrow$ iii):}
If $\M'\cap(\M\cros\A)=\one_\M\cros\A_R$ then, identifying
$\M\equiv\M\cros\one_\A$,
$$
C(\M)=\M\cap(\one_\M\cros\A_R)=\one_\M\cros(\A_L\cap\A_R)
=(\A_L\cap\A_R) \re\one_\M
$$
and by \thm{4.9} the action is outer, whence regular.
\qed

\bsn
Note that the action in our example of Sect. 2.5 is
regular, iff\\
1.) $\M$ is a factor (since $\A_L\cap\A_R=\CC$),
\\
2.) $H=\Ker(\pr_{\Out\M}\circ\al)$ (outerness),
\\
3.) the implementers $u(h),\ h\in H$, are linearly
independent in $\M$ (standardness).
\\[.2cm]
\prop{4.17} then implies that
a $G$-action $\al$ on a factor $\M$ is Galois, if and
only if for
$H:=\Ker(\pr_{\Out\M}\circ\al)$
the above ``standardness--condition" 3.) holds%
\footnote
{The case where $u:H\to\M$ is just a projective representation
is treated in Appendix B.
}.
In fact, we have
$\M\cros\A\equiv\M\cros_\al G=\M\o\CC G$ as a
linear space and therefore the elements
$u(h)h^{-1}\in\M\cros_\al G,\ h\in H$, are always linearly independent.
Hence, by \no{2.5.8} and \no{0.7c}, if $\M$ is a factor the
conditions
$\M'\cap(\M\cros\A)=\one_\M\cros\A_R\cong\A_R$
of \prop{4.17}i) always holds in this example.


\begin{appendix}

\renewcommand{\thesection}{\,}

\renewcommand{\theequation}{\mbox{\Alph{subsection}.\arabic{equation}}}
\renewcommand{\thesubsection}{\Alph{subsection}}
\renewcommand{\thetheorem}{\Alph{subsection}\arabic{theorem}}
\ncm{\subsec}{\setc{0}\subsection}

\sec{Appendix}

\subsec{Galois actions}

In this Appendix we generalize the notion of a Galois
(co-)action introduced in the case of ordinary Hopf algebras by [CS].
We then show that an action $\lef :\A\o\M\to\M$ is Galois if and only
if $\M\cros\A=\M\,h\,\M$. Hence, under the setting of \thm{4.3} this
is further equivalent to $p=\one$ and therefore to
$\M^\A\subset \M\subset\M\cros\A$ being a Jones triple.
However we emphasize, that the methods and results of this Appendix
apply to arbitrary $\A$-module algebras $\M$.

\bigskip
We start with considering the linear space
$\M\o_{\A_L}\hat\A$, where $a\in\A_L$ acts on $\M$ from the right by right
multiplication with $a\lef \one_\M$ and on $\hat\A$ from the left by right
multiplication with $\hat S^{-1}(a\arr\hat\one)$.
Then $\M\o_{\A_L} \hat \A$ naturally becomes a cyclic left
$(\M\o\hat\A)$-module, which in fact is isomorphic to the module
$(\M\o\hat\A)\hat\rho(\one_M)$.

\bsn
\Lem{C1}
{The map $f:\M\o_{\A_L}\hat\A\to(\M\o\hat\A)\hat\rho(\one_\M)$
$$
f(m\o_{\A_L}\varphi):=(m\o\varphi)\hat\rho(\one_M)
$$
provides a well defined isomorphism of left $(\M\o\hat\A)$-modules.
}
\proof
To show that $f$ is well defined we have to check
\beq\lb{4.15}
(a\lef\one_M\o\hat\one)\hat\rho(\one_\M)=
\big(\one_\M \o\hat S^{-1} (\e_R(a)\big)
\hat\rho(\one_\M)
\eeq
for all $a\in \A_L$. Equivalently, this means
\beq\lb{4.16}
(a\lef\one_\M)(b\lef\one_\M) =\big(b\arl \hat
S^{-1}(\e_R(a))\big) \lef\one_\M
\eeq
for all $a\in \A_L$ and all $b\in\A$. Now, by \lem{2.16}i)
the l.h.s is equal
to $ab\lef\one_\M$ and applying \lem{2.16}i)
to the left action of $(\hat\A)_{op}$ on $\A$ we get
$$
b\arl \hat S^{-1}\big(\e_R(a)\big)= \big( \one\arl
\e_R(a)\big)b=ab
$$
where we have used $\hat\e_L\e_R|\A_L=id$ and the fact that the antipode on
$(\hat\A)_{op}$ is given by $\hat S^{-1}$.
This proves \no{4.16} and therefore $f$ is well defined.
Clearly, $f$ is a surjective left $(\M\o\hat\A)$-module map. To
show that $f$
is an isomorphism let $\bar f:\M\o\hat\A\to\M\o_{\hat\A_L}\hat\A$ be the
canonical projection, which is also a $(\M\o\hat\A)$-left module map. Hence,
putting $\bar\rho=\bar f\circ \hat \rho$ we get
$$
(\bar f\circ f)(m\o_{\A_L}\varphi)=(m\o\varphi)\bar\rho(\one_\M)
$$
and therefore we may conclude $\bar f\circ f=id$ provided
\beq\lb{4.17}
\bar\rho(\one_\M)=\one_\M\o_{\A_L} \hat\one
\eeq
To prove \no{4.17} we write
$\hat\rho(\one_\M)=\sum m_i\o\xi^i\in \M\o\hat\A_L$
by \lem{2.12} and use \no{2.26} and \no{2.25b} to compute
\beanon
\bar\rho(\one_\M) &=& \sum m_i\o_{\A_L} (\e_L\hat\e_L)(\xi^i)\\
&=& \sum m_i\o_{\A_L}(\hat S^{-1} \e_R\hat\e_L)(\xi^i)\\
&=& \sum m_i(\e_L(\xi^i)\lef\one_\M)\o_{\A_L} \hat\one\\
&=& (\one_{(1)} \lef\one_\M)(\one_{(2)}\lef\one_\M)\o_{\A_L}\hat\one\\
&=& \one_\M\o_{\A_L} \hat\one
\eeanon
Thus indeed $\bar f\circ f=id$ which proves that $f$ is an isomorphism.
\qed

\bsn
Inspired by [CS] we now introduce the map
$\gamma:\M\o_\N\M\to\M\o_{\A_L}\hat\A$
given by
\beq\lb{4.18}
\gamma(m\o_\N m'):=(m\o \hat\one)\bar\rho(m')
\eeq
where $\N=\M^\A$ as before. Similarly as in Definition 7.3 of [CS] we
then define

\Definition{C2}
{Let $\hat\rho:\M\to\M\o\hat\A$ be the coaction corresponding to $\lef :\A
\o\M\to\M$ and put $\N\equiv \M^\A$. Then $\N\subset\M$ is called an
{\em $\hat\A$-Galois extension} (equivalently, $\hat\rho$ and $\lef$ are
called {\em Galois (co-) actions}, respectively), if the map $\gamma$ in
\no{4.18} is bijective.
}
Let us check this definition for the example of the partly
inner group action in Sect. 2.5. In this example we have a
$*$-algebra embedding $i:\hat\G\to\hA$, where $\hat\G$ is the
dual of the group algebra $\CC G$ (i.e. the abelian function
algebra on $G$), given by
$$
[i(\varphi)](h,g):=|H|\delta(h)\varphi(g),\quad h\in H,\,g\in G.
$$
Let $\rho:\M\to\M\o\hA$ be the coaction dual to the $\A$-action
\no{2.5.13}.
Then $\rho(\one_\M)=|H|^{-1}\sum_{h\in H}u(h)\o\chi_h$, where
$\chi_h(h',g):=|H|\delta(h^{-1}h')$ is the basis of
$\CC H\subset\hA$. Using this, one straightforwardly verifies
that
$$
{\bf i}: \M\o\hat\G\ni(m\o\varphi)\mapsto
(m\o i(\varphi))\rho(\one_\M)\in(\M\o\hA)\rho(\one_\M)
$$
provides  a linear isomorphism.
Moreover, using \no{2.2.5a} we have for all $m\in\M$
\beanon
\rho(m) &=&
\sum_{(h,g)\in H\x G}u(h)\al_g(m)\o\delta_{(h,g)}
\\
&=&\sum_{(h,g)\in H\x G}\al_g(m)u(h)\o\delta_{(h,h^{-1}g)}
\\
&=&{1\over |H|}\sum_{(h,g)\in H\x G}
\al_g(m)u(h)\o i(\delta_g)*\chi_h
\\
&=&{\bf i}(\sigma_\al(m))
\eeanon
where
$\s_\al:\M\ni m\mapsto\sum_{g\in G}(\al_g(m)\o\delta_g)\in\M\o\hat\G$
is the coaction dual to the $G$-action $\al$.
Together with \lem{C1} this proves that the $\A$-action $\re$
associated with $(\al,u)$ is Galois in the sense of
\defi{C2}, if and only if the $G$-action $\al$ is
Galois in the usual sense.

Let us now show that -- as for ordinary Hopf algebras --
$\re:\A\o\M\to\M$ is Galois, if and
only if
$\N\equiv\M^\A\subset\M\subset\M\cros \A$ is a Jones triple.

\Proposition{C3}
{{\rm (Galois actions and Jones triples)\\}
Let $\hat\rho : \M\to\M\o\hat\A$ be a coaction with dual left
$\A$-action $\lef$ and
let $\lambda\in\L(\hat\A)$ be the left integral dual to the normalized
Haar integral $h\in\A$. Then
\\
i) The maps
$id_\M\o h_L:\M\o\hat\A\to\M\o\A$ and
$id_\M\o (\hat S^{-1} \circ\lambda_R):\M\o\A\to\M\o\hat\A$
are inverses of each other and pass to well
defined left $\M$-module isomorphisms $ F_\M:\M\o_{\A_L} \hat\A \to \M\cros\A$
and $ F_\M^{-1} :\M\cros\A\to \M\o_{\A_L}\hat\A$, respectively.
\\
ii) $( F_\M\circ \gamma)(m\o_\N m')=mhm'$
\\
iii) If the conditional expectation $E_h:\M\to\M^\A$ is of
index-finite type then the map $\gamma$ is injective.
\\
iv) $\gamma$ is bijective
(i.e. $\hat\rho$ is Galois) if and only if
$\M\cros \A=\M h\M$.
}
\proof
i) By \prop{2.23} $h_L^{-1}=\hat S^{-1} \circ \lambda_R$.
To show that $ F_\M$
is well defined we compute for $a\in\A_L$ and $m\o\varphi
\in \M\o\hat\A$
\beanon
 F_\M (m(a\lef\one_M)\o_{\A_L}\varphi) &=& m a(h\arl \varphi)\\
&=& m(h\arl(\varphi\arl S(a)))\\
&=& m(h\arl \varphi\e_L(S(a)))\\
&=&  F_\M(m\o_{\A_L} \varphi\hat S^{-1}(\e_R(a)))
\eeanon
Here we have used \no{1.18} in the second line, \lem{2.16}ii) for the left
action of $\A_{op}$ on $\hat\A_{op}$ in third line and the identity $\e_L
\circ S=\hat S^{-1} \circ \e_R$ in the last line.
Conversely, $ F_\M^{-1}$ is
also well defined, since for $a\in \A_L$ and $(m\o b) \in \M\o \A$
\beanon
 F_\M^{-1} (m(a\lef\one_M)b)
&=& m(a\lef\one_M)\o_{\A_L}\hat S^{-1} (b\arr \lambda)\\
&=& m\o_{\A_L} \hat S^{-1} ((a\arr\hat\one)(b\arr\lambda))\\
&=& m\o_{\A_L}\hat S^{-1}((ab))\arr\lambda)\\
&=&  F_\M^{-1} (mab)
\eeanon
where in the third line we have used \lem{2.16}i). This proves
part i).
Using the
notation $\hat\rho(m')=m'_{(0)} \o m'\1$ part ii) follows from
\beanon
( F_\M\circ \gamma)(m\o_\N m') &=& m\,m'_{(0)} (h\arl m\1')\\
&=& m(h\1\lef m')h\2\\
&=& mhm'
\eeanon
The injectivity of $\gamma$ in part iii) follows from the well
known fact that
the map $ F_\M\circ\gamma : \M\o_\N\M\ni (m\o_\N m')\mapsto me_1 m'\in
\M e_1\M\equiv \M_1$,
provides an isomorphism whenever $e_1$ is the Jones projection
associated with
a conditional expectation $E:\M\to \N$ of index-finite type
[Wa, Prop.1.3.3].
Specifically, using our identification $\M_1=\pi_\omega(\M\cros\A)$,
under the setting of Theorem 4.3 the inverse of this map is given by
$$
\M_1 \ni \pi_\omega (ma)\mapsto\sum_i m(a\lef u_i)\o_\N u_i^*\in \M
\o_\N \M
$$
Finally, to prove part iv) we note that
$ F_\M$ being bijective $\gamma$ is surjective if and only if
$\M\cros\A=\M h\M$. But under this condition there also exist
$u_i,v_i\in\M$ such that $\one_{\M\cros\A}=\sum u_ihv_i$.
Consequently, $(u_i,v_i)$ provide a quasi-basis for $E_h$ and
therefore, by part iii), in this case $\gamma$ is also injective.
\qed

\bsn
\prop{C3} justifies our use of the terminology ``Galois action"
in Section 4, since under the setting of \thm{4.3}
it just means $p=\one$ and therefore $\M\cros\A\cong\M_1$.



\subsec{Partly inner group actions}

In this Appendix we analyze the example of a partly inner group
action in more generality.
So, let
$\al :G\to\Aut \M$ be an action of a finite group $G$ on a factor
$\M$ and let
$H:=\Ker (\prr \circ \al)\subset G$, where $\prr:\Aut \M\to \Out \M$ is the
canonical projection onto the group of outer automorphisms of $\M$
(i.e. $\Out\M=\Aut\M/ \Ad U(\M)$, where $U(\M)$ denotes the group of
unitaries in $\M$). Then $H\subset G$ is a normal subgroup and there
exists a section
$
u:H\to U(\M)
$
such that
$
\Ad u(h) = \al_h,\ \forall h\in H.
$
Associated with the section $u$ there
exists a 2-cocycle $z:H\times H\to U(1)$
such that for all $h,h'\in H$
\beq\lb{B3}
u(h) u(h') = z(h,h')u(hh')
\eeq
Moreover, without loss of generality
we may assume $u(h^{-1}) = u(h)^{-1} \equiv u(h)^*$ and
$u(\one_H)=\one_\M$, implying
$$
z(h,h^{-1})=z(\one,h)=z(h,\one)=1,\quad\forall h\in H.
$$
We then define the {\em $z$-twisted group algebra}
 $\CC H_z$ as the abstract $C^*$-algebra generated by the relations
\no{B3}, i.e. $\CC H_z=\CC H$ as a linear space with
$*$-algebra structure given on
the basis $h\in H$ by
\bea
\lb{B4} h *_z h' &:=& z(h,h') hh'\\
\lb{B5} h^* &:=& h^{-1}
\eea
Then, $u:H\to U(\M)$ extends linearly to a unital $*$-algebra homomorphism
$u:\CC H_z\to \M$.
Since $\M$ is a factor and $\alpha_h = \Ad u(h)$ for all $h \in H$ we must
have
\beq\lb{B6}
\alpha_g(u(h)) = c(g,h) u(ghg^{-1}),~~\forall g\in G
\eeq
 for some function $c:G\times H\to U(1)$ satisfying
\beq\lb{B7} c(g,\one)= 1 = c(\one,h)
\eeq
\beq\lb{B8}
c(g_1 g_2,h) = c(g_1,g_2 h g_2^{-1} ) c(g_2,h)
\eeq
\beq\lb{B9}
z(h_1,h_2) c(g,h_1 h_2) = c(g,h_1) c(g,h_2) z(gh_1g^{-1}, g h_2 g^{-1})
\eeq
\beq\lb{B10}
c(h_1,h_2)=z(h_1,h_2) z(h_1 h_2,h_2^{-1}) \equiv z(h_2,h_1^{-1})z(h_1,h_2h_1^{-1})
\eeq
for all $g,g_1,g_2 \in \G$ and all $h,h_1,h_2 \in H$.
Any such function $c$ allows to define a {\em twisted adjoint
action} $\beta:G\to \Aut (\CC H_z)$  given on the basis $h\in H$ by
\beq\lb{B11}
\beta_g(h) := c(g,h) ghg^{-1}
\eeq
and extended linearly to $\CC H_z$.
Clearly, for $h,h' \in H$ we have $\beta_h(h')=
h*_z h' *_z h^{-1}$ by \no{B10}, implying also $\beta_h(h)=h$.

With these data the crossed product
\beq\lb{B12}
\A:= \CC H_z \>cros_\beta G
\eeq
carries a weak $C^*$-Hopf algebra structure given similarly as
in Sect. 2.5 by
\bea
\lb{B20}
\Delta (h,g) &:=&
{1\over |H|}\sum_{\tilde h\in H}
(h *_z\tilde h^{-1},\ \tilde h g) \o (\tilde h,g)
\\
\lb{B21}
\e (h,g) &:=& |H| \delta (h)
\\
\lb{B22}
S(h,g) &:=& (\beta_g^{-1} (h), g^{-1} h^{-1})
\equiv c(g^{-1},h)\,(g^{-1}hg,\,g^{-1}h^{-1})
\eea
As in the untwisted case of Sect. 2.5 we have
\bea\lb{B2.5.7}
\A_L&=& \span\{(h,\one_G)\mid h\in H\} \cong \CC H_z
\\
\lb{B2.5.8}
\A_R&=& \span\{(h,h^{-1})\mid h\in H\} \cong (\CC H_z)_{op}
\eea
and for $m\in\M$ the definition
$$
(h,g)\re m:=u(h)\al_g(m)
$$
provides an $\A$-module algebra structure on $\M$.
As in Sect. 4.2 this action is outer iff
$\Ker (\prr\circ\alpha)=H$ and in this case it is Galois if and
only if it is standard, i.e. iff the implementers $u(h),\
h\in H$, are linearly independent in $\M$, see the remarks at
the end of Sect. 4.4. Hence, in this case $\A$ acts regularly
on $\M$ and we have the same scenario as for the untwisted case.

Let us conclude with mentioning that also in the twisted case
we have $S^2=\id$ and for $x=(h,g)\in\A$ one easily verifies
the formulas
\bea
\lb{B23} \e(\one\1x)\one\2 &=& (h,\one) = x\1 S(x\2)\\
\lb{B24} \one\1\e (x\one\2)
&=& (\beta_g^{-1} (h), g^{-1} h^{-1} g) = S(x\1)x\2.
\eea
These imply, that the normalized Haar integrals $e_{Haar}\in\A$
and $\l_{Haar}\in\hA$ are given by the same formulas as in
\no{Haar} and \no{2.6.4}, i.e.
\beanon
e_{Haar} &=& \frac{1}{|G|}\sum_{g\in G}(\one_H,g)
\\
\l_{Haar} (h,g) &=& |H|\delta(h)\delta(g),
\quad h\in H,\,g\in G.
\eeanon

\end{appendix}

\end{document}